\documentclass[preprint,3p,12pt,sort&compress]{elsarticle}

\pdfoutput=1

\usepackage{graphicx}
\usepackage{amsmath}
\usepackage{amsfonts}
\usepackage{enumitem} 							
\usepackage[labelfont=bf,tableposition=top]{caption} 	
\usepackage{subfig}									
\usepackage{epstopdf}								
\usepackage{floatrow}
\newfloatcommand{capbtabbox}{table}[][\FBwidth]
\newsavebox{\bigleftbox}
\floatstyle{plaintop}
\restylefloat{table}
\usepackage[table,xcdraw]{xcolor}
\usepackage{soul}
\usepackage{mathrsfs}  
\usepackage{multirow}
\usepackage{amsmath}
\usepackage{amssymb}
\usepackage[capitalise]{cleveref}
\crefname{rmk}{Remark}{Remarks}
\usepackage{titlesec}
\usepackage{etoolbox}
\usepackage{lipsum}
\usepackage{chngpage}
\usepackage{wrapfig}
\usepackage{subfiles}

\usepackage{tikz}

\newtheorem{remark}{Remark}[section]

\usepackage{dsfont}
\usepackage{empheq}
\usepackage{booktabs}

\renewcommand{\vec}[1]{\boldsymbol{#1}}

\renewcommand{\u}{\vec{u}}

\newcommand{\Ndofs}{{N_{\text{dofs}}}}
\newcommand{\Nelem}{{N_{\text{elem}}}}

\newcommand{\Dom}{\Omega}

\newcommand{\bdy}[1]{\partial #1}
\newcommand{\bdyDom}{{\bdy{\Dom}}}

\makeatletter
\def\H{\@ifnextchar\bgroup\H@arg\H@noarg}
\newcommand{\H@arg}[1]{H^1( #1 )}
\newcommand{\H@noarg}{H^1}
\def\Ho{\@ifnextchar\bgroup\Ho@arg\Ho@noarg}
\newcommand{\Ho@arg}[1]{H^1_0( #1 )}
\newcommand{\Ho@noarg}{H^1_0}
\def\BH{\@ifnextchar\bgroup\BH@arg\BH@noarg}
\newcommand{\BH@arg}[1]{\vec{H}^1( #1 )}
\newcommand{\BH@noarg}{\vec{H}^1}
\def\BHo{\@ifnextchar\bgroup\BHo@arg\BHo@noarg}
\newcommand{\BHo@arg}[1]{\vec{H}^1_{\vec{0}}( #1 )}
\newcommand{\BHo@noarg}{\vec{H}^1_{\vec{0}}}
\def\Hdiv{\@ifnextchar\bgroup\Hdiv@arg\Hdiv@noarg}
\newcommand{\Hdiv@arg}[1]{\vec{H}( \text{div}; #1 )}
\newcommand{\Hdiv@noarg}{\vec{H}(\text{div})}
\def\Hdivo{\@ifnextchar\bgroup\Hdivo@arg\Hdivo@noarg}
\newcommand{\Hdivo@arg}[1]{\vec{H}_{\vec{0}}( \text{div}; #1 )}
\newcommand{\Hdivo@noarg}{\vec{H}_{\vec{0}}(\text{div})}
\def\L{\@ifnextchar\bgroup\L@arg\L@noarg}
\newcommand{\L@arg}[1]{L^2( #1 )}
\newcommand{\L@noarg}{L^2}
\def\Lp{\@ifnextchar\bgroup\Lp@arg\Lp@noarg}
\newcommand{\Lp@arg}[1]{L^p( #1 )}
\newcommand{\Lp@noarg}{L^p}
\def\Cz{\@ifnextchar\bgroup\Cz@arg\Cz@noarg}
\newcommand{\Cz@arg}[1]{\mathcal{C}^0( #1 )}
\newcommand{\Cz@noarg}{\mathcal{C}^0}
\def\Cinf{\@ifnextchar\bgroup\Cinf@arg\Cinf@noarg}
\newcommand{\Cinf@arg}[1]{\mathcal{C}^\infty( #1 )}
\newcommand{\Cinf@noarg}{\mathcal{C}^\infty}
\def\Cinfo{\@ifnextchar\bgroup\Cinfo@arg\Cinfo@noarg}
\newcommand{\Cinfo@arg}[1]{\mathcal{C}^\infty_0( #1 )}
\newcommand{\Cinfo@noarg}{\mathcal{C}^\infty_0}
\def\Cneg{\@ifnextchar\bgroup\Cneg@arg\Cneg@noarg}
\newcommand{\Cneg@arg}[1]{\mathcal{C}^{-1}( #1 )}
\newcommand{\Cneg@noarg}{\mathcal{C}^{-1}}
\def\CI{\@ifnextchar\bgroup\CI@arg\CI@noarg}
\newcommand{\CI@arg}[1]{\mathcal{C}^{1}( #1 )}
\newcommand{\CI@noarg}{\mathcal{C}^{1}}
\def\CII{\@ifnextchar\bgroup\CII@arg\CII@noarg}
\newcommand{\CII@arg}[1]{\mathcal{C}^{2}( #1 )}
\newcommand{\CII@noarg}{\mathcal{C}^{2}}
\makeatother

\makeatletter
\def\Proj{\@ifnextchar\bgroup\Proj@arg\Proj@noarg}
\newcommand{\Proj@arg}[1]{\mathscr{P}_{#1}}
\newcommand{\Proj@noarg}{\mathscr{P}}
\makeatother

\newcommand{\jump}[1]{\mbox{$[\![ #1 ]\!]$}}

\renewcommand{\d}[1]{\,\mathrm{d}#1}
\newcommand{\dDom}{\d{V}}

\newcommand{\dbdy}{\d{S}}

\newcommand{\SF}{\textsc{sf}}
\newcommand{\DA}{\textsc{12}}
\newcommand{\SD}{\textsc{s1}}
\newcommand{\SA}{\textsc{s2}}
\newcommand{\DD}{\textsc{1}}
\renewcommand{\AA}{\textsc{2}}
\newcommand{\EQ}[1]{(\ref{eq:#1})}

\newcommand{\aGN}{\alpha_{\text{GN}}}
\newcommand{\aDW}{\alpha_{\text{DW}}}
\graphicspath{{figures/}}

\usepackage{subfiles} 

\title{Stabilized immersed isogeometric analysis for the Navier-Stokes-Cahn-Hilliard equations, with applications to binary-fluid flow through porous media}

\author[add1]{Stein K.F. Stoter}
\ead{k.f.s.stoter@tue.nl}
\author[add1]{Tom B. van Sluijs}
\ead{t.b.v.sluijs@tue.nl}
\author[add1]{Tristan H.B. Demont}
\ead{t.h.b.demont@tue.nl}
\author[add1]{E. Harald van Brummelen}
\ead{e.h.v.brummelen@tue.nl}
\author[add1]{Clemens V. Verhoosel\corref{cor}}
\ead{c.v.verhoosel@tue.nl}

\cortext[cor]{Corresponding author: C.V.\,Verhoosel ({c.v.verhoosel@tue.nl})}
\address[add1]{Eindhoven University of Technology, PO Box 513, 5600 MB Eindhoven, The Netherlands}

\date{May 2023}

\begin{document}

\begin{abstract}
    Binary-fluid flows can be modeled using the Navier-Stokes-Cahn-Hilliard equations, which represent the boundary between the fluid constituents by a diffuse interface. The diffuse-interface model allows for complex geometries and topological changes of the binary-fluid interface. In this work, we propose an immersed isogeometric analysis framework to solve the Navier-Stokes-Cahn-Hilliard equations on domains with geometrically complex external binary-fluid boundaries. The use of optimal-regularity B-splines results in a computationally efficient higher-order method. The key features of the proposed framework are a generalized Navier-slip boundary condition for the tangential velocity components, Nitsche's method for the convective impermeability boundary condition, and skeleton- and ghost-penalties to guarantee stability. A binary-fluid Taylor-Couette flow is considered for benchmarking. Porous medium simulations demonstrate the ability of the immersed isogeometric analysis framework to model complex binary-fluid flow phenomena such as break-up and coalescence in complex geometries.
\end{abstract}


\maketitle

\section{Introduction}

Binary-fluid flows, in which the two fluid components are separated by a molecular transition layer, are omnipresent in science and engineering. This article primarily focuses on imbibition processes in porous media -- which occur in such diverse applications as inkjet printing, groundwater flows and reservoir engineering -- but the same flow physics is relevant in numerous other contexts, including free-surface flows (\emph{e.g.}, waves, jets) and bubbly flows (\emph{e.g.}, reactor cooling, steam generation). The complex physical behavior of such flows makes the use of computational methods for their study indispensable \cite{prosperetti2009computational}.

The complexity of the domain and the loading conditions makes high-fidelity simulations of binary-fluid flows often extremely challenging, even when considering state-of-the-art computational techniques. The challenges related to such simulations are of both physical and numerical nature, \emph{viz.}: \emph{(i)} The interface layer separating the two fluids is subject to break-up and coalescence, changing not only the shape, but also the topology of the domain of each constituent; \emph{(ii)} The dynamics of the contact lines, \emph{i.e.}, the intersection of the interface layers with the boundary of the solid domain, is of essential importance and must be emulated accurately by the computational model; \emph{(iii)} The geometry of the binary-fluid flow domain can be very complex (\emph{e.g.}, scan-based micro-structures), rendering automatic high-quality mesh generation for the construction of (higher-order) stable approximation spaces nearly impossible.

To capture topological changes of the interface layer due to, \emph{e.g.}, break-up and coalescence, we employ a diffuse (or smeared) interface model. Diffuse-interface models describe the fluid components by a continuous phase field, which varies gradually over a thin-but-finite transition zone representing the interface layer. In contrast to sharp (or discrete) interface models, diffuse models have the intrinsic ability to accommodate topological changes without the need for explicit interface reconstruction \cite{Sun2007}. Diffuse-interface models for two immiscible incompressible fluid species are generally described by the Navier-Stokes-Cahn-Hilliard equations \cite{Hohenberg:1977hh,Lowengrub:1998uq,Simsek:2018gb,Abels:2012vn}. In this work we build upon the model by Abels, Garcke and Gr\"{u}n \cite{Abels:2012vn}, on account of its thermodynamic consistency, incompressibility, and consistent reduction to the underlying single-fluid Navier--Stokes equations in pure species.

To properly model the dynamics of contact lines, we make use of the generalized Navier boundary condition~\cite{Navier1823,Gerbeau:2009jx} for the velocity components tangent to the contact surface. This boundary condition circumvents the traction singularities~\cite{Cox:1986jq} which are associated with the no-slip condition in sharp interface models~\cite{huh1971hydrodynamic}. The Navier-Stokes-Cahn-Hilliard model implicitly introduces a slip length \cite{Jacqmin:2000kx} related to the diffuse interface thickness parameter (via the mobility), thereby in principle avoiding stress singularities. This intrinsic slip length is hard to control, however, and results in increasing stress concentrations and, ultimately, unbounded shear forces as the diffuse-interface thickness vanishes. A Navier-Stokes-Cahn-Hilliard model augmented with the generalized Navier boundary condition conjecturally permits a physically relevant sharp interface limit, although the nature and behavior of this limit are still open research questions~\cite{Abels:2018ly,Yue:2010hq,Brummelen:2021aw}.

To construct a suitable approximation space on domains that are complex both in terms of geometry and topology, we base our work on immersed isogeometric analysis. This analysis technique combines the favorable approximation properties of splines inherent to isogeometric analysis \cite{Hughes:2005it}, with the topological flexibility of immersed methods (\emph{e.g.}, the Finite Cell Method \cite{parvizian_finite_2007,duster_finite_2008,schillinger_finite_2015} and CutFEM \cite{burman_ghost_2010,burman_fictitious_2012,burman_cutfem_2015}). The combination of isogeometric analysis with immersed methods has been explored for a wide range of applications \cite{rank_geometric_2012,schillinger_isogeometric_2012,ruess_weakly_2013,kamensky_immersogeometric_2015,hsu_direct_2016,verhoosel_image-based_2015,ruess_finite_2012,de_prenter_multigrid_2020}. Compared to boundary-fitted isogeometric analysis (and finite element techniques in general), the immersed analysis approach requires dedicated integration techniques for cut cells \cite{abedian_performance_2013,divi_error-estimate-based_2020}, special treatment of essential boundary conditions \cite{burman_cutfem_2015}, which is particularly non-trivial in the context of non-linear and multi-field equations \cite{Douglas1975,Burman2015}, and stabilization techniques to circumvent problems associated with small volume-fraction cut cells \cite{deprenter2023stability,Stoter2023a}. In this work we adopt the integration technique developed in Refs.~\cite{verhoosel_image-based_2015,divi_topology-preserving_2022}, making it possible to numerically evaluate integrals over both the immersed domain and its boundaries. In the context of multi-phase problems, immersed techniques have been considered previously in the restricted scope of the Cahn-Hilliard equations separately~\cite{karatzas2021reduced}. The Navier-Stokes-Cahn-Hilliard system considered in this work requires dedicated treatment of boundary-condition imposition and stabilization. For this, we draw inspiration from the skeleton-stabilized immersogeometric analysis framework for the (Navier-)Stokes equations developed in Refs.~\cite{hoang_skeleton-stabilized_2018,hoang_skeleton-stabilized_2019}. 

In this contribution we propose a novel immersed isogeometric analysis formulation for the simulation of binary-fluid flows governed by the Navier-Stokes-Cahn-Hilliard equations. To impose boundary conditions on the immersed boundary, we derive a Nitsche formulation to control the normal component of the velocity. For the tangential components we propose the generalized Navier boundary condition, to regularize the traction singularity associated with interface pinning in the sharp-interface limit. We propose to discretize the formulation using equal-order optimal-regularity splines for all field variables. To stabilize this formulation, the skeleton- and ghost-stabilization techniques are extended to the problem considered here, focusing on minimizing the required number of penalty parameters. To focus on the novelties in terms of stabilized formulation and immersed isogeometric finite element discretization, our presentation is restricted to two-dimensional scenarios. The presented novelties are, to a large extent, not specific to this setting, but 
the adaptivity and high-performance computing aspects that are necessary to enable three-dimensional simulations are beyond the scope of this work.

We dedicate this contribution to Thomas J.R.\ Hughes for his lifetime achievements in computational mechanics. This contribution is strongly influenced by the seminal work conducted by Tom on all key aspects, \emph{viz}. isogeometric analysis, multi-phase flows and phase-field modeling, and stabilization techniques. We also draw inspiration from the elegance with which Tom develops and combines sophisticated computational techniques to solve complex, multifaceted, problems, acknowledging the importance of the usability of advanced methods in engineering workflows.

The remainder of this article is structured as follows. In \cref{sec:diffInt} we establish the diffuse-interface model equations appropriate for describing moving contact lines in binary-fluid flows. In \cref{sec:immersed} we introduce the necessary technology components to perform immersed finite element simulations, and in \cref{sec:WeakFormulation} we develop an immersed isogeometric analysis formulation for the diffuse-interface model. We demonstrate our immersed isogeometric analysis framework in \cref{sec:numExp} for a number of numerical experiments, including a binary-fluid Taylor-Couette flow benchmark case, and the binary-fluid flow through porous media. Conclusions are finally presented in \cref{sec:conclusion}.

\section{The diffuse-interface model}\label{sec:diffInt}

\begin{figure}[!t]
\begin{center}
\includegraphics[width=0.85\textwidth]{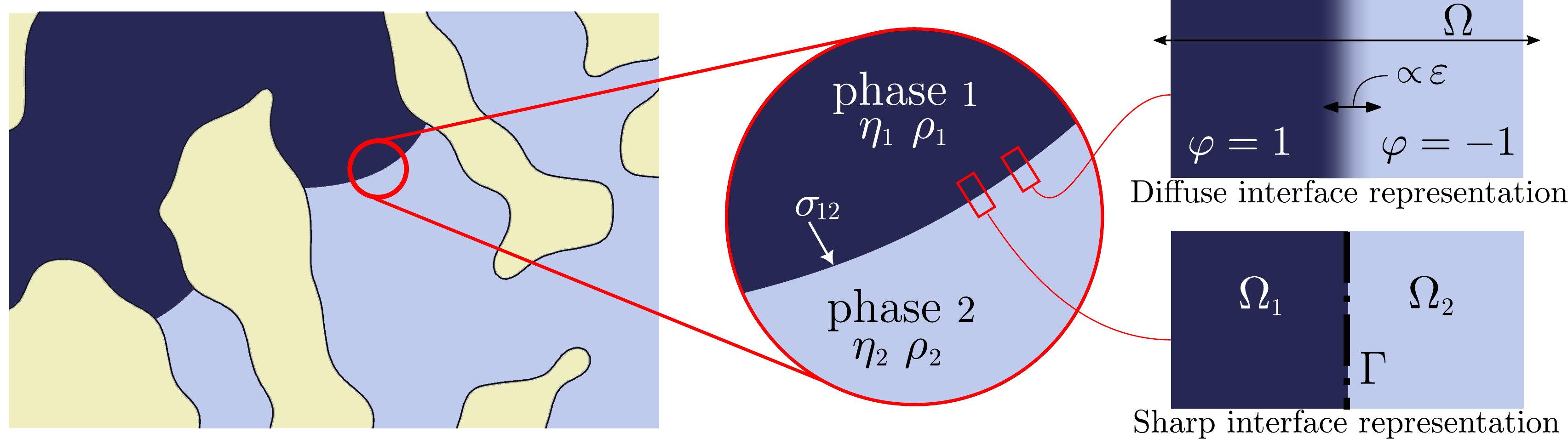}
\end{center}
\caption{Schematic of a porous medium binary-fluid system, modeled with either a diffuse interface representation or a sharp-interface representation.}
\label{fig:droplet_model}
\end{figure}

We examine porous media imbibition of an isothermal binary-fluid system consisting of two immiscible incompressible Newtonian fluids, which we label as $\DD$ and $\AA$. The two prevalent approaches for modeling binary-fluid systems are by a sharp-interface representation and by a diffuse-interface representation, both illustrated in \cref{fig:droplet_model}. On account of the prominent occurrence of interface coalescence and break-up in porous media imbibition processes, we consider a Navier-Stokes-Cahn-Hilliard diffuse-interface model.

\subsection{Navier-Stokes-Cahn-Hilliard equations}
\label{ssec:NSCH}
\noindent
In diffuse-interface models, the interface between the two immiscible fluids is characterized by a layer of finite thickness consisting of a mixture of both fluids. This layer represents a gradual transition between fluid $\DD$ and fluid $\AA$. The fluid-constituents in a diffuse-interface model are typically described by an order, or phase-field, parameter $\varphi \in [-1,1]$, where $\varphi=1$ signifies pure species $\DD$ and $\varphi=-1$ pure species $\AA$. We utilize the Navier-Stokes-Cahn-Hilliard equations to describe the evolution of the mixture in terms of the volume-averaged velocity, $\u$, the pressure, $p$, and the newly introduced order parameter, $\varphi$. Specifically, we adopt the model developed by Abels, Garcke, and Gr\"{u}n \cite{Abels:2012vn}, on the basis of its thermodynamic consistency and its capacity to reduce to the single-fluid Navier-Stokes equations in scenarios involving pure species. The balance laws for this model are given by
\begin{subequations} \label{eq:strong}
\begin{empheq}[right={\empheqrbrace \textrm{ in } \Omega}]{align}
\partial_t \left( \rho \u \right) + \nabla \cdot \left( \rho \u \otimes \u \right) + \nabla \cdot \left( \u \otimes \vec{J} \right) - \nabla \cdot \vec{\tau} - \nabla \cdot \vec{\zeta}  + \nabla p & = 0,\\
\nabla \cdot \u & = 0,\\
\partial_t \varphi + \nabla \cdot \left( \varphi \u \right) - \nabla \cdot \left( m \nabla \mu \right) & = 0.
\end{empheq}
\end{subequations}
The closure relations for the relative mass flux $\vec{J}$, the viscous stress $\vec{\tau}$, the capillary stress $\vec{\zeta}$ and the chemical potential $\mu$ are given as
\begin{subequations}\label{eq:closureeqs}
\begin{alignat}{3}
& \vec{J} \coloneqq - \frac{ \rho_{\DD} - \rho_{\AA} } { 2 } m \nabla \mu  \,, \label{eq:closureeqsJ}\\
& \vec{\tau} \coloneqq \eta \left( \nabla \u + ( \nabla \u )^T \right) \,, \\  
& \vec{\zeta} \coloneqq - \sigma \varepsilon \nabla \varphi \otimes \nabla \varphi + \vec{I} \left( \frac{ \sigma \varepsilon } { 2 } | \nabla \varphi |^2 + \frac{ \sigma } { \varepsilon } \Psi \right) \,, \label{eq:closureeqszeta}\\
& \mu = -\sigma \varepsilon \Delta \varphi + \frac{ \sigma } { \varepsilon } \Psi' \,, \label{eq:chempot}
\end{alignat}
\end{subequations}
in which $\Psi=\Psi(\varphi)$ is the mixture energy density, for which we make use of the typical double well function
\begin{align}
\Psi ( \varphi ) \coloneqq \frac{1}{4} \left( \varphi^2 -1 \right)^2 \,.    
\end{align}
The model parameters in \cref{eq:strong,eq:closureeqs} involve the interface thickness parameter $\varepsilon>0$ and the mobility parameter $m>0$. The former is the diffuse interface length-scale, and the latter controls the diffusive time scale as well as the model-intrinsic wall-slip length 
scale~\cite{Demont:2022dk,Jacqmin:2000kx}.

The material parameters involved in these model equations are the fluid-fluid surface tension $\sigma_\DA=\frac{2\sqrt{2}}{3}\sigma$, the mixture density $\rho$, and the mixture viscosity $\eta$. The mixture density and viscosity generally depend on~$\varphi$. To ensure positive densities even for the non-physical scenario $\varphi \notin [-1,1]$, we adopt the density extension~\cite{Bonart:2019re}
\begin{equation}\label{eq:densityextension}
\rho(\varphi) = \left\{
\begin{tabular}{ll}
$\frac{ 1 } { 4 } \rho_{\AA},$ & $\varphi \leq - 1 - 2 \lambda \,,$\\
$\frac{ 1 } { 4 } \rho_{\AA} + \frac{ 1 } { 4 } \rho_{\AA} \lambda^{-2} \left( 1 + 2 \lambda + \varphi \right)^2,$ & $\varphi \in ( - 1 - 2 \lambda , - 1 - \lambda) \,,$\\
$\frac{ 1 + \varphi } { 2 } \rho_{\DD} + \frac{ 1 - \varphi } { 2 } \rho_{\AA},$ & $\varphi \in [ - 1 - \lambda, 1 + \lambda ] \,,$\\
$\rho_{\DD} + \frac{ 3 } { 4 } \rho_{\AA} - \frac{ 1 } { 4 } \rho_{\AA} \lambda^{-2} \left( 1 + 2 \lambda - \varphi \right)^2,$ & $\varphi \in ( 1 + \lambda, 1 + 2 \lambda ) \,,$\\
$\rho_{\DD} + \frac{ 3 } { 4 } \rho_{\AA},$ & $\varphi \geq 1 + 2 \lambda \,,$
\end{tabular}
\right.
\end{equation}
with $\lambda = \rho_{\AA} / \left( \rho_{\DD} - \rho_{\AA} \right)$. Such a density extension is particularly important when considering constituents of widely varying densities, such as water and air.

For the viscosity interpolation, we make use of the Arrhenius mixture-viscosity model~\cite{Arrhenius:1887xr}
\begin{equation}\label{eq:arrhenius}
\log \eta( \varphi ) = \frac{ \left( 1 + \varphi \right) \Lambda \log \eta_{\DD} + \left( 1 - \varphi \right) \log \eta_{\AA} } { \left( 1 + \varphi \right) \Lambda + \left( 1 - \varphi \right) } \,,
\end{equation}
where $\Lambda = \frac{ \rho_{\DD} M_{\AA} } { \rho_{\AA} M_{\DD} }$ is the intrinsic volume ratio (with $M_{\DD}$ and $M_{\AA}$ the molar masses). For all future results, we use $\Lambda=1$, for which the denominator in~\cref{eq:arrhenius} is guaranteed to remain positive irrespective of $\varphi$ (see \emph{Remark 2} in Ref.~\cite{Brummelen:2021aw} for further details).

\subsection{Fluid-solid boundary condition}
\label{ssec:BCs}
The behavior of a binary fluid in contact with a solid surface and, in particular, of the contact line corresponding to the intersection of the fluid-fluid interface with the solid surface, has been the subject of extensive investigation over several decades, and contemporary understanding is still incomplete. For sharp-interface models, the standard no-slip condition at fluid-solid interfaces leads to a non-integrable stress singularity at the contact line~\cite{huh1971hydrodynamic}. This stress singularity is removed by the introduction of slip~\cite{Cox:1986jq}. Diffuse-interface models provide an intrinsic slip mechanism and, hence, regularization of the traction singularity~\cite{Jacqmin:2000kx}. However, in combination with a no-slip condition, the traction at the contact line still degenerates in the sharp-interface limit $\varepsilon\to{}+0$.

To avoid the degeneration of the fluid traction at the contact line in the sharp-interface limit, in the present work we employ a generalized Navier boundary condition~\cite{Gerbeau:2009jx}. The generalized Navier boundary condition is an extension of the classical Navier slip condition~\cite{Navier1823}, including capillary effects. In the diffuse-interface setting, the generalized Navier boundary condition is given by
\begin{equation}
\label{eq:GNBC}
\mathbf{P}_{\Gamma}
\big(
 \aGN (\u-\u_\text{wall})+( \vec{\tau}\vec{n} + \vec{\zeta}\vec{n} )\big)-\nabla_{\Gamma}\sigma_{\SF}(\varphi)=0\,,
\end{equation}
with the generalized Navier model parameter $\aGN >0$ as the relaxation coefficient, $\mathbf{P}_{\Gamma}(\cdot)=\vec{n}\times(\cdot)\times\vec{n}$ the tangential projection onto the solid surface, $\nabla_{\Gamma}(\cdot)=\mathbf{P}_{\Gamma}\nabla(\cdot)$ the surface gradient, and $\sigma_{\SF}$ the solid-fluid surface tension according to
\begin{equation}
\label{eq:sigmasf}
\sigma_{\SF}(\varphi)=\frac{1}{4}(\varphi^3-3\varphi)(\sigma_{\SA}-\sigma_{\SD})+\frac{1}{2}(\sigma_{\SD}+\sigma_{\SA})\,,
\end{equation}
where $\sigma_{\SD}\geq{}0$ and $\sigma_{\SA}\geq{}0$ 
denote the solid-fluid surface tensions of fluid species~1 and~2, respectively; 
see also 
Refs.~\cite{Jacqmin:2000kx,Yue:2011uq,Shokrpour-Roudbari:2016dp,Brummelen:2016qa}. Essentially, the generalized Navier boundary condition insists that the tangential slip velocity, $\mathbf{P}_{\Gamma}(\u-\u_{\text{wall}})$, is proportional to the tangential component of the fluid traction including the capillary stress, according to $\mathbf{P}_{\Gamma}( \vec{\tau}\vec{n} + \vec{\zeta}\vec{n} )$, and the Marangoni traction emerging from the non-uniform fluid-solid surface tension, $\nabla_{\Gamma}\sigma_{\SF}(\varphi)$.

The description of the wetting behavior of the two fluid components at the solid surface provided by~\EQ{GNBC} is complemented by the dynamic contact angle condition 
\begin{equation}
\label{eq:DCA}
\aDW (\partial_t\varphi+\u\cdot\nabla_{\Gamma}\varphi)
+
\sigma\varepsilon\,\nabla \varphi \cdot \vec{n}
+
\sigma'_{\SF}(\varphi)
=0\,,
\end{equation}
with the dynamic wetting model parameter $\aDW\geq{}0$ as a relaxation coefficient. For $\aDW=0$, the boundary condition~\EQ{DCA} is referred to as the static contact-angle condition, and one can infer that~\EQ{DCA} imposes the equilibrium contact angle~$\theta_{\textsc{e}}$ between the diffuse interface and the solid surface (interior to the liquid) in accordance with $\sigma_{\DA}\cos(\theta_{\textsc{e}})+\sigma_{\SD}-\sigma_{\SA}=0$. In the present work we restrict ourselves to stationary neutral wetting scenarios, \emph{i.e.}, $\aDW=0$ and~$\theta_{\textsc{e}}=\pi/2$. The latter implies that $\sigma_{\SD}=\sigma_{\SA}=\sigma_{\SF}$. The contact angle condition \eqref{eq:DCA} then reduces to the homogeneous Neumann condition
\begin{equation}
\label{eq:dnphi=0}
\nabla \varphi \cdot \boldsymbol{n}=0\,.
\end{equation}
We introduce this provision to avoid the trace terms in~\EQ{DCA}, which would severely complicate the weak formulation, requiring its own dedicated study. It is to be noted that the neutral wetting assumption also implies that the Marangoni term~$\nabla_{\Gamma}\sigma_{\SF}(\varphi)$ in~\EQ{GNBC} vanishes.

In conjunction with the generalized Navier boundary condition and the contact angle condition, we consider the impermeability conditions
\begin{align}
\vec{u}\cdot\vec{n}&=0\,,
\label{eq:noconv}
\\
-m\nabla\mu\cdot\vec{n}&=0\,.
\label{eq:nodiff}
\end{align}
These conditions respectively impose that the convective and diffuse transport into the solid surface vanish.

\section{Immersed isogeometric analysis}\label{sec:immersed}

To evaluate the Navier-Stokes-Cahn-Hilliard model on complex domains, such as our prototypical porous medium domain shown in Fig.~\ref{fig:fcmdomain}, we propose an immersed isogeometric analysis approach. In this section, we summarize the key technology components required to perform such an analysis. These are the construction of analysis-suitable spline spaces over non-boundary-fitted domains, discussed in Section~\ref{sec:immersedsplines}, and the algorithms for evaluating volumetric and surface integrals on cut elements, discussed in \cref{sec:integration}. 

\subsection{Non-boundary-fitted B-spline basis functions}\label{sec:immersedsplines}

\begin{figure}[!b]
    \centering
    \subfloat[Ambient and physical domain]{\includegraphics[width=0.49\linewidth]{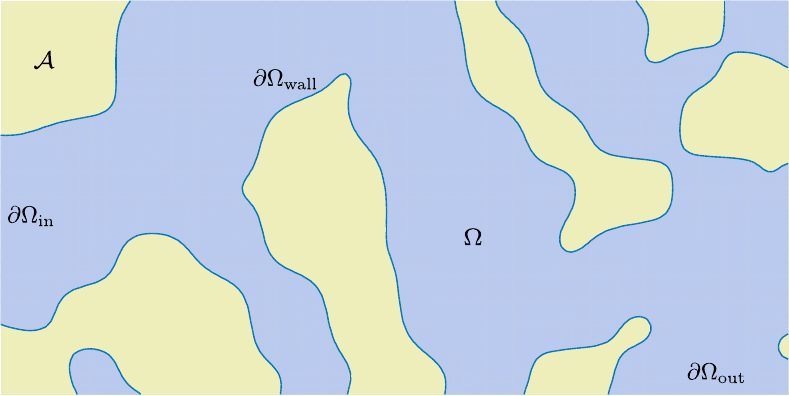}\label{fig:fcmdomaindef}}\hfill
    \subfloat[Ambient and background mesh]{\includegraphics[width=0.49\linewidth]{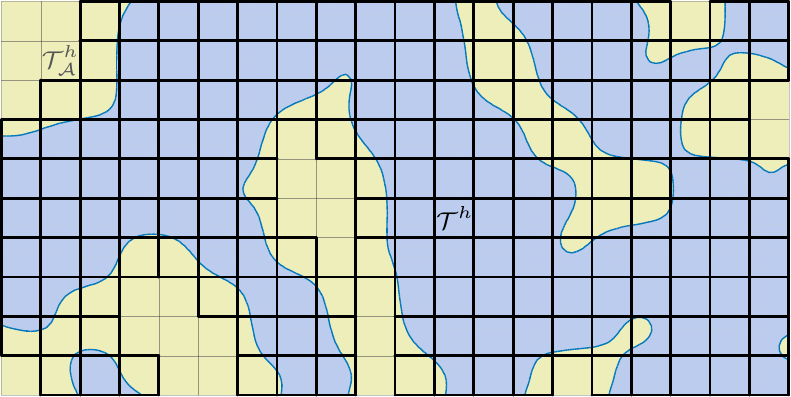}\label{fig:fcmmeshdef}}\\
    \subfloat[Skeleton faces]{\includegraphics[width=0.49\linewidth]{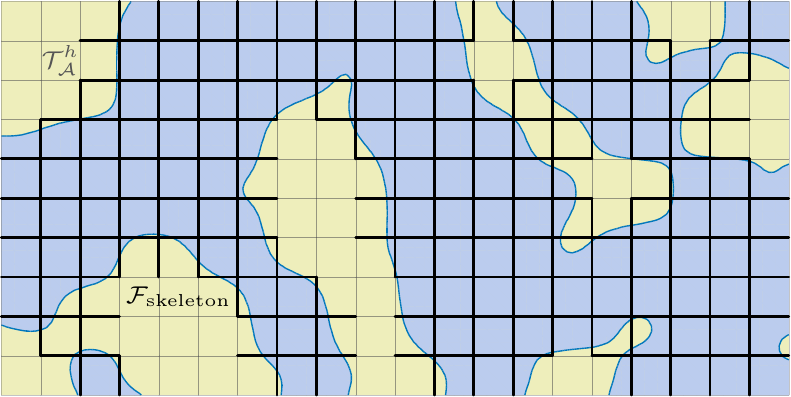}\label{fig:fcmskeltonfacets}}\hfill
    \subfloat[Ghost faces]{\includegraphics[width=0.49\linewidth]{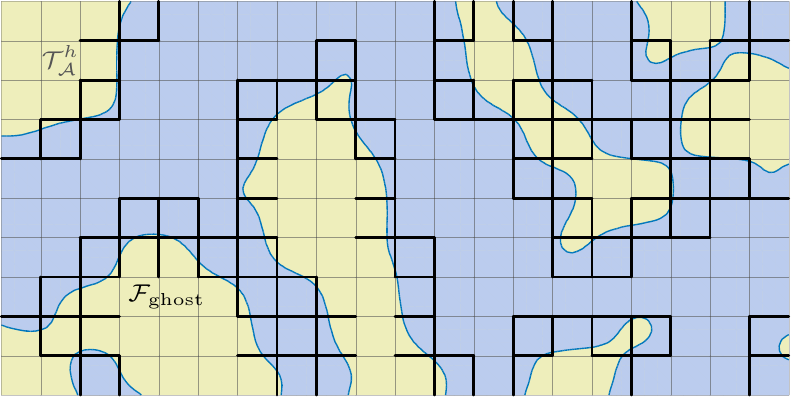}\label{fig:fcmghostfacets}}
    \caption{Illustration of the different domains and meshes used in the immersed setting, considering an artificially constructed porous medium as a prototypical example of a complex geometry.}
	\label{fig:fcmdomain}
\end{figure}

The physical domain, $\Omega$, with non-overlapping boundary segments $\partial \Omega_{\rm in} \cup \, \partial \Omega_{\rm out} \cup \, \partial \Omega_{\rm wall} = \partial \Omega$ is immersed in an ambient domain, $\mathcal{A} \supset \Omega$, which is typically cuboid; see Fig.~\ref{fig:fcmdomaindef}. The ambient domain is partitioned by a rectilinear mesh, $\mathcal{T}^h_{\mathcal{A}}$, with elements $K$, where the superscript $h$ refers the size of the elements. Elements that do not intersect with the physical domain can be omitted, resulting in the (active) background mesh
\begin{equation}\label{eq:background}
\mathcal{T}^{h} := \{K \, | \, K \in \mathcal{T}^h_{\mathcal{A}}, K \cap \Omega \neq \emptyset \}.
\end{equation}
The ambient domain mesh and (active) background mesh are illustrated in Fig.~\ref{fig:fcmmeshdef}. The mesh consisting of elements that are trimmed to the physical domain is denoted by
\begin{equation}\label{eq:cutmesh}
\mathcal{T}^h_\Omega := \{K \cap \Omega \, | \, K \in \mathcal{T}^h \}
\end{equation}
and the corresponding mesh for the immersed boundary by
\begin{equation}\label{eq:boundaryedges}
\mathcal{T}^h_{\partial \Omega} := \{ E \subset \partial \Omega \, | \, E = \partial K \cap \partial \Omega,\,  K \in \mathcal{T}^h_\Omega \}.
\end{equation}

The considered formulation, which we will detail in Section~\ref{sec:WeakFormulation}, incorporates stabilization terms formulated on the edges of the background mesh, which we refer to as the skeleton mesh
\begin{equation}\label{eq:skeletoninterfaces}
\mathcal{F}_{\rm skeleton}^h  := \{ F = \partial K \cap \partial K' \,|\ K,K'\in \mathcal{T}^h,  K \neq K' \}.
\end{equation}
Note that the faces, $F \in \mathcal{F}_{\rm skeleton}^h$, which are illustrated in Fig.~\ref{fig:fcmskeltonfacets}, can be partially outside of the domain, $\Omega$, and that the boundary of the background mesh is not part of the skeleton mesh.

We also define the ghost mesh, illustrated in Fig.~\ref{fig:fcmghostfacets}, as the subset of the skeleton mesh composed of faces that belong to an element intersected by the domain boundary, \emph{i.e.},
\begin{equation}\label{equation:ghostinterfaces}
\mathcal{F}_{\rm ghost}^h := \{ F \cap \partial K \,|\, F \in \mathcal{F}^h_{\rm skeleton}, K \in \mathcal{G} \},
\end{equation}
where $\mathcal{G} := \{ K \in \mathcal{T}^h \mid K \cap \partial \Omega \neq \emptyset \}$ is the collection of elements in the background mesh that are crossed by the immersed boundary.

\begin{figure}[!b]
    \centering
    \subfloat[$k=1$]{\includegraphics[width=0.33\linewidth,trim={0.25cm 0 0.75cm 0},clip]{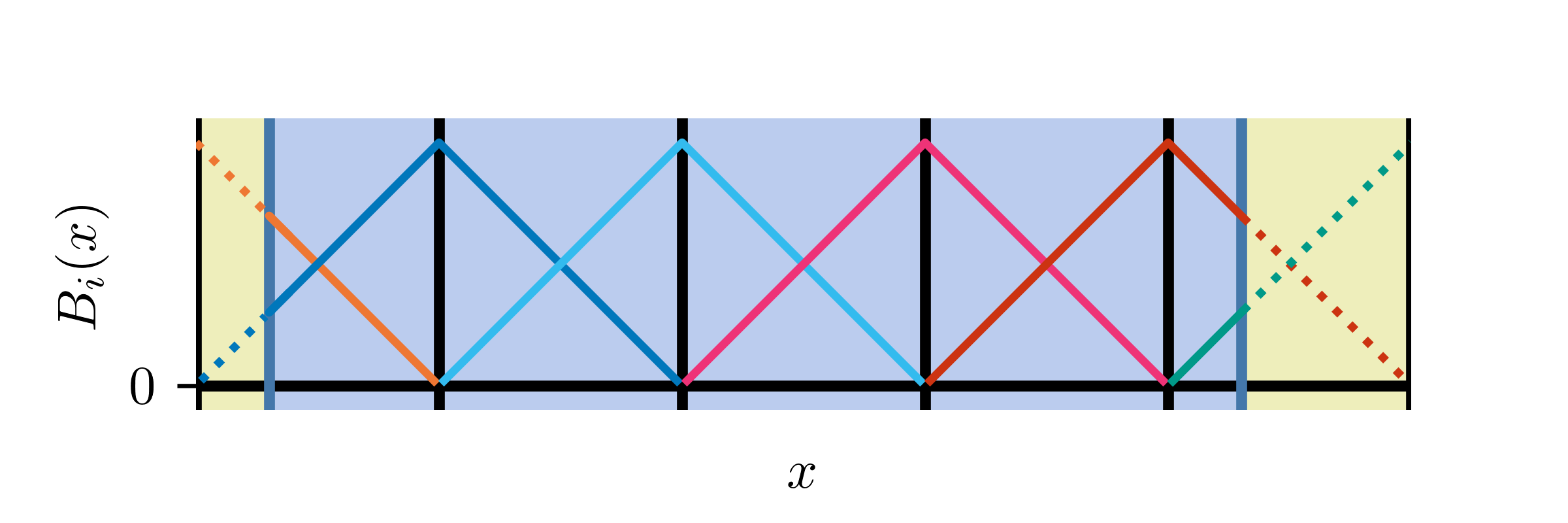}} \hfill
    \subfloat[$k=2$]{\includegraphics[width=0.33\linewidth,trim={0.25cm 0 0.75cm 0},clip]{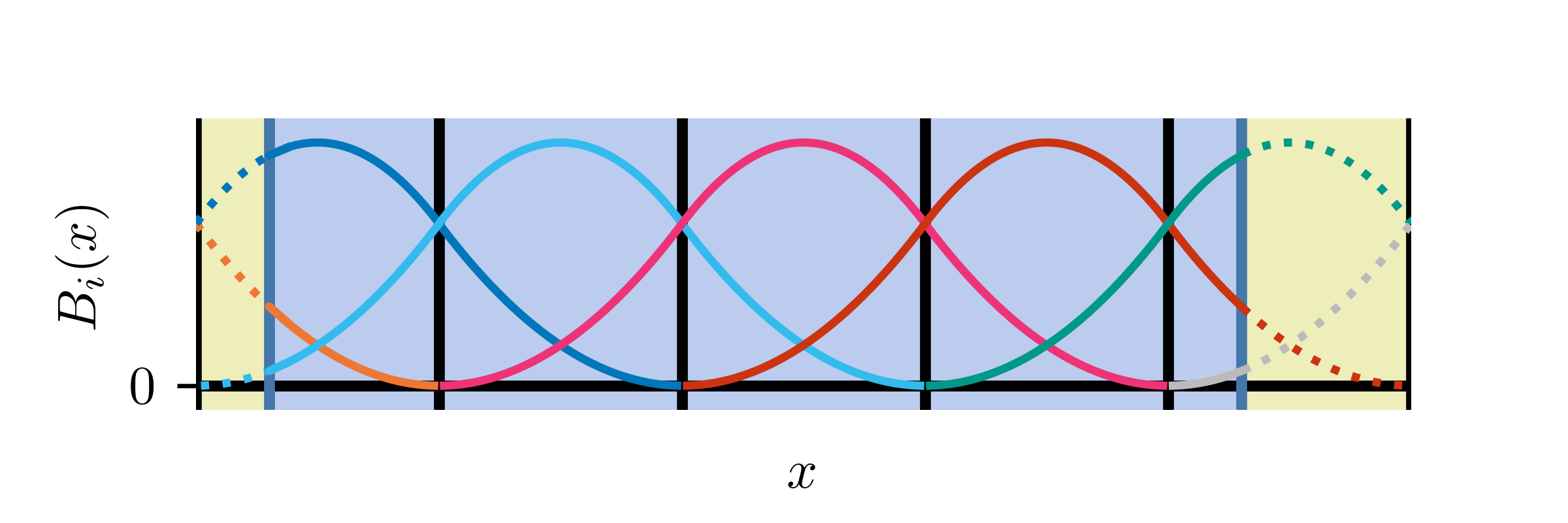}} \hfill
    \subfloat[$k=3$]{\includegraphics[width=0.33\linewidth,trim={0.25cm 0 0.75cm 0},clip]{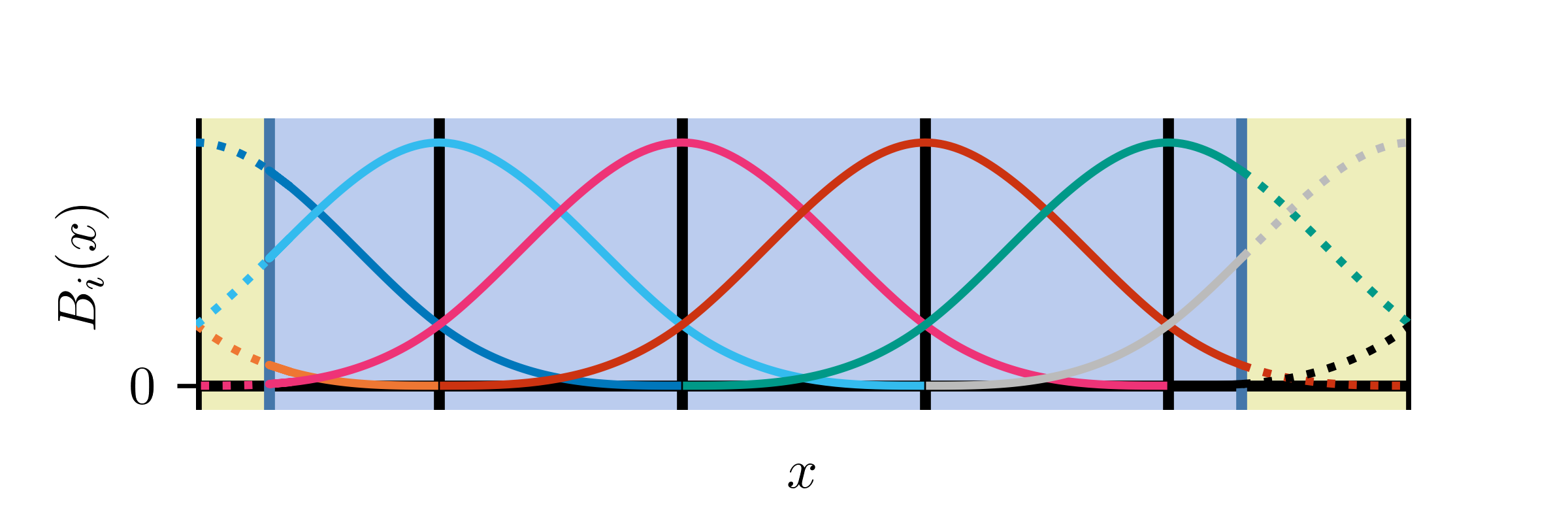}} 

    \caption{Optimal-regularity B-splines of various orders, $k$, constructed over an ambient mesh with five elements. The physical domain is shown in blue and the part of the basis functions without support in the physical domain is printed by dotted lines.}
	\label{fig:basisfunc}
\end{figure}

The rectilinear mesh constructed over the ambient domain, $\mathcal{T}_{\mathcal{A}}^h$, allows for the construction of multivariate B-spline basis functions as tensor products of univariate B-spline functions \cite{Cottrell:2009ad}. We herein consider optimal-regularity B-splines of order $k$ and regularity $k-1$, as illustrated in Fig.~\ref{fig:basisfunc} for univariate B-splines, $\{ B_{i} \in  C^{k-1}\}$, of various orders defined over an ambient domain with $\Nelem=5$ elements. Due to the $C^{k-1}$-continuity of optimal-regularity B-splines, the number of basis functions is generally substantially smaller than that of the corresponding Lagrange basis. For the illustrated univariate case, the number of optimal-regularity B-splines is equal to $\Nelem+k$, whereas the corresponding $C^0$-continuous Lagrange basis has $k\cdot \Nelem+1$ basis functions.

To restrict the approximation space to the physical domain, $\Omega$, only the basis functions with support over the physical domain are considered. We denote this set of functions by
\begin{align}
    \mathcal{S}^k_{k-1} = \big\{ N \in \{B_{i}\}  : \text{supp}(N) \cap \Omega \neq \emptyset \big\}
\end{align}
and the corresponding $\Ndofs$-dimensional approximation space by
\begin{equation}
\mathcal{V}^h = \text{span}\left( \mathcal{S}^k_{k-1} \right).\label{eq:Vh}
\end{equation}

\begin{remark}[Local refinements]
    A limitation of tensor-product B-splines, to which the discussions in this article are restricted, is that these cannot be refined locally. In the context of immersed isogeometric analysis we have found the extension with local refinement capabilities through (truncated) hierarchical B-splines particularly convenient \cite{Giannelli:2012rr,Brummelen:2021aw}.
    Locally refined meshes can then be constructed by sequential bisectioning of selections of elements in the mesh, starting from a rectilinear mesh, after which a truncated hierarchical B-spline basis can be constructed over the immersed domain. This procedure is, \emph{e.g.}, detailed in Ref.~\cite{divi_residual-based_2022} in the context of error estimation and adaptivity. 
\end{remark}

\subsection{Quadrature rules for non-boundary-fitted elements}\label{sec:integration}
\begin{figure}[!b]
    \centering
    \includegraphics[width=0.63\textwidth]{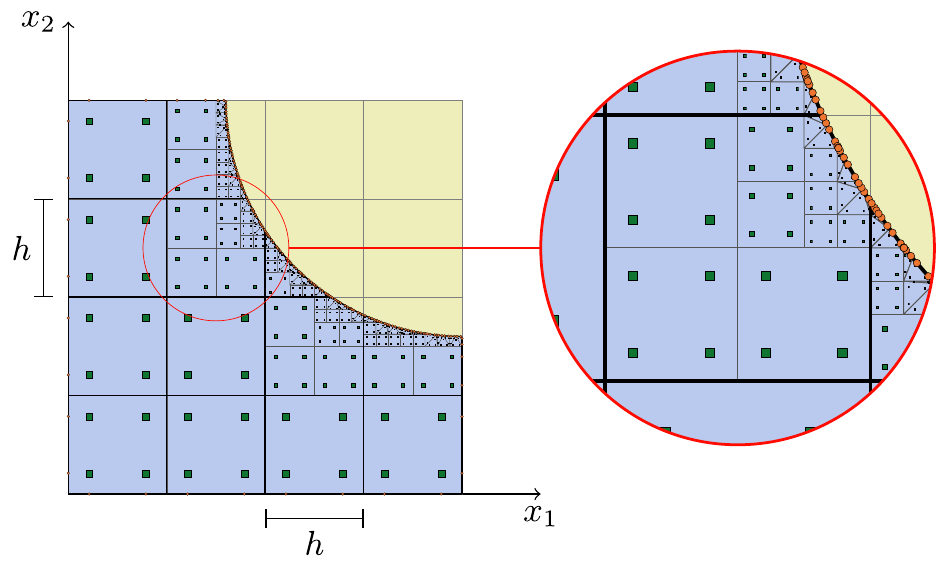}
    \caption{Volumetric (green squares) and surface (orange circles) quadrature rules obtained by the octree integration procedure with tessellation at the lowest bisectioning level.\\}
    \label{fig:octreeintegration}
\end{figure}

\noindent
The procedure used to construct integration quadrature on (a polygonal approximation of) the cut elements and their boundaries is illustrated in Fig.~\ref{fig:octreeintegration}. This procedure builds on the octree subdivision integration strategy described in Ref.~\cite{duster_finite_2008}, which is a widely used approach due to its simplicity and robustness. In the octree procedure, elements in the background mesh that intersect the boundary of the computational domain are bisected into $2^d$ integration sub-cells. If a sub-cell lies entirely within the domain, it is retained in the partitioning of the cut-element, whereas it is discarded if it lies entirely outside the domain. This bisectioning procedure is recursively applied to all the sub-cells that intersect the boundary. This recursion is terminated at a specified recursion depth by using the tessellation procedure proposed in Ref.~\cite{verhoosel_image-based_2015} (see Refs.~\cite{divi_error-estimate-based_2020,verhoosel2022scan} for implementation details). Through agglomeration of quadrature points on all of the sub-cells, cut-element integration rules can be constructed for both volumetric integrals (green squares in Fig.~\ref{fig:octreeintegration}) and immersed boundary integrals (orange spheres in Fig.~\ref{fig:octreeintegration}).

\section{Stabilized finite element formulation}\label{sec:WeakFormulation}

\noindent
The development of a Galerkin formulation of the model equations from \cref{sec:diffInt}, suitable for the discretization framework of \cref{sec:immersed}, requires careful consideration on two fronts:\\[-0.8cm]
\begin{enumerate} \setlength{\itemsep}{0pt}%
    \setlength{\parskip}{0pt}%
    \item Boundary condition enforcement: strong imposition is no-longer an option at immersed boundaries, where boundary conditions must instead be enforced weakly.
    \item Stabilization of the equal-order velocity-pressure pair and of basis functions with small support: when left untreated, small-cut elements lead to severely ill-conditioned (tangent) stiffness matrices.\\[-0.8cm]
\end{enumerate}
In \cref{ssec:weakform,ssec:nitsche} we propose a weak formulation where all boundary conditions are treated weakly, and in \cref{ssec:skelpen} this formulation is supplemented with edge stabilization terms to ensure inf-sup stability and small-cut insensitivity.

\subsection{Weak formulation}\label{ssec:weakform}
\noindent
To arrive at a weak formulation where the majority of the boundary conditions can be enforced naturally, we choose to use a mixed formulation of the Cahn-Hilliard equations, introducing the chemical potential $\mu$, as defined in equation \eqref{eq:chempot}, as an additional unknown field. If we assume that the solution fields that satisfy the strong form of \cref{eq:strong,eq:chempot} are sufficiently smooth, then they also satisfy the following weighted residual statement, obtained by multiplying each of \cref{eq:strong,eq:chempot} by a test function, integrating over the domain and applying integration by parts wherever appropriate:
\begin{subequations}\label{eq:weak1}
\begin{alignat}{3}
\begin{split}
    & \text{For a.e. }t\in (0,T)\text{ and } \forall \,\big(\vec{v},q,\omega,\lambda\big) \in \vec{H}^1(\Omega) \times L^2(\Omega) \times H^1(\Omega) \times H^1(\Omega): \nonumber
\end{split}\nonumber\\
\begin{split}&\int\limits_\Dom \Big\{  \partial_t ( \rho \u ) \cdot \vec{v} + \nabla \cdot( \rho \u\otimes \u )\cdot\vec{v} - \vec{J} \cdot \nabla \vec{v} \cdot \u+ \vec{\tau}:\nabla^s \vec{v} + \vec{\zeta}:\nabla^s\vec{v}  - p \nabla\cdot\vec{v} \Big\}\dDom \\[-0.5cm]
&\qquad\qquad + \int\limits_{\bdyDom} \Big\{\vec{J}\cdot\vec{n} \, \vec{u} \cdot\vec{v} + \big( - \vec{\tau}\vec{n}  - \vec{\zeta}\vec{n}  + p\vec{n} \big) \cdot\vec{v} \Big\} \dbdy = 0\,, 
\end{split}\\
& \int\limits_\Dom q \nabla\cdot \u \dDom  = 0\,,\\
& \int\limits_\Dom \Big\{ \partial_t \varphi \, \lambda + \nabla\cdot( \varphi \u ) \nabla \lambda  + m\nabla\mu\cdot\nabla\lambda \Big\} \dDom - \int\limits_{\bdyDom}m\nabla\mu \cdot\vec{n} \, \lambda \dbdy = 0\,,\\ 
& \int\limits_\Dom \Big\{ \mu\,\omega - \sigma\varepsilon \nabla\varphi\cdot\nabla \omega - \frac{\sigma}{\varepsilon} \Psi' \omega  \Big\} \dDom + \int\limits_{\bdyDom} \sigma \varepsilon \nabla\varphi\cdot\vec{n} \, \omega  \dbdy = 0\,.
\end{alignat}%
\end{subequations}

We may now substitute the generalized Navier boundary condition \cref{eq:GNBC}, the contact angle condition \cref{eq:DCA}, and the diffusive impermeability condition \cref{eq:nodiff} directly in the boundary integrals. The impermeability condition $\u\cdot\vec{n} = 0$ from \cref{eq:noconv} remains an essential condition, for which weak imposition requires specialized treatment (to be discussed in \cref{ssec:nitsche}). In the discrete form this impermeability condition will then not be satisfied pointwise along the immersed boundaries. To avoid spurious in- or outflow of mass and kinetic energy ensuing from discretization errors, some manipulation of the weak formulation is still in order. We negate the spurious kinetic energy flux through impermeable domain boundaries by introducing the skew-symmetric form of the non-linear advective term:
\begin{align}
\begin{split}
&\int\limits_\Dom  \nabla \cdot( \rho \u\otimes \u )\cdot\vec{v} \dDom = \int\limits_\Dom \Big\{  \frac{1}{2} \nabla \cdot( \rho \u\otimes \u \cdot\vec{v} ) - \frac{1}{2} \rho \u \cdot \nabla \vec{v} \cdot \u + \frac{1}{2} \nabla \rho \cdot \u \, \u \cdot\vec{v} \\
& \qquad + \frac{1}{2} \rho \nabla \cdot \u  \,  \u \cdot\vec{v} + \frac{1}{2} \rho \u \cdot \nabla \u \cdot\vec{v} \Big\} \dDom = \int\limits_{\partial\Omega_\text{out}} \u\cdot\vec{n} \, \frac{1}{2} \rho  \u \cdot\vec{v} \dbdy \\
& \qquad + \int\limits_\Dom \Big\{  - \frac{1}{2} \rho \u \cdot \nabla \vec{v} \cdot \u + \frac{1}{2} \nabla \rho \cdot \u \, \u \cdot\vec{v} + \frac{1}{2} \rho \u \cdot \nabla \u \cdot\vec{v} \Big\} \dDom \,. \label{skewsym}
\end{split}
\end{align}
In the last line, the divergence free condition $\nabla\cdot\u=0$ is substituted, as well as $\vec{v}=\vec{0}$ on (conforming) inflow boundaries and $\u\cdot\vec{n}=0$ on (immersed) domain boundaries representing impermeable walls.

Similarly, spurious mass flux through immersed boundaries can be negated by manipulation of the convective term in the order-parameter transport equation
\begin{align}
&\int\limits_\Dom \nabla\cdot( \varphi \u ) \lambda \dDom = -\int\limits_\Dom  \varphi \u \cdot \nabla \lambda \dDom + \!\!\!\!\!\! \int\limits_{\partial\Omega_\text{in}\cup\partial\Omega_\text{out}}  \!\!\!\!\!\!  \u\cdot\vec{n} \, \varphi \lambda \dbdy \,,\label{masscons}
\end{align}
where $\u\cdot\vec{n}$ is also substituted on impermeable domain boundaries.

Combining \cref{eq:weak1,skewsym,masscons} and substituting all appropriate boundary conditions then leads to:
\begin{subequations}\label{eq:weak}
\begin{alignat}{3}
\begin{split}
    & \text{For a.e. }t\in (0,T) \text{, find }U = (\u,p,\varphi,\mu)(t)\in \vec{V}_{\!\!\u_\text{in}}(\Omega) \times L^2(\Omega) \times V_{\varphi_\text{in}}(\Omega) \times H^1(\Omega) \\
    &\qquad \text{ s.t. } \forall \,(\vec{v},q,\omega,\lambda) \in \vec{V}_{\!\!\vec{0}}(\Omega) \times L^2(\Omega) \times V_{0}(\Omega) \times H^1(\Omega) :
\end{split}\nonumber\\
\begin{split}&r_{\vec{v}}(U,\vec{v}) = \int\limits_\Dom \Big\{  \partial_t ( \rho \u ) \cdot \vec{v} - \frac{1}{2} \rho \u\cdot \nabla\vec{v} \cdot \u  + \frac{1}{2}  \nabla\rho \cdot \u  \, \u\cdot \vec{v} + \frac{1}{2} \rho \u\cdot \nabla\u \cdot \vec{v} \\[-0.2cm]
&\qquad - \vec{J} \cdot \nabla \vec{v} \cdot \u  + \vec{\tau}:\nabla^s \vec{v} + \vec{\zeta}:\nabla^s\vec{v}  - p \nabla\cdot\vec{v} \Big\}\dDom + \!\! \int\limits_{\partial\Omega_\text{out}} \!\! \u\cdot\vec{n} \, \frac{1}{2} \rho  \u \cdot\vec{v} \dbdy \\[-0.2cm]
&\qquad  + \int\limits_{\bdyDom_\text{wall}} \Big\{  \aGN \, (\u-\u_s)\cdot \vec{v} - \nabla_\Gamma \sigma_{\SF}(\varphi) \cdot \vec{v} \Big\} \dbdy = 0 \end{split}\\[-0.1cm]
& r_{q}(U,q) = \int\limits_\Dom q \nabla\cdot \u \dDom  = 0 \label{eq:weak_mass}\\
& r_{\omega}(U,\omega) = \int\limits_\Dom \Big\{ \partial_t \varphi \, \lambda - \varphi \u \cdot\nabla \lambda  + m\nabla\mu\cdot\nabla\lambda \Big\} \dDom + \!\!\!\!\!\! \int\limits_{\partial\Omega_\text{in}\cup\partial\Omega_\text{out}}  \!\!\!\!\!\!  \u\cdot\vec{n} \, \varphi \lambda \dbdy = 0 \label{eq:WFphitrans}\\
& r_{\lambda}(U,\lambda) = \int\limits_\Dom \Big\{ \mu\,\omega - \sigma\varepsilon \nabla\varphi\cdot\nabla \omega - \frac{\sigma}{\varepsilon} \Psi' \omega  \Big\} \dDom - \int\limits_{\bdyDom_\text{wall}} \sigma'_{\SF}(\varphi) \, \omega  \dbdy = 0 \label{eq:WFchempot}
\end{alignat}
\end{subequations}
Let us note that in~\EQ{weak1} we have adopted a formal functional setting, to avoid the many technical complications associated with a rigorous formulation.

In the above formulation, the tangential projection operators are removed from the generalized Navier boundary terms since the normal-flow essential condition is imposed directly on the function spaces, as is the inflow condition,  per
\begin{subequations}
\begin{alignat}{3}
& \vec{V}_{\!\!\vec{g}}(\Omega) = \big\{ \vec{v} \in \vec{H}^1(\Omega) : \vec{v}=\vec{g} \text{ on }\partial\Omega_\text{in},\,  \vec{v}\cdot\vec{n} =0 \text{ on } \partial\Omega_\text{wall}\big\} \,, \\
& V_{g}(\Omega) = \big\{ \omega \in H^1(\Omega) : \omega = g  \text{ on }\partial\Omega_\text{in} \big\} \,.
\end{alignat}
\end{subequations}

\begin{remark}\label{rmk:testfuncs}
It should be noted that the test function~$\lambda$ for the transport equation for the phase field~\EQ{WFphitrans}, pairs with the trial function for the chemical potential, $\mu$, and the test function~$\omega$ for the chemical-potential closure relation~\EQ{WFchempot}, in fact pairs with the trial function for the order parameter, $\varphi$. This apparent asymmetry in the formulation is consistent with the fact that~\EQ{WFphitrans} (resp.~~\EQ{WFchempot}) contains the principal part of the operator acting on~$\mu$ (resp.~$\varphi$). This becomes important when considering the spaces which accommodate the variables for which essential conditions are imposed strongly: $\varphi,\omega \in V_{\bullet}(\Omega)$ and $\mu,\lambda\in H^1(\Omega)$. These choices of spaces then permit to select $\lambda = 1$, for which \cref{eq:WFphitrans} reduces to $\partial_t \int_\Omega \varphi \dDom = - \int_{\partial\Omega_\text{in}\cup\partial\Omega_\text{out}} \u\cdot\vec{n} \, \varphi \dbdy$, which, even in discrete form, signifies exact conservation.
\end{remark}

\subsection{Nitsche's method for imposition of convective impermeability}\label{ssec:nitsche}
\noindent
We discretize the weak formulation \eqref{eq:weak} in space using generally non-boundary-fitted optimal-regularity B-splines for all field variables, that is
\begin{align}
    \u^h &\in [\mathcal{V}^h]^d_{\u_\text{in}}, & p^h &\in \mathcal{V}^h, & \varphi^h &\in \mathcal{V}^h_{\varphi_\text{in}}, & \mu^h &\in \mathcal{V}^h,
    \label{eq:trialspaces}
\end{align}
with the space $\mathcal{V}^h$ as defined in equation \eqref{eq:Vh}. We assume that the inflow boundary does, in fact, align with the background mesh, whereby the essential inflow conditions $\vec{u}=\vec{u}_\text{in} $ and $\varphi = \varphi_\text{in}$ \textit{can} be imposed strongly:
\begin{subequations}
\begin{alignat}{3}
[\mathcal{V}^h]^d_{\vec{g}} &= \big\{ \vec{v} \in [\mathcal{V}^h] : \vec{v}=\vec{g} \text{ on }\partial\Omega_\text{in} \big\} \,, \\
\mathcal{V}^h_{g} &= \big\{ \omega \in \mathcal{V}^h : \omega = g  \text{ on }\partial\Omega_\text{in} \big\}\,.
\end{alignat}
\end{subequations} 
For the remaining essential condition, $\vec{u}\cdot\vec{n} =0$ on  $\partial\Omega_\text{wall}$, we propose a variant of Nitsche's method.

Nitsche's method combines penalty enforcement of a condition with a consistency term and an appropriate symmetry term. For multi-field non-linear equations, care must be taken to assure that the penalty term provides a sufficient bound on the consistency and symmetry terms, relating to the inf-sup stability of the resulting form \cite{Burman2015} as well as energy dissipation rates in the thermodynamic analysis \cite{Douglas1975,Simsek:2018gb}. The consistency term arises naturally in the weak form when the test functions do not vanish on the immersed boundary where the constraint condition is prescribed. From \cref{eq:weak1}, it may be inferred that the consistency term comprises normal components of the traction vector due to the diffusive, capillary, and pressure stresses, and the relative mass flux. The penalty and consistency terms then become:
\begin{subequations}
\begin{alignat}{3}
    s^{\rm pen}_{\vec{v}} ((\u^h,p^h,\mu^h,\varphi^h), \vec{v}) &:= \!\! \int\limits_{\bdyDom_\text{wall}} \!\! \beta \, h^{-1}  \eta(\varphi) \, (\u^h\cdot\vec{n})\, (\vec{v}\cdot\vec{n}) \dbdy \,,\label{eq:nitsche_pen}\\ 
    \tilde{s}^{\,\rm cons}_{\vec{v}} ((\u^h,p^h,\mu^h,\varphi^h), \vec{v}) &:= \!\!\int\limits_{\bdyDom_\text{wall}} \!\! (\vec{v}\cdot\vec{n}) \big\{ \u^h\cdot \vec{n}  \vec{J}\cdot \vec{n} - \vec{n} \cdot \vec{\tau} \vec{n} - \vec{n} \cdot \vec{\zeta} \vec{n}  + p^h \big\} \dbdy \,.
\label{eq:nitsche_cons}
\end{alignat}
\end{subequations} 
The penalty on the normal component of the velocity scales with the $\varphi$-dependent viscosity, $\eta$, and the inverse of the (background) mesh size, $h$. The penalty parameter, $\beta$, should be selected large enough to ensure stability of the formulation, but not so large as to induce locking-type effects (see, \emph{e.g.}, Refs.~\cite{de_prenter_note_2018,badia_mixed_2018,deprenter2023stability}). 

Upon substitution of the closure relations \eqref{eq:closureeqs}, the consistency term of the Nitsche operator can be simplified to
\begin{equation}
\begin{aligned}
    s^{\rm cons}_{\vec{v}}((\u^h,p^h,\mu^h,\varphi^h), \vec{v}) &:= \!\! \int\limits_{\bdyDom_\text{wall}} \!\! (\vec{v}\cdot\vec{n}) \big\{ p^h - \vec{n}\cdot (\eta \nabla^s \u^h ) \cdot\vec{n} - \frac{\sigma\epsilon}{2}|\nabla \varphi^h|^2 - \frac{\sigma}{\epsilon}\Psi\big\} \dbdy \, ,
\end{aligned}    
\label{eq:nitscherewritten}
\end{equation}
where the term $\sigma \varepsilon \, \vec{n}\cdot\nabla \varphi \otimes \nabla \varphi \cdot \vec{n} = \sigma \varepsilon \big( \nabla \varphi \cdot \vec{n} \big)^2$ emerging from $\vec{n}\cdot\boldsymbol{\zeta}\vec{n}$ per \cref{eq:closureeqszeta} has vanished on account of the stationary neutral wetting  
condition~\cref{eq:dnphi=0}, and the normal component of the relative mass flux, $\boldsymbol{J} \cdot \vec{n}$, has vanished on account of the homogeneous Neumann boundary condition for the chemical potential \eqref{eq:nodiff}.

The terms that remain from the capillary stress contribution precisely represent the mixture energy: $E_\text{mix} = \frac{\sigma\epsilon}{2}|\nabla \varphi|^2 + \frac{\sigma}{\epsilon}\Psi$. In a discrete dissipation analysis, following from $\vec{v}=\u^h$, the last two terms in \cref{eq:nitscherewritten} cancel out the in- or outflux of mixture energy at impermeable walls and thus avoid spurious inflow of mixture energy. Symmetrizing these terms would reintroduce spurious influx terms. Instead, we choose to only add symmetry terms for the pressure and viscous stress tensor contributions:
\begin{subequations}
\begin{alignat}{3}
    s^{\rm sym}_{\vec{v}}((\u^h,p^h,\mu^h,\varphi^h), \vec{v}) &:=  - \int\limits_{\bdyDom_\text{wall}}\!\! (\u^h\cdot\vec{n}) \, \vec{n}\cdot(\eta \nabla^s \vec{v} ) \cdot \vec{n} \dbdy\,, \\
    s^{\rm sym}_{q}(\u^h, q) &:=  \int\limits_{\bdyDom_\text{wall}} \!\! \u^h\cdot\vec{n} \, q \dbdy \,.
\end{alignat} \label{nitsche_sym}%
\end{subequations} 
Combining \cref{eq:nitsche_pen,eq:nitscherewritten,nitsche_sym} yields the final form of the Nitsche operator as
\begin{align}
\begin{split}
   s^{\rm nitsche}(U^h, (\vec{v},q)) & = \underbrace{s^{\rm pen}_{\vec{v}} (U^h, \vec{v}) + s^{\rm cons}_{\vec{v}} (U^h, \vec{v}) + s^{\rm sym}_{\vec{v}} (U^h, \vec{v})}_{s^{\rm nitsche}_{\vec{v}}(U^h, \vec{v})} + \underbrace{s^{\rm sym}_{q}(\u^h, q)}_{s^{\rm nitsche}_{q}(\u^h, q)}\,,
\end{split}
\end{align}
with $U^h = (\u^h,p^h,\mu^h,\varphi^h)$ the complete solution tuple.

\subsection{Ghost- and skeleton-penalty stabilization}\label{ssec:skelpen}
\begin{figure}[!b]
    \centering
    \subfloat[First-order derivative.]{\includegraphics[width=0.33\linewidth,trim={0.25cm 0 0.75cm 0},clip]{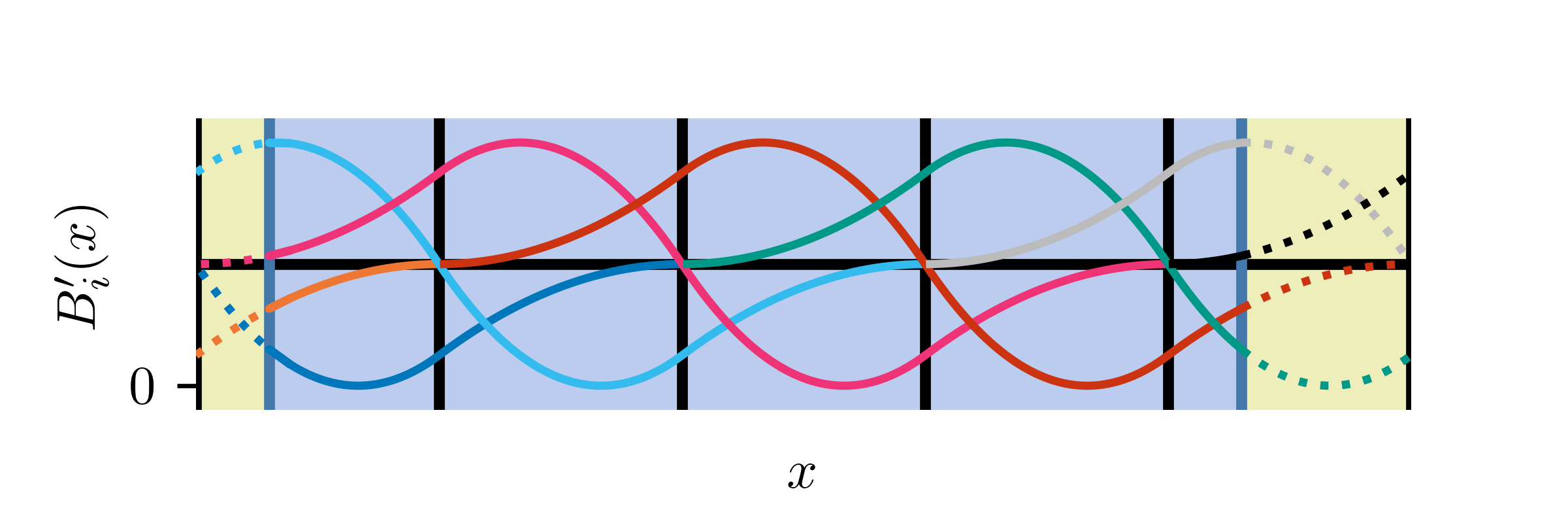}} \hfill
    \subfloat[Second-order derivative.]{\includegraphics[width=0.33\linewidth,trim={0.25cm 0 0.75cm 0},clip]{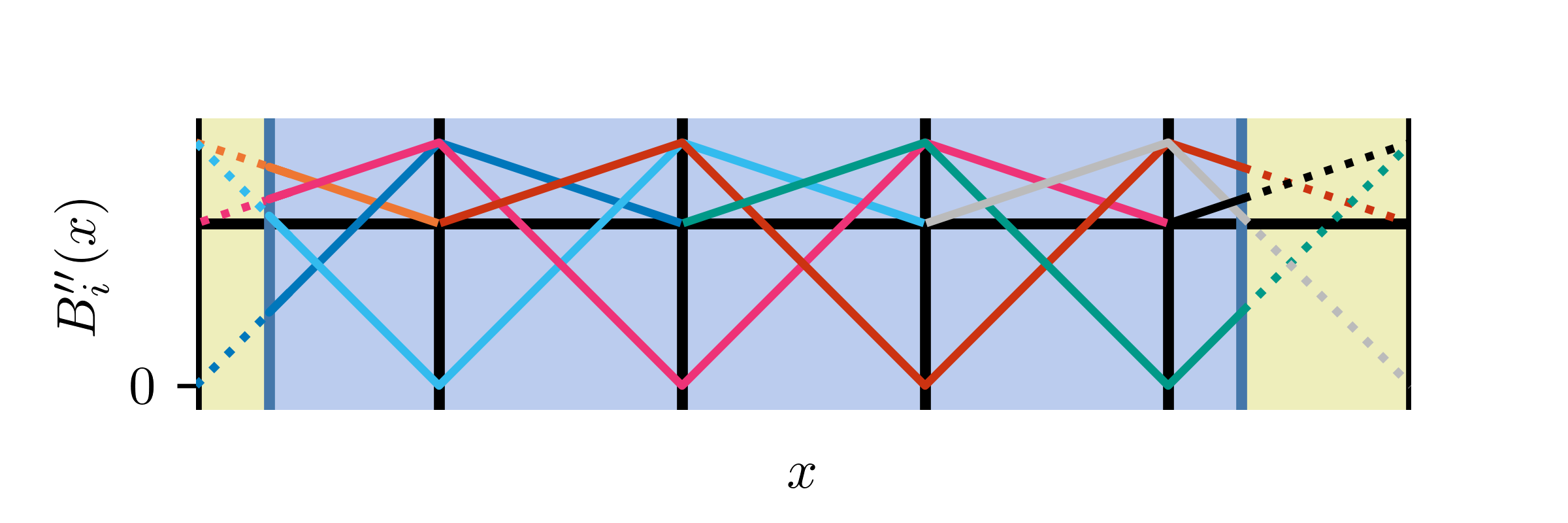}} \hfill
    \subfloat[Third-order derivative.]{\includegraphics[width=0.33\linewidth,trim={0.25cm 0 0.75cm 0},clip]{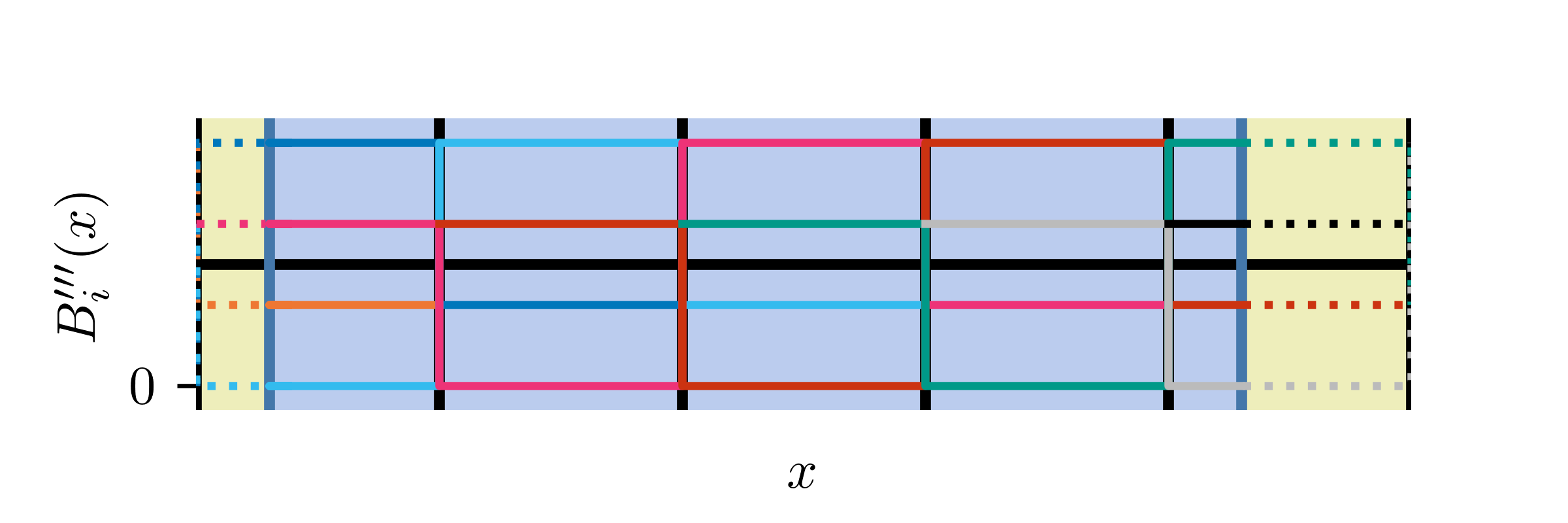}} \hfill
    \caption{Derivatives of the cubic ($k=3$) optimal-regularity B-spline basis shown in Fig.~\ref{fig:basisfunc}. Only the highest-order derivative is discontinuous over the skeleton mesh, which in this univariate illustration corresponds to the interior element boundaries.}
	\label{fig:derivatives}
\end{figure}

\noindent
Since the equal-order spaces \eqref{eq:trialspaces} are not inf-sup stable for the velocity-pressure pair, without further stabilization, spurious pressure oscillations will occur. We therefore augment the Galerkin formulation with the skeleton-stabilization operator proposed in Ref.~\cite{hoang_skeleton-stabilized_2019}, which reads
\begin{equation}
 s^{\rm skeleton}_{q}(p^h, q) = \int\limits_{\mathcal{F}^h_{\rm skeleton}}  \!\!\!  \gamma_{\rm skeleton} \, h^{2k+1}\eta^{-1} \jump{ \partial^k_n p^h  } \jump{ \partial^k_n q } \dbdy .\label{eq:skeletonstabilization}
\end{equation}
In this expression, $\partial_n^k(\cdot)$, represents the $k$-th order normal derivative, and $\jump{\cdot}$ is the interface-jump operator. Note that, for optimal-regularity B-splines, all derivative-orders lower than $k$ are continuous across the skeleton interfaces; see Fig.~\ref{fig:derivatives}. The pressure stabilization \eqref{eq:skeletonstabilization}, which acts on the complete skeleton, $\mathcal{F}_{\rm skeleton}^h$, penalizes jumps in higher-order pressure gradients, and can be regarded as the higher-order continuous version of the interior penalty method proposed in Ref.~\cite{burman_edge_2006}. To ensure stability and optimality, the operator \eqref{eq:skeletonstabilization} must scale with the size of the faces as $h^{2k+1}$. The parameter $\gamma_{\rm skeleton}$ should be set large enough to suppress pressure oscillations, and small enough to limit its influence on the accuracy of the solution.

To avoid stability problems associated with small or unfavorably-cut elements, ghost-penalty stabilization is applied to all field variables, that is,
\begin{subequations}
\begin{align}
 s^{\rm ghost}_{\vec{v}}(\u^h, \vec{v}) &:=  \int\limits_{\mathcal{F}^h_{\rm ghost}} \!\!\!  \gamma_{\rm ghost}\, h^{2k-1}\eta \,\jump{\partial^k_n \u^h} \cdot \jump{\partial^k_n \vec{v}} \, \text{d} S,\\
 s^{\rm ghost}_{\omega}(\varphi^h, \omega) &:=  \int\limits_{\mathcal{F}^h_{\rm ghost}}  \!\!\!  \gamma_{\rm ghost} \, h^{2k-1}  \sigma\epsilon \jump{\partial^k_n \varphi^h} \jump{\partial^k_n \omega} \dbdy,\\
 s^{\rm ghost}_{\lambda}(\mu^h, \lambda) &:=  \int\limits_{\mathcal{F}^h_{\rm ghost}}  \!\!\! \gamma_{\rm ghost} \, h^{2k-1}\, m\, \jump{ \partial^k_n \mu^h } \jump{\partial^k_n \lambda  } \dbdy.
\end{align}
\label{eq:ghostoperators}%
\end{subequations}
Note that, since the ghost mesh is a subset of the skeleton mesh, the pressure field is already stabilized through the operator \eqref{eq:skeletonstabilization}. 

The ghost-penalty operators \eqref{eq:ghostoperators} control the $k^{\rm th}$-order normal derivative jumps over the interfaces of the elements which are intersected by the domain boundary. Since we consider splines of degree $k$ with $C^{k-1}$ continuity, only the jump in the $k^{\rm th}$ normal derivative is non-vanishing at the ghost mesh. This leads to a single penalization term, contrasting the case of Lagrange elements, where all derivatives are non-vanishing and require penalization. The ghost-stabilization terms are scaled with the size of the faces as $h^{2k-1}$. They are also scaled with the physical parameters of their corresponding (vector)-Laplace terms in the Galerkin problem \eqref{eq:abstractgalerkin}, which is also the reason for the $\varphi$-ghost-operator to be subtracted in \eqref{eq:abstractgalerkin} and relates to \cref{rmk:testfuncs}. The pairing with the Laplace operators allows for the application of the same ghost-penalty parameter, $\gamma_{\rm ghost}$, for all fields. Appropriate selection of this parameter assures the stability of the formulation independent of the cut-cell configurations. To avoid loss of accuracy, the ghost-penalty parameter, $\gamma_{\rm ghost}$, should also not be too large \cite{badia_linking_2022}.

\subsection{Concluding immersed isogeometric finite element formulation}
\noindent
In summary, we propose the following stabilized immersed isogeometric finite element formulation:
\begin{align}
\!\!\!\!\!\!\left\{ \begin{aligned} 
&\text{For a.e. }t\in (0,T) \text{, find } U^h = (\u^h,p^h,\varphi^h,\mu^h) \in [\mathcal{V}^h]^d_{\u_\text{in}} \times \mathcal{V}^h \times \mathcal{V}^h_{\varphi_\text{in}} \times \mathcal{V}^h  \text{ s.t.:} & &  \\
&r_{\vec{v}}(U^h, \vec{v}) + s^{\rm nitsche}_{\vec{v}}(U^h, \vec{v}) + s^{\rm ghost}_{\vec{v}}(\u^h, \vec{v})=0 & & \hspace{-3cm}\forall \boldsymbol{v} \in [\mathcal{V}^h]^d_{\vec{0}} \\
&r_{q}(p^h, q) - s^{\rm sym}_{q}(\u^h, q) + s^{\rm skeleton}_{q}(p^h, q)=0 & & \hspace{-3cm}\forall q \in \mathcal{V}^h \\
&r_{\lambda}(U^h, \lambda) + s^{\rm ghost}_{\lambda}(\mu^h, \lambda) =0 & & \hspace{-3cm}\forall \lambda \in \mathcal{V}^h  \\
&r_{\omega}(U^h, \omega) - s^{\rm ghost}_{\omega}(\varphi^h, \omega)=0 & & \hspace{-3cm}\forall \omega \in \mathcal{V}^h_0
\end{aligned}\right. \label{eq:abstractgalerkin}
\end{align}
In this formulation, the (generally nonlinear) operators $r_{\vec{v}}$, $r_{q}$, $r_{\lambda}$ and $r_{\omega}$ follow directly from the weak formulation \eqref{eq:weak}, the operators  $s^{\rm nitsche}_{\vec{v}}$ and $s^{\rm sym}_{q}$ from \cref{eq:nitsche_pen,eq:nitscherewritten,nitsche_sym}, and the stabilization operators $s^{\rm ghost}_{\vec{v}}$, $s^{\rm skeleton}_{q}$, $s^{\rm ghost}_{\lambda}$ and $s^{\rm ghost}_{\omega}$ from \cref{eq:skeletonstabilization,eq:ghostoperators}.

\section{Numerical experiments}\label{sec:numExp}

In this section we study the proposed stabilized immersed isogeometric analysis formulation for a series of test cases. In Section~\ref{sec:taylorcouette} we consider a binary-fluid Taylor-Couette flow, which we use to benchmark the framework against a conventional boundary-fitted finite element analysis. In Section~\ref{sec:inclusion} an idealized porous medium with a periodic microstructure is considered to demonstrate the modeling capabilities of the proposed framework. Finally, in Section~\ref{sec:porousmedium} a porous medium application is considered to demonstrate its ability to handle complex geometries.

\begin{table}
    
    \centering
    \small
    \begin{tabular}{lcrl}
        \toprule \multicolumn{4}{c}{\it Fluid-flow parameters} \\
        \toprule Mass densities & $\rho_{1}$ & 1000 & ${\rm kg/m^3}$ \\
         & $\rho_{2}$ & 1.3 & ${\rm kg/m^3}$ \\
        Dynamic viscosities & $\eta_{1}$ & $1 \cdot 10^{-3}$ & ${\rm Pa\, s}$ \\
         & $\eta_{2}$ & $1.813 \cdot 10^{-5}$ & ${\rm Pa\, s}$ \\ 
        Surface tension & $\sigma_{12}$ & $72.8\cdot10^{-3}$ & N/m\\
        Generalized Navier parameter & $\aGN$ & 100 & ${\rm Pa\,  s/m}$\\
        \bottomrule
            & & & \\[0.75em]
    \toprule \multicolumn{4}{c}{\it Phase-field parameters} \\
    \toprule Interface thickness & $\varepsilon$ & $0.78125\cdot10^{-6}$ & ${\rm m}$ \\
     Mobility & $m$ & $3.0487\cdot10^{-10}$ & ${\rm m\, s^2 /kg}$\\ 
    \bottomrule
        & & & \\[0.75em]
    \toprule \multicolumn{4}{c}{\it Stabilization parameters} \\
    \toprule 
     Nitsche penalty & $\beta$ & 100 & --\\
     Skeleton penalty & $\gamma_{\rm skeleton}$ & 0.01 & -- \\
     Ghost penalty & $\gamma_{\rm ghost}$ & 0.01 & -- \\
    \bottomrule
    \end{tabular}
    \caption{Properties of the Navier-Stokes-Cahn-Hilliard model used for the numerical experiments.}
    \label{tab:properties}
\end{table}

Unless specified otherwise, the parameters of the Navier-Stokes-Cahn-Hilliard model as listed in Table~\ref{tab:properties} are used for all test cases. These parameters -- which represent a water-air flow -- have a strong influence on the mesh and time resolution required to accurately evaluate the model. In particular, the interface thickness, $\varepsilon$, has a strong impact on the element size, and the mobility, $m$, has a strong influence on the time step size. The model parameters are selected such that stable and accurate results can be obtained on uniform background meshes with a moderate number of elements and with a moderate number of time steps. This enables studying the proposed immersed isogeometric analysis framework at an acceptable computational expense using various mesh and time step sizes, as well as studying its robustness with respect to cut-element configurations. The adaptive solution strategy for the NSCH equations presented in Ref.~\cite{Demont:2022dk} can be tailored to our immersed framework (a theoretical basis for this is provided in Ref.~\cite{divi_residual-based_2022}) to substantially reduce the computational effort associated with the uniform discretization considered here, but this extension is beyond the scope of this work. The test cases are restricted to the two-dimensional case and are implemented in the Python-based (isogeometric) finite element framework Nutils \cite{nutils}. Although there are no fundamental obstacles in extending the work to three dimensions -- both in terms of the formulation and in terms of the discretization method -- the computational burden associated with three-dimensional simulations would necessitate the usage of an adaptive solution strategy and the implementation in a high-performance computing framework (iterative solvers with preconditioners, scalable parallel implementation, \emph{etc.}).

The computational domains for the upcoming test cases are constructed by trimming the domain through the octree-tessellation procedure discussed in Section~\ref{sec:integration} based on analytical level set functions. The octree depth for all cases is set to 3, meaning that the octree recursion is terminated after 3 element bisections. Fifth-order Gaussian quadrature rules are selected on the integration sub-cells, as they are found to yield a suitable balance between minimizing computational expense and the impact of integration errors on the presented results. It is noted, however, that optimized integration rules can be used to enhance the performance of the framework \cite{divi_error-estimate-based_2020}.

For all immersed isogeometric simulations presented below, cubic ($k=3$) B-splines with optimal regularity are employed for all field variables. For the selection of the various stabilization parameters, we follow the empirical rules proposed in Ref.~\cite{hoang_skeleton-stabilized_2019}. The Nitsche parameter is set to 100 and both the skeleton-penalty parameter and the ghost-penalty parameter are set to 0.01 (see Table~\ref{tab:properties}). These parameters are found to provide a good balance between model stability (large enough) and impact on the accuracy of the approximation (small enough), for all considered simulations.

\subsection{Taylor-Couette flow}
\label{sec:taylorcouette}

\begin{figure}[!b]
    \centering
    \includegraphics[width=0.8\textwidth]{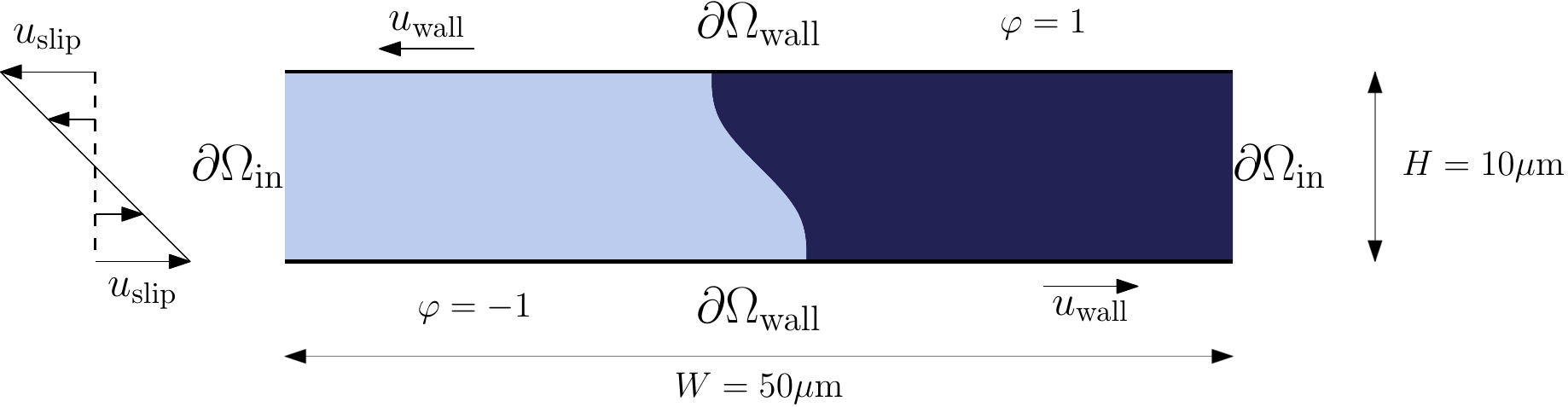}
    \caption{Illustration of the domain and boundary conditions for the binary-fluid Taylor-Couette benchmark case.}
    \label{fig:tcdomain}
\end{figure}

\noindent
We consider the binary-fluid Taylor-Couette flow between two parallel plates moving in opposite directions, as illustrated in Fig.~\ref{fig:tcdomain}, and studied in Refs.~\cite{Jacqmin:1999fk,Gerbeau:2009jx}. The fluid constituents are separated by a vertical interface in the center of the domain. To maintain symmetry, for this benchmark the mass density and viscosity are taken the same for both constituents, \emph{i.e.}, $\rho_1=\rho_2=\rho=1000{\rm kg/m^3}$ and $\eta_1=\eta_2=\eta=1\cdot 10^{-3}{\rm Pa\,s}$. Moreover, the mobility is decreased to $m=3.0487\cdot10^{-11}{\rm m\,s^2/kg}$. The speed of the plate is increased gradually over the first second per $u_{\rm wall}=\frac{1}{2}\big(1-\cos(\pi\,t)\big) 10\,{\rm m/s}$, after which it is kept constant at $u_{\rm wall}=10{\rm m/s}$ until a steady state solution is obtained. In the initial state, the vertical diffuse interface is prescribed using an analytical approximation in accordance with the selected model 
parameters. On the left and right (far field) boundaries, a linear velocity profile  is imposed, where $u_{\rm slip}=(1+\frac{2\eta}{\aGN H })^{-1}u_{\rm wall}$. Such a profile corresponds to the far-field (pure species) Taylor-Couette steady state solution compatible with the generalized Navier boundary condition. The phase field is set to +1 on the left boundary, and to -1 on the right boundary.

\begin{figure}
    \centering
\subfloat[The boundary-fitted mesh used to compute the finite element analysis reference result.]{\includegraphics[width=0.48\textwidth]{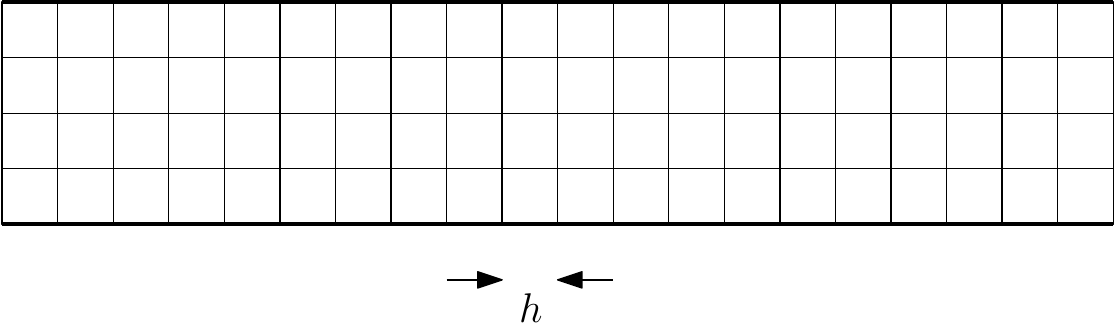}\label{fig:tcmesha}}\hfill
\subfloat[The background mesh for the immersed isogeometric analysis simulations. Rotated by an angle~$\theta$.]{\includegraphics[width=0.48\textwidth]{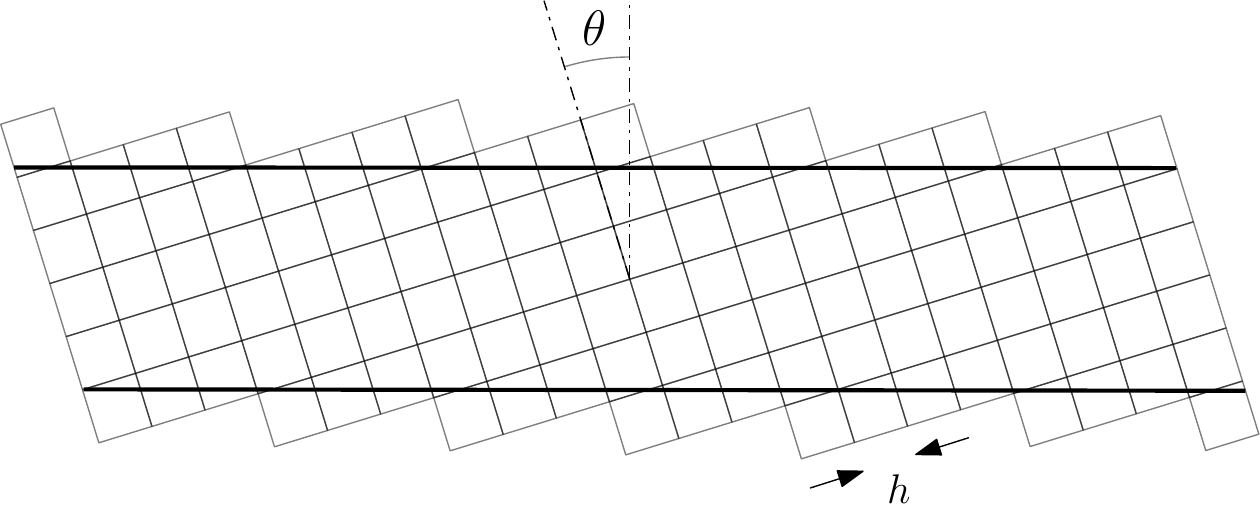}\label{fig:tcmeshb}}
    \caption{Illustration of the computational meshes used for the binary-fluid Taylor-Couette benchmark.}
    \label{fig:tcmesh}
\end{figure}

A boundary-fitted mesh can straightforwardly be constructed on the rectangular domain considered here (Fig.~\ref{fig:tcmesha}), making this test case suitable for benchmarking against a classical finite element discretization (Section~\ref{sec:tcbenchmark}). To study the proposed immersed isogeometric analysis method, the rectangular domain is immersed in a rotated ambient domain, as illustrated in Fig.~\ref{fig:tcmeshb}. By varying the rotation of the ambient domain, $\theta$, the robustness with respect to cut-element configurations can be studied (Section~\ref{sec:tcorientation}).

\begin{figure}
\begin{tabular}{cc}
\centering
\subfloat[$t=0$s]{\includegraphics[width=0.49\textwidth]{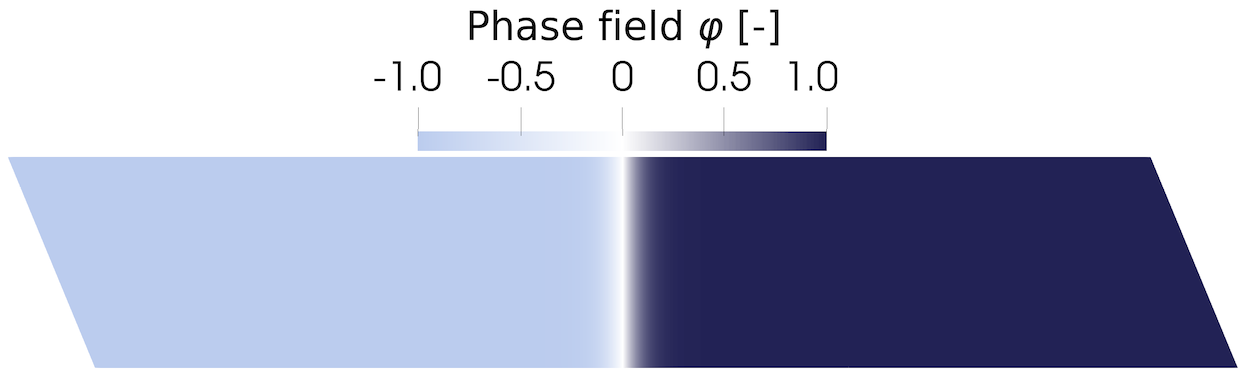}}
&
\subfloat[$t=0$s]{\includegraphics[width=0.49\textwidth]{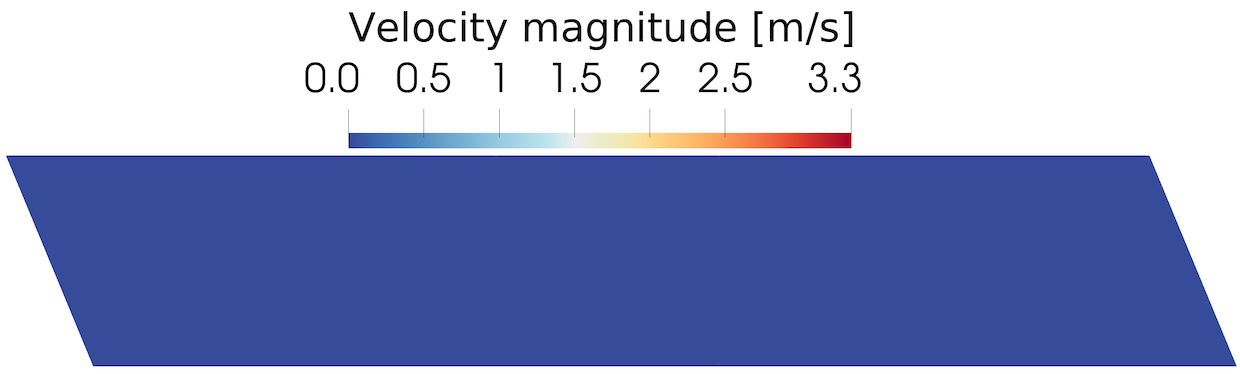}}
\\
\subfloat[$t=0.5$s]{\includegraphics[width=0.49\textwidth]{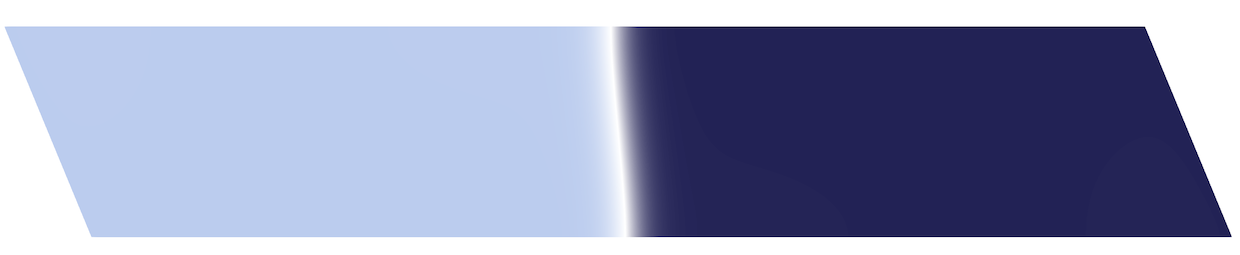}}
&
\subfloat[$t=0.5$s]{\includegraphics[width=0.49\textwidth]{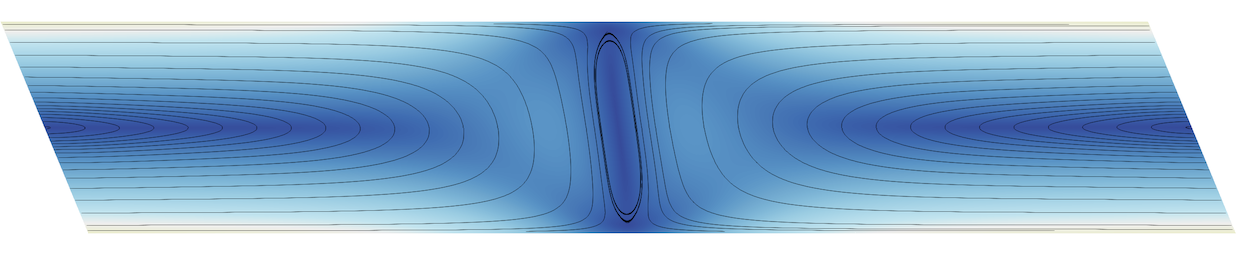}}
\\
\subfloat[$t=1$s]{\includegraphics[width=0.49\textwidth]{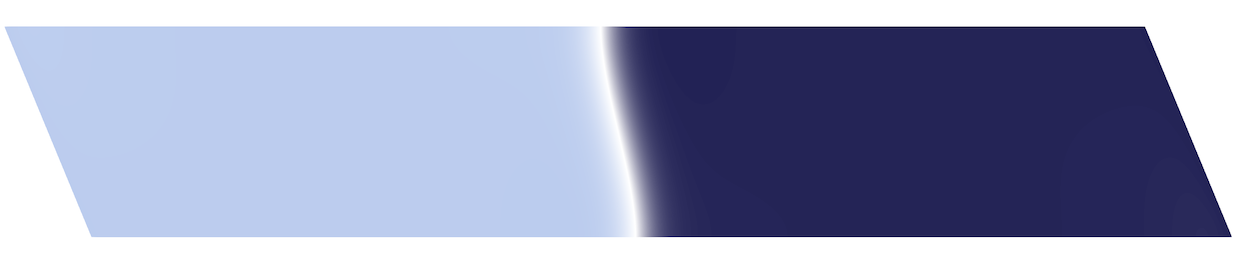}}
&
\subfloat[$t=1$s]{\includegraphics[width=0.49\textwidth]{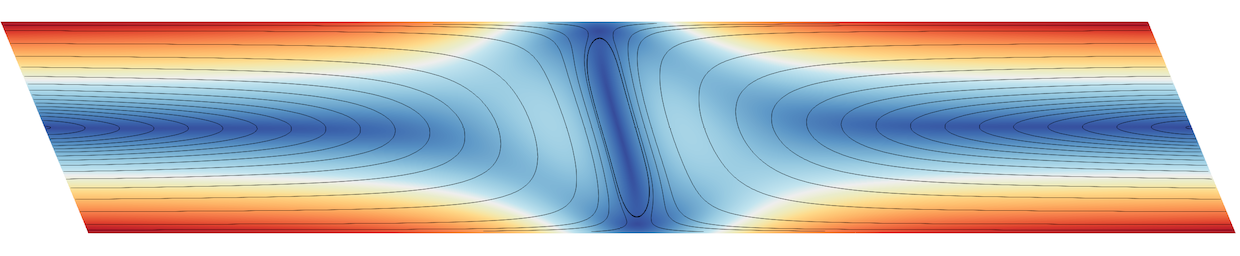}}
\\
\subfloat[$t=\infty$]{\includegraphics[width=0.49\textwidth]{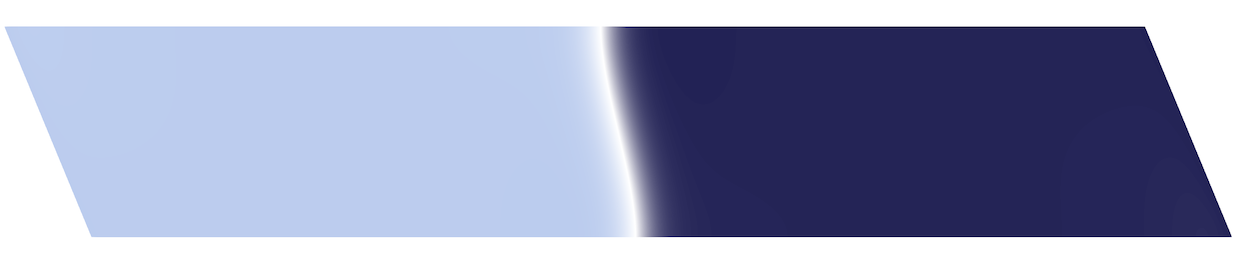}}
&
\subfloat[$t=\infty$]{\includegraphics[width=0.49\textwidth]{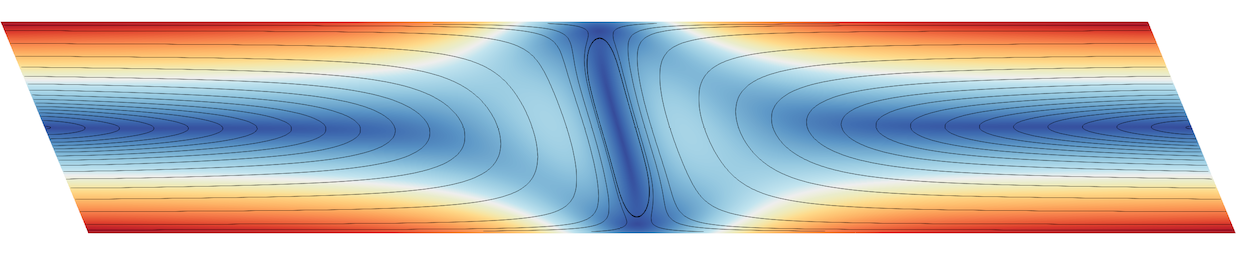}}
\\
\end{tabular}
    \caption{Time evolution of the phase field $\varphi$ (left) and the velocity field $\u$ (right) for the Taylor-Couette test case. The presented results are obtained using immersed isogeometric analysis with a $h=0.3125\mu\rm m$ ambient domain mesh rotated by $\theta = \frac{\pi}{8}{\rm rad}$.}
    \label{fig:tcsnapshots}
\end{figure}

To provide intuition for the studies presented in the remainder of this section, in Fig.~\ref{fig:tcsnapshots} we illustrate the evolution of the phase field (left column) and velocity field (right column) obtained using an ambient domain mesh with element size $h=0.3125{\rm \mu\rm m}$, which is rotated by $\theta=\frac{\pi}{8}{\rm rad}$ in counter clock-wise direction. The top row shows the solution in the initial state, where the diffuse interface is vertical and the velocity is zero. At $t=0.5{\rm s}$ (second row), the plates are moving, which results in the interface being dragged along with the plates, resulting in a (mildly) curved interface. The generalized Navier boundary condition is observable through the mismatch between the speed of the fluid at the plate boundaries (approximately $1.7{\rm m/s}$ away from the interface) and that of the plates ($u_{\rm wall}=5{\rm m/s}$). As time progresses, and the plate velocity is increased up to the maximum speed at $t=1.0{\rm s}$ (third row), the curvature of the interface increases, until a steady-state solution is obtained (bottom row). Note that the slip-velocity-based Reynolds number ${\rm Re} = 33$ is sufficiently small to allow for a steady state solution, and that convective effects are small enough as not to require the use of convection-stabilization. In this steady state, the maximum rotation of the interface relative to its (vertical) initial position is $0.23{\rm rad}$, and the tangential stresses due to capillary effects induce sufficient slip to equilibrate the meniscus-to-wall contact point. Away from this contact point, the classical wedge-shape Huh-Scriven velocity profiles are retrieved \cite{huh1971hydrodynamic}.

\subsubsection{Benchmarking against boundary-fitted finite elements}
\label{sec:tcbenchmark}
\noindent
To benchmark the proposed immersed isogeometric formulation, we compare the discretization described in \cref{sec:immersed,sec:WeakFormulation} to a
boundary-fitted simulation, for which all essential boundary conditions can be imposed strongly: since the plates are aligned with the horizontal direction, the normal component of the velocity can be constrained by removal of the basis functions corresponding to the vertical velocity component. Taylor-Hood $C^0$-continuous finite elements with cubic velocity functions and quadratic pressure functions are considered to obtain a stable velocity-pressure discretization. The phase field and chemical potential field are discretized using cubic finite elements. In this boundary-fitted finite element setting, no additional stabilization of the weak form \eqref{eq:weak} is required.

\begin{figure}
\begin{tabular}{cc}
\subfloat[Coarse mesh ($h=0.625{\rm \mu m}$) boundary fitted FEA with 51,778 dofs.]{\includegraphics[width=0.49\textwidth,trim={0 0 0 2.6cm},clip]{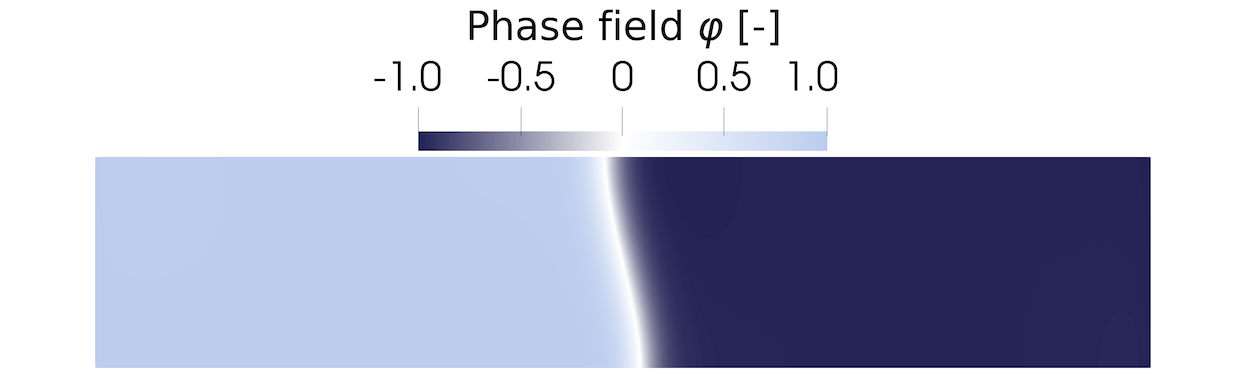}}
&
\subfloat[Coarse mesh ($h=0.625{\rm \mu m}$) immersed IGA with 9,384 dofs.]{\includegraphics[width=0.49\textwidth,trim={0 0 0 2.6cm},clip]{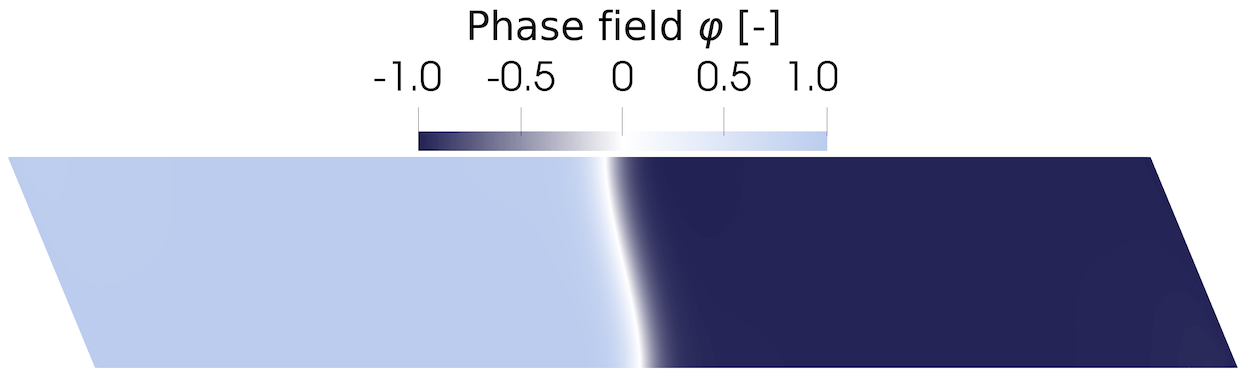}}
\\
\subfloat[Fine mesh ($h=0.3125{\rm \mu m}$) boundary fitted FEA with 205,954 dofs.]{\includegraphics[width=0.49\textwidth,trim={0 0 0 2.6cm},clip]{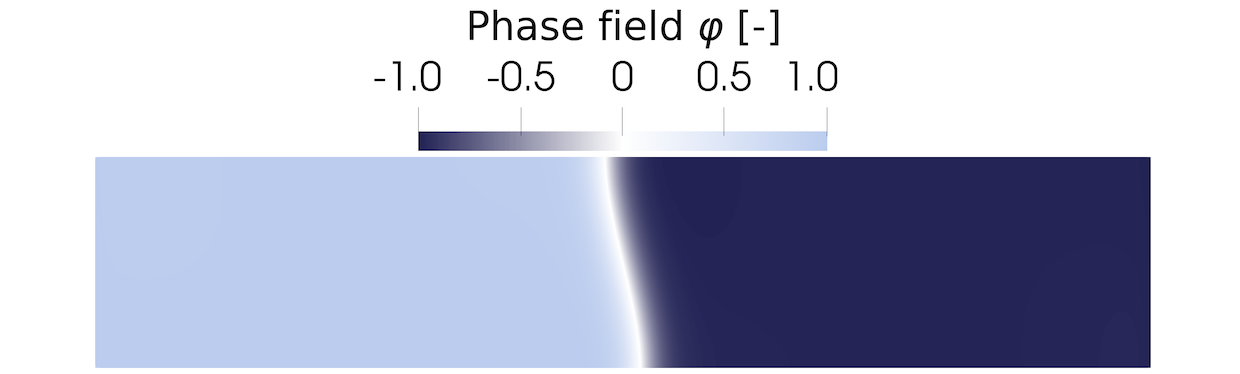}}
&
\subfloat[Fine mesh ($h=0.3125{\rm \mu m}$) immersed IGA with 32,582 dofs.]{\includegraphics[width=0.49\textwidth,trim={0 0 0 2.6cm},clip]{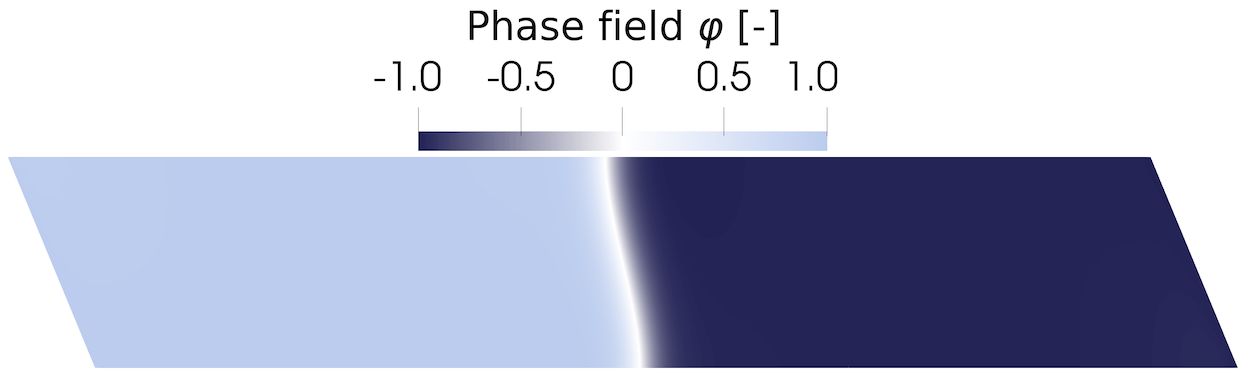}}
\\
\end{tabular}
    \caption{Comparison between a boundary fitted FEA and an immersed IGA of the steady state phase-field solution, $\varphi(t=\infty)$, for a coarse and a fine mesh.
    }\label{fig:tcbenchmarkphase}
\end{figure}

In Fig.~\ref{fig:tcbenchmarkphase}, the steady-state phase field obtained using traditional boundary-fitted finite element analysis (FEA) (left column) and the proposed immersed isogeometric analysis (IGA) (right column) are compared. In the top row, an element size of $h=0.625{\rm \mu m}$ is used, resulting in a mesh with $80 \times 16=1280$ elements for the boundary-fitted case. The bottom row presents the results for an element size of $h=0.3125{\rm \mu m}$, resulting in a boundary-fitted mesh with $160 \times 32=5120$ elements. For both meshes, the obtained phase-field solutions for the boundary-fitted and immersed simulations are virtually indistinguishable, demonstrating the consistency of the stabilized immersed IGA formulation. The maximum interface rotation angle for the immersed case is within a few percent of the FEA reference result. This difference can be attributed to the fact that the far-field in/out-flow conditions are imposed at a finite distance from the interface. The error associated with this inconsistent application of the far-field condition differs between the FEA case and IGA case on account of the size and rotation of the ambient domain (\emph{i.e.}, the boundaries are located in a different position). Due to the optimal-regularity B-splines used in the isogeometric analysis, the number of degrees of freedom associated with these simulations is substantially (approximately 6 times) smaller than that for the FEA case.

\begin{figure}
\begin{tabular}{cc}
\subfloat[Coarse mesh ($h=0.625{\rm \mu m}$) boundary fitted FEA with 51,778 dofs.]{\includegraphics[width=0.49\textwidth]{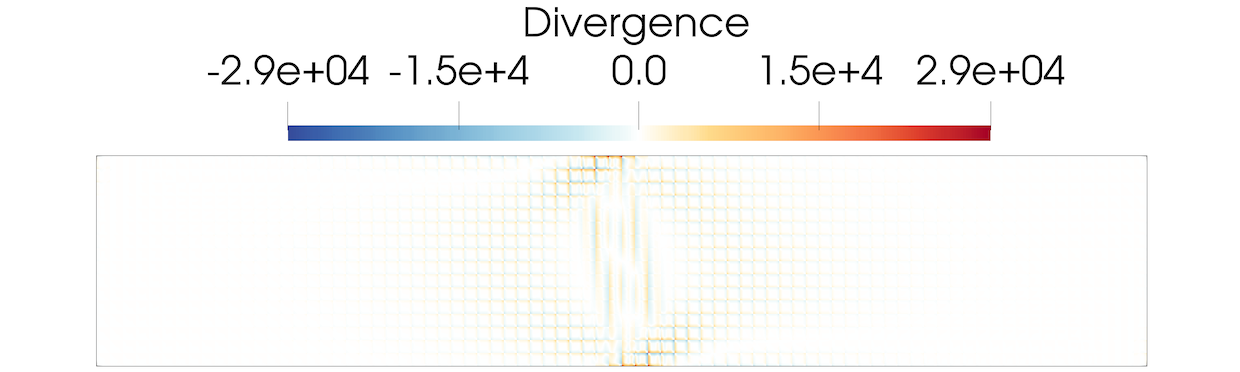}}
&
\subfloat[Coarse mesh ($h=0.625{\rm \mu m}$) immersed IGA with 9,384 dofs.]{\includegraphics[width=0.49\textwidth]{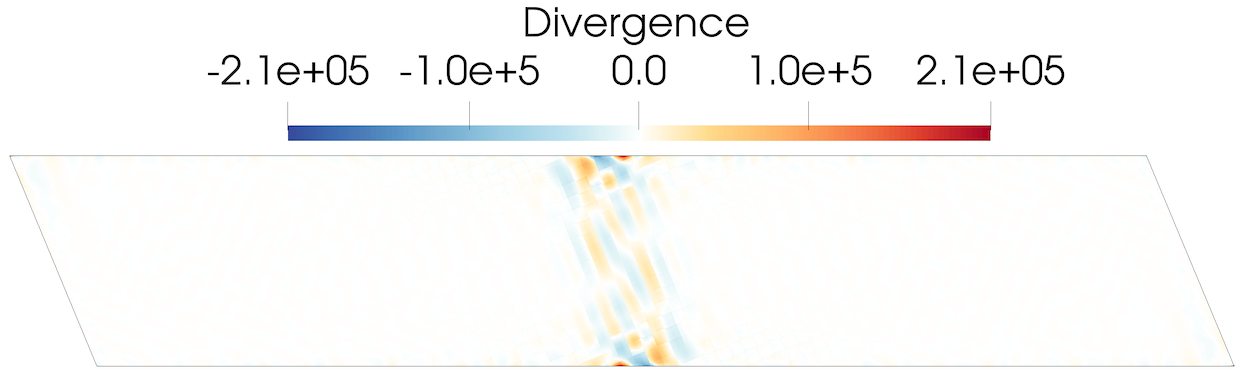}}
\\
\subfloat[Fine mesh ($h=0.3125{\rm \mu m}$) boundary fitted FEA with 205,954 dofs.]{\includegraphics[width=0.49\textwidth]{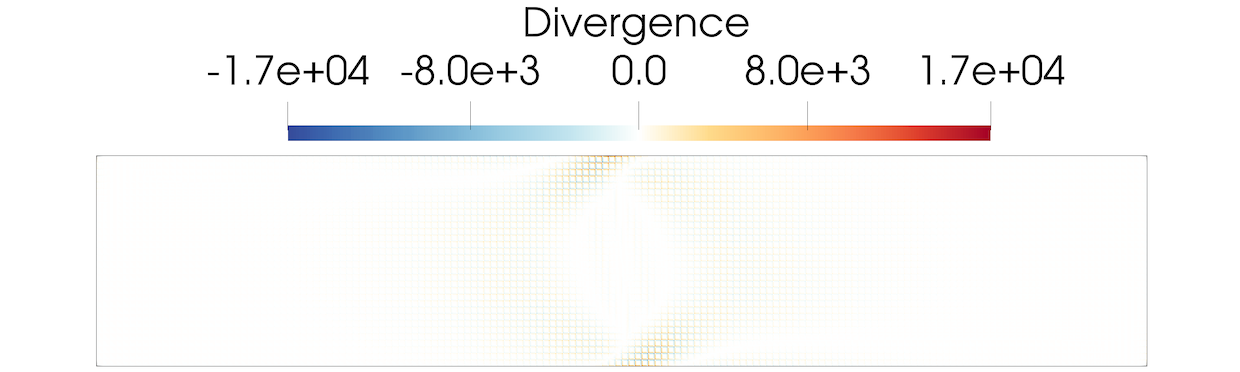}}
&
\subfloat[Fine mesh ($h=0.3125{\rm \mu m}$) immersed IGA with 32,582 dofs.]{\includegraphics[width=0.49\textwidth]{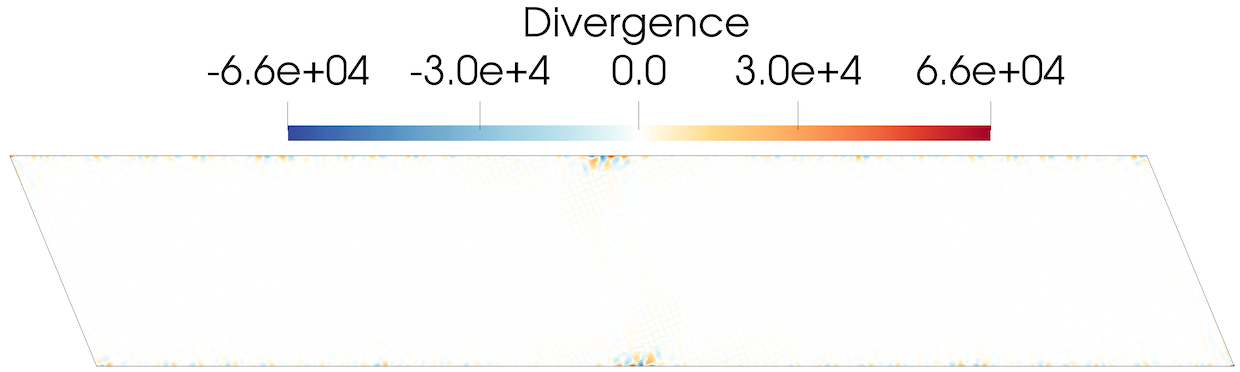}}
\\
\end{tabular}
    \caption{Comparison between a boundary fitted FEA and an immersed IGA of the steady state velocity-divergence solution, $\nabla\cdot\u(t=\infty)$, for a coarse and a fine mesh. Note that the color bars have different ranges.}
    \label{fig:tcbenchmarkdiv}
\end{figure}

Fig.~\ref{fig:tcbenchmarkdiv} compares the divergence of the velocity field between the traditional FEA case and the immersed IGA, considering the same meshes as for the phase-field solutions in Fig.~\ref{fig:tcbenchmarkphase}. Although, like the phase field, the velocity fields (not shown here) are virtually indistinguishable between the two methods, notable differences can be observed in the divergence of the velocity field. Since the mass conservation balance \eqref{eq:weak_mass} is solved weakly for both the FEA case and the immersed IGA case, discretization errors result in a non-zero divergence for both analyses, in particular in the vicinity of the interface. For both FEA and immersed IGA the error in the mass conservation decreases under mesh refinement, although the errors observed using traditional FEA are substantially smaller than for the immersed IGA case. We attribute this increased error in the velocity-divergence for the immersed IGA case to the reduced number of degrees of freedom in combination with the skeleton-penalty term required for stabilization. Although this reduction does not noticeably affect the phase field and the velocity field, it does impact the local conservation properties. Furthermore, it is observed that while the errors in the divergence are essentially restricted to the interface in the boundary-fitted case, in the immersed case errors are also observed along the immersed boundaries. This is due to the penalized and conflicting nature of the divergence-free condition and the convective impermeability condition in the cut elements. In our numerical experiments we have not experienced any negative consequences due to this phenomenon, but it should be recognized as a potential source of instabilities.

\subsubsection{Robustness with respect to cut-element configurations}
\label{sec:tcorientation}
\begin{figure}
\begin{tabular}{cc}
\centering
\subfloat[$\theta=0.001$rad]{\includegraphics[width=0.36\textwidth]{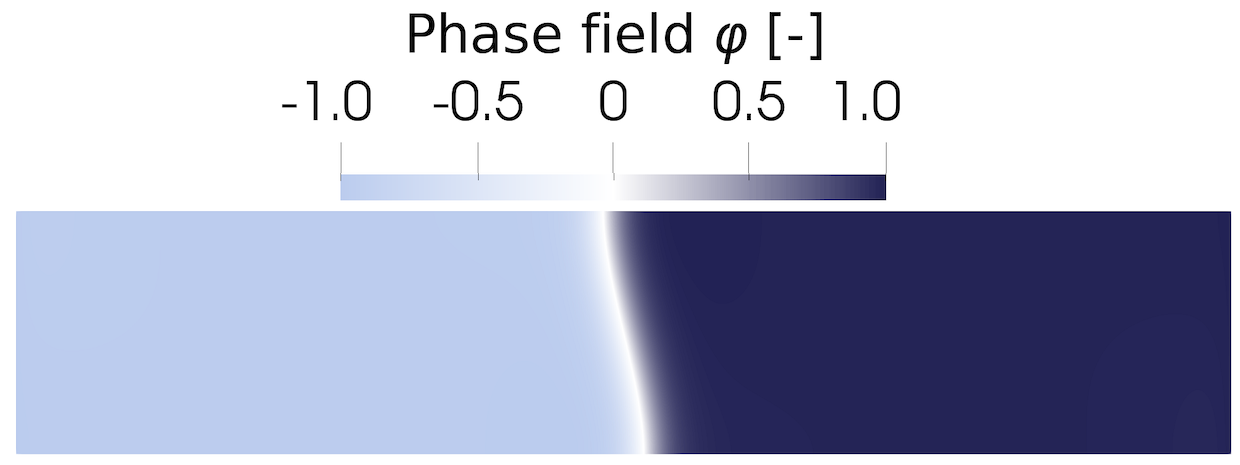}}
&
\subfloat[ $\theta=0.001$rad]{\includegraphics[width=0.36\textwidth]{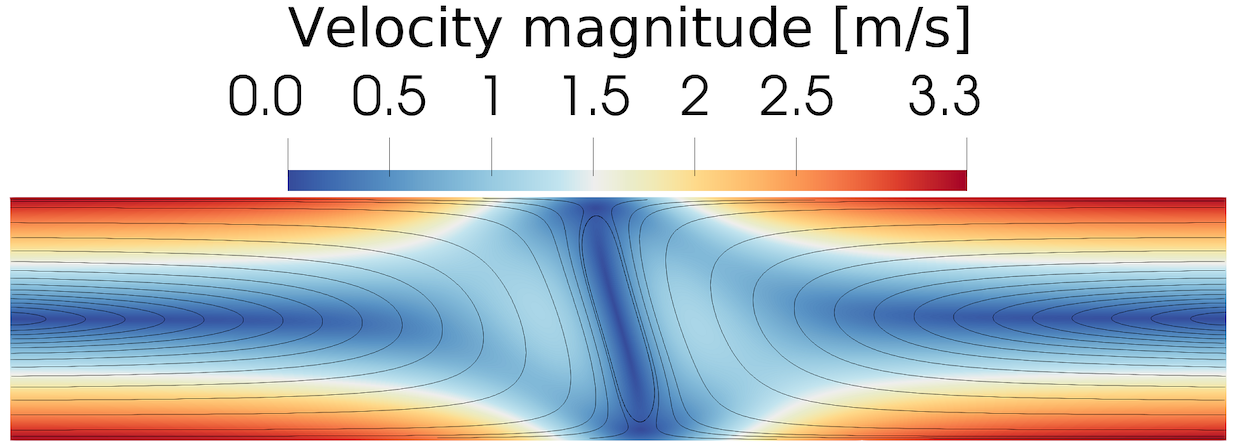}}
\\
\subfloat[$\theta=\frac{\pi}{8}$rad]{\includegraphics[width=0.41\textwidth]{figures/time_series/phi_TC_im_nx160_t1_25.png}}
&
\subfloat[$\theta=\frac{\pi}{8}$rad]{\includegraphics[width=0.41\textwidth]{figures/time_series/velocity_TC_im_nx160_t1_25.png}}
\\
\subfloat[$\theta=\frac{\pi}{4}$rad]{\includegraphics[width=0.49\textwidth]{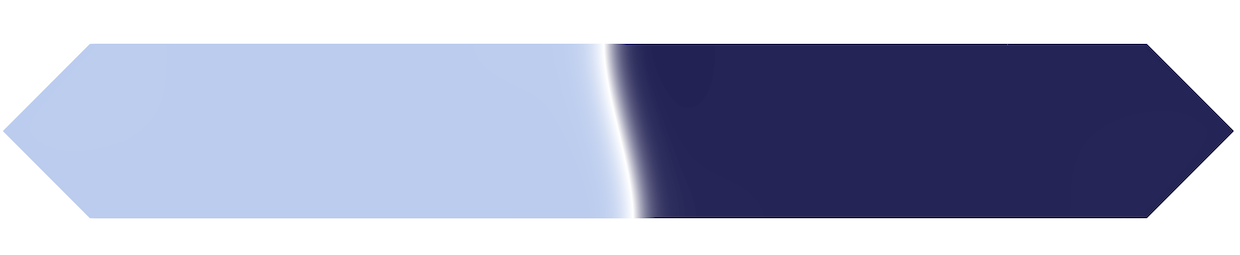}}
&
\subfloat[$\theta=\frac{\pi}{4}$rad]{\includegraphics[width=0.49\textwidth]{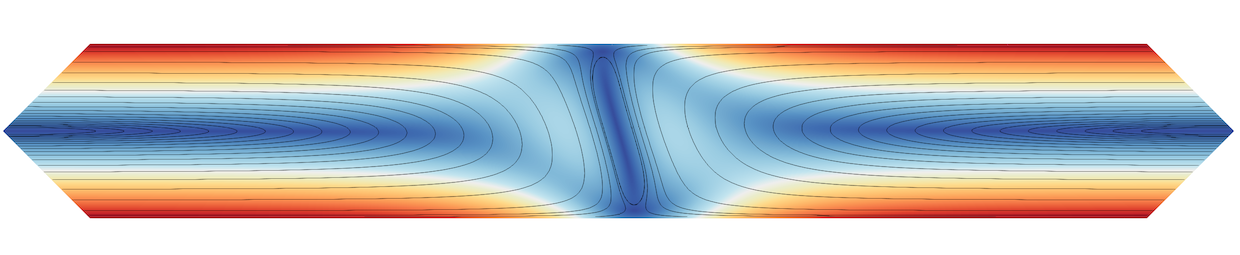}}
\\
\end{tabular}
    \caption{Steady-state phase-field, $\varphi(t=\infty)$, and velocity field, $\boldsymbol{u}(t=\infty)$, solutions computed with the immersed IGA framework on ambient domain meshes with a mesh size of $h=0.3125\mu\rm m$ and various mesh rotations $\theta$.}
    \label{fig:tcorientation}
\end{figure}

\noindent
To study the robustness of the proposed formulation in terms of how elements are cut, various rotations of the ambient domain mesh are considered. Fig.~\ref{fig:tcorientation} shows the phase field (left column) and velocity field (right column) for three different orientations of the ambient mesh. The considered mesh size is $h=0.3125\mu\rm m$. The top row concerns a very small rotation angle of 0.001rad. This mesh rotation leads to very thin sliver cuts along the boundaries. The presented results are virtually indistinguishable from the reference solution discussed above. The absence of instabilities associated with the Nitsche boundary condition and ill-conditioning problems \cite{deprenter2023stability} conveys that the ghost-penalty stabilization is effective. The middle and bottom rows of Fig.~\ref{fig:tcorientation} pertain to the case of a $\frac{\pi}{8}$rad and $\frac{\pi}{4}$rad rotation of the ambient mesh, respectively. Also for these orientations, the obtained solutions are indistinguishable from the reference solution, despite the different cut-element shapes caused by the different rotations.

\subsection{Lattice of circular inclusions}
\label{sec:inclusion}
\begin{figure}
    \centering
    \includegraphics[width=0.99\textwidth]{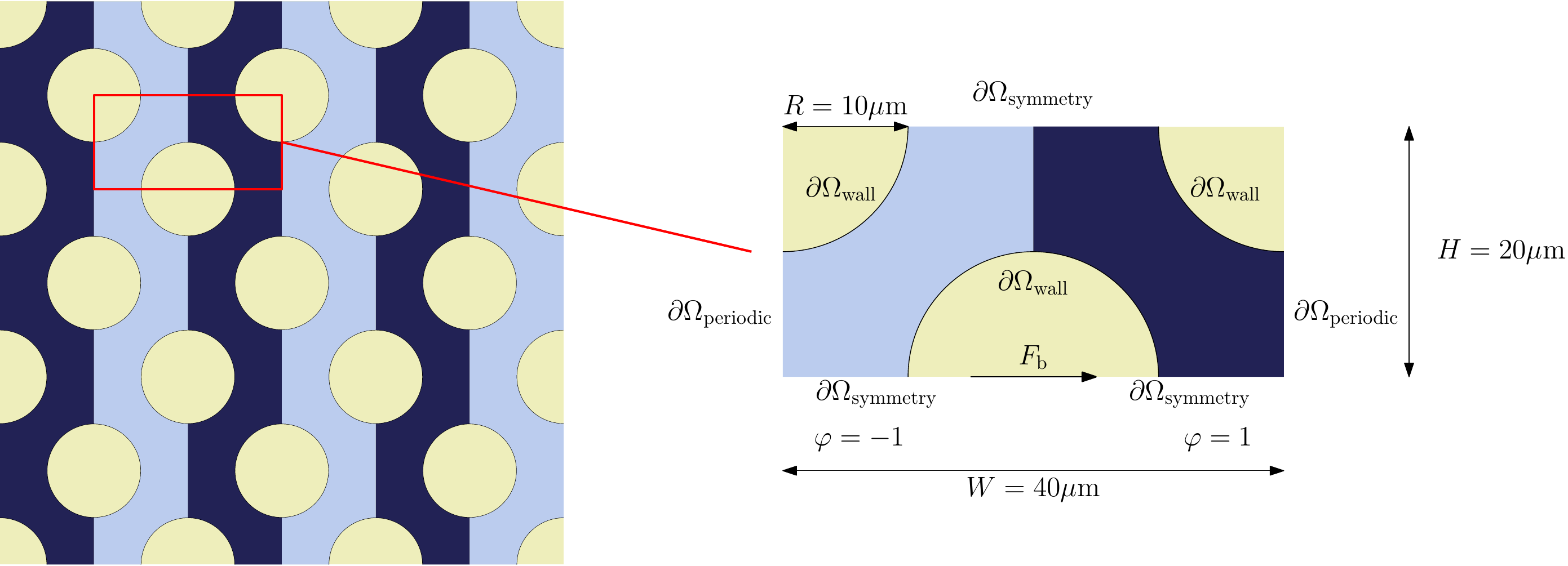}
    \caption{Illustration of the domain and boundary conditions for the binary-fluid flow through a lattice of circular inclusions.}
    \label{fig:latticedomain}
\end{figure}

\noindent
To demonstrate the capability of the immersed isogeometric analysis framework to handle arbitrary geometries, we consider a water-air binary-fluid system (\cref{tab:properties}) flowing through a lattice of circular inclusions (Fig.~\ref{fig:latticedomain}). The simulated unit cell has a width of $40{\rm \mu m}$ and a height of $20{\rm \mu m}$. The circular inclusions have a radius of $10{\rm \mu m}$ and an offset of $40{\rm \mu m}$ in both the horizontal and vertical direction. In the immersed isogeometric analysis framework, this problem requires a rectangular ambient domain corresponding to the bounding box of the unit cell. Construction of a boundary-fitted isogeometric analysis mesh is possible, but would require a multi-patch description.

The flow through the lattice is forced by means of a constant, horizontally oriented, body force $F_{\rm b}=1\cdot 10^9 \rm N/m^{3}$. In line with the lattice interpretation of this problem, we consider periodic boundary conditions for the problem on the left and right boundaries, meaning that the field variables attain identical values on corresponding points on the inflow and outflow boundaries. In our simulations, this periodic behavior is enforced by constructing an ambient domain mesh with periodic B-splines in the horizontal direction. Note that this construction is trivial on account of the rectangular ambient domain. On the top and bottom boundaries of the domain, symmetry boundary conditions are considered, meaning that the vertical velocity vanishes, as do the horizontal 
traction and the normal gradients of the phase field and chemical potential.

We discretize the ambient domain with elements of size $h=0.25{\rm \mu m}$. After trimming, this results in a discrete problem with 43,085 degrees of freedom. The initial time step size is set to $\Delta t=0.5{\rm \mu s}$. With interface velocities as high as approximately $70{\rm m/s}$, the use of time step adaptation is required. In the presented results this is implemented by dividing the time step by two when the nonlinear problem fails to converge, and restoring to the original time step after 8 converged time steps using the smaller step size; see also Ref.~\cite{Demont:2022dk}.
In the initial state we consider the two phases to be distributed as vertically oriented lamellae, with the interfaces positioned at the shortest distances between the obstacles, as illustrated in Fig.~\ref{fig:latticedomain}. The horizontal body force drives these lamellae past the obstacles, requiring them to break up and eventually to realign with the flow direction.

\begin{figure}
\begin{tabular}{cc}
\centering
\subfloat[$t=4.25\mu \rm s$]{\includegraphics[width=0.49\textwidth]{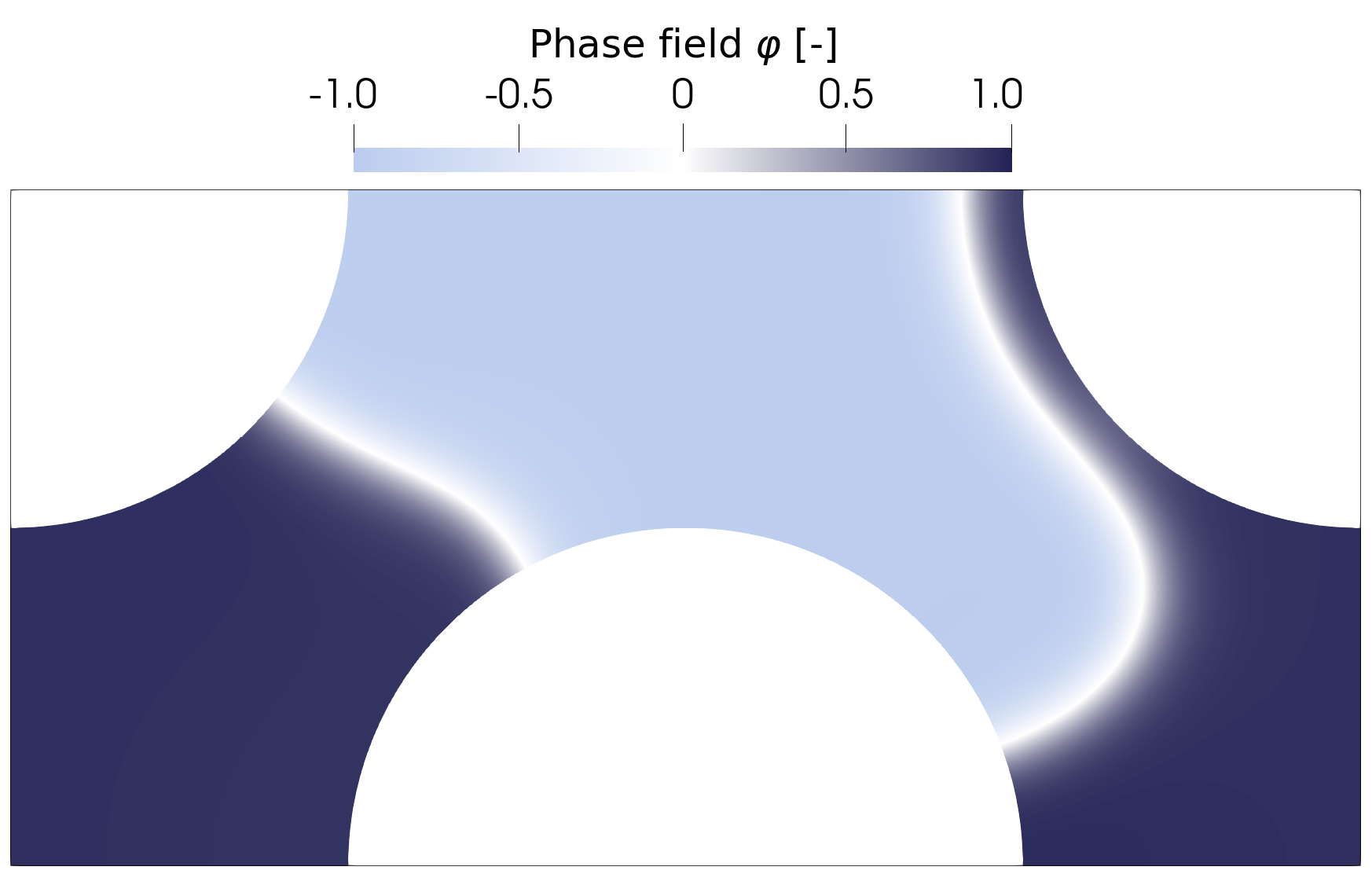}\label{fig:latticefirstcyclea}}
&
\subfloat[$t=4.25\mu \rm s$]{\includegraphics[width=0.49\textwidth]{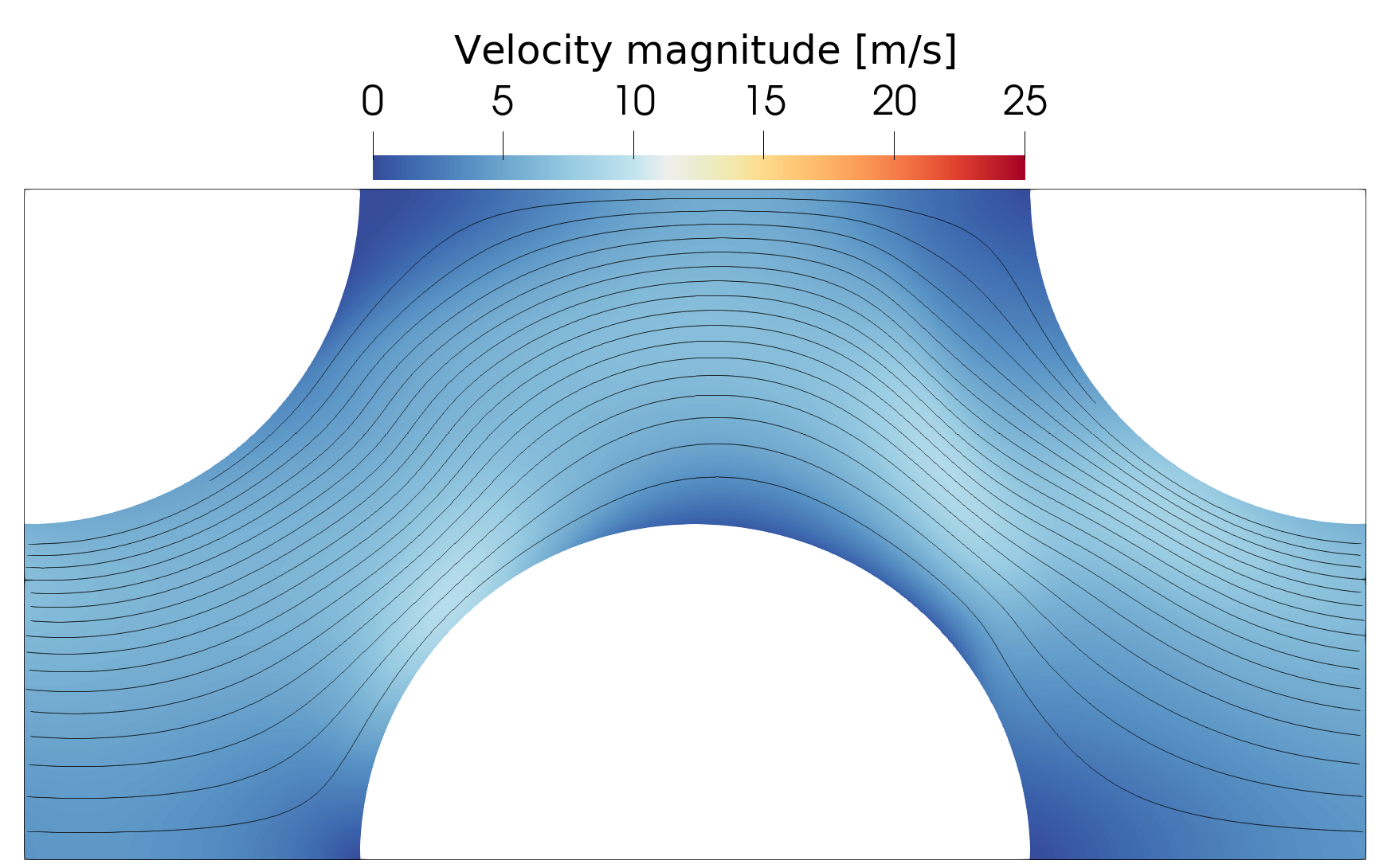}\label{fig:latticefirstcycleb}}
\\
\subfloat[$t=5.625\mu \rm s$]{\includegraphics[width=0.49\textwidth]{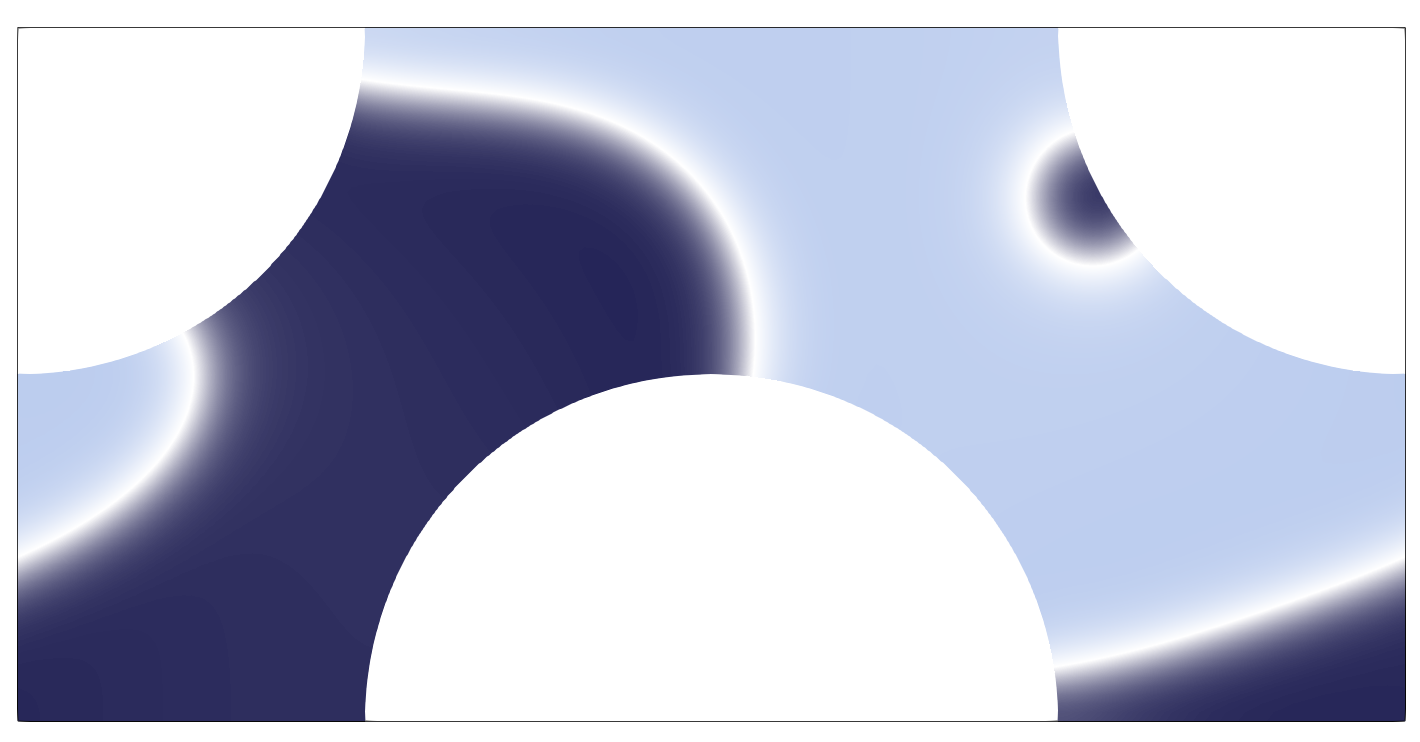}\label{fig:latticefirstcyclec}}
&
\subfloat[$t=5.625\mu \rm s$]{\includegraphics[width=0.49\textwidth]{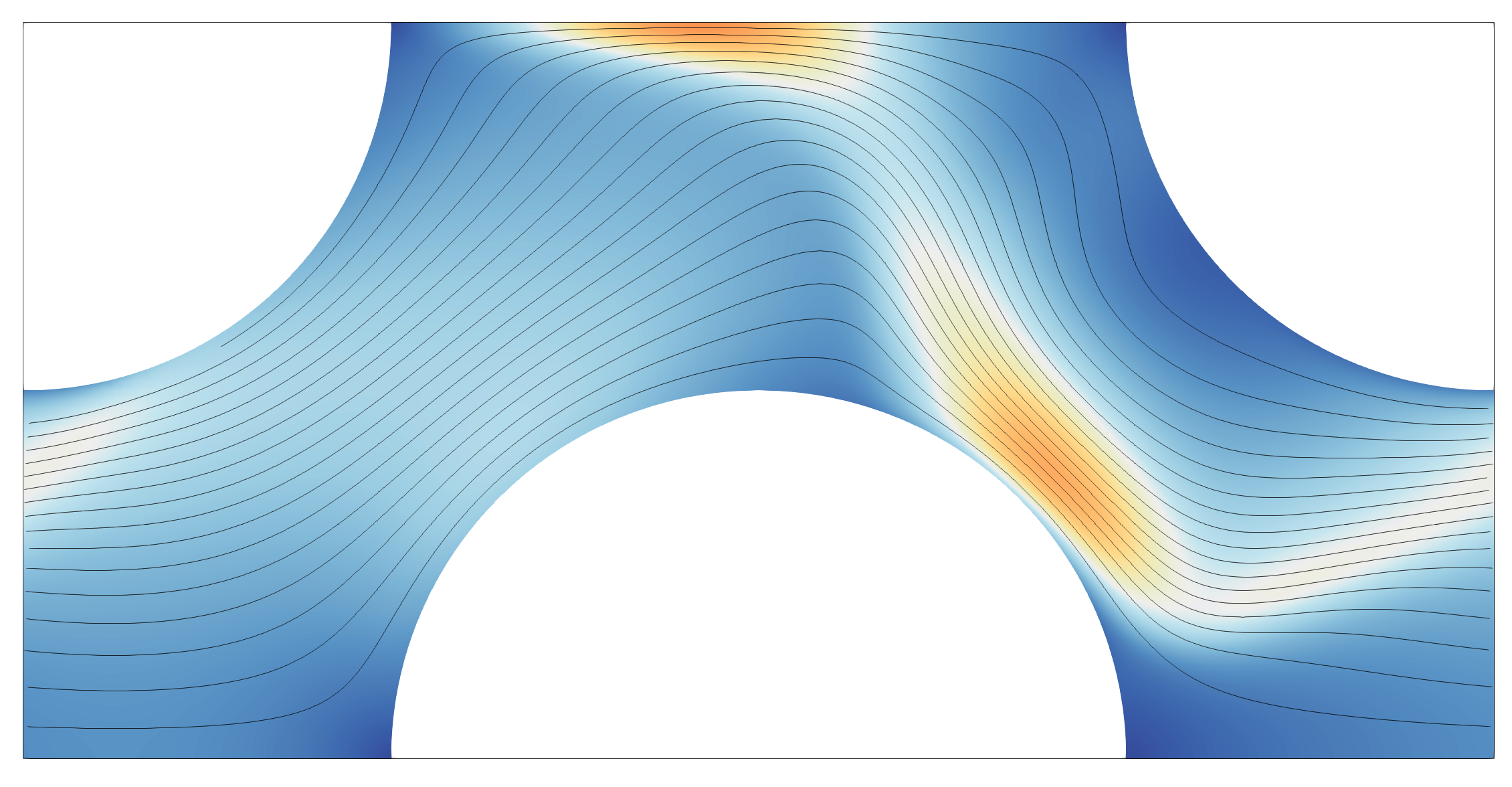}\label{fig:latticefirstcycled}}
\\
\subfloat[$t=6.25\mu \rm s$]{\includegraphics[width=0.49\textwidth]{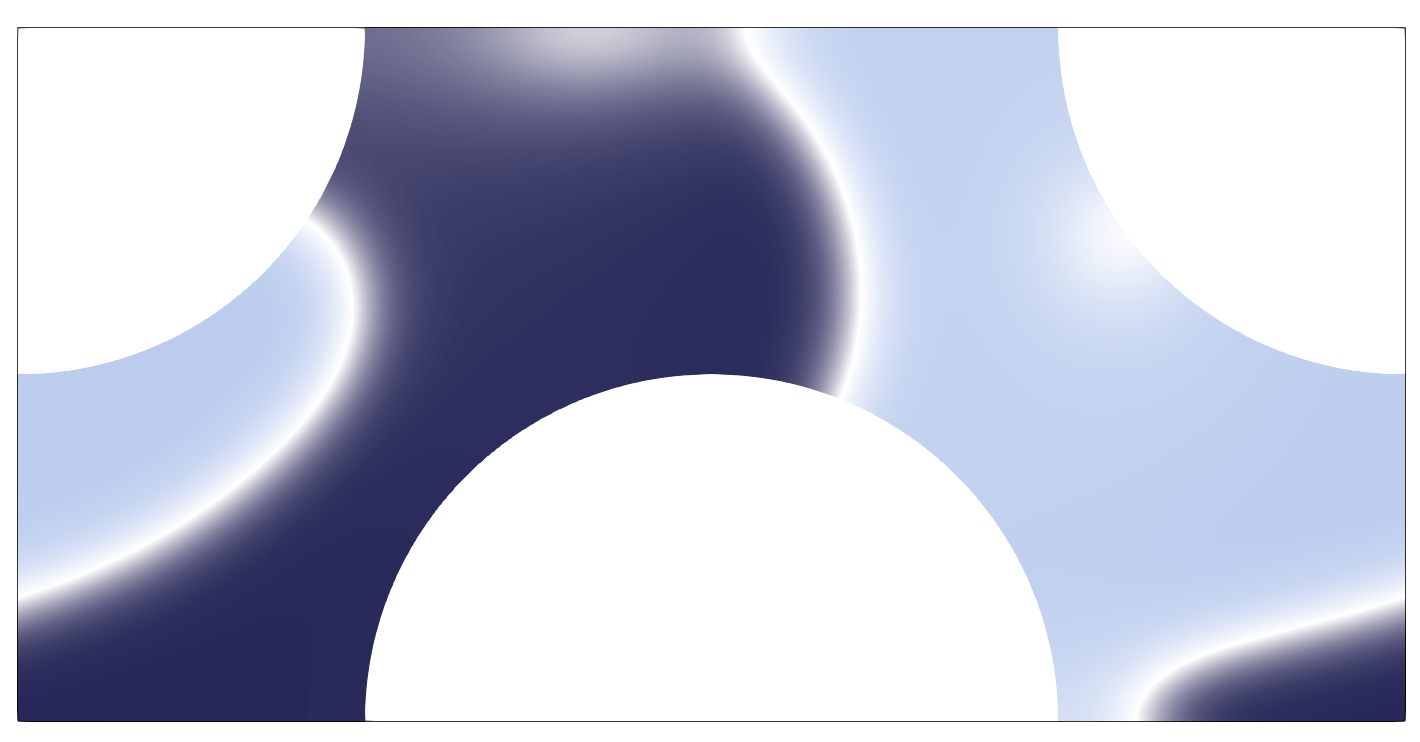}\label{fig:latticefirstcyclee}}
&
\subfloat[$t=6.25\mu \rm s$]{\includegraphics[width=0.49\textwidth]{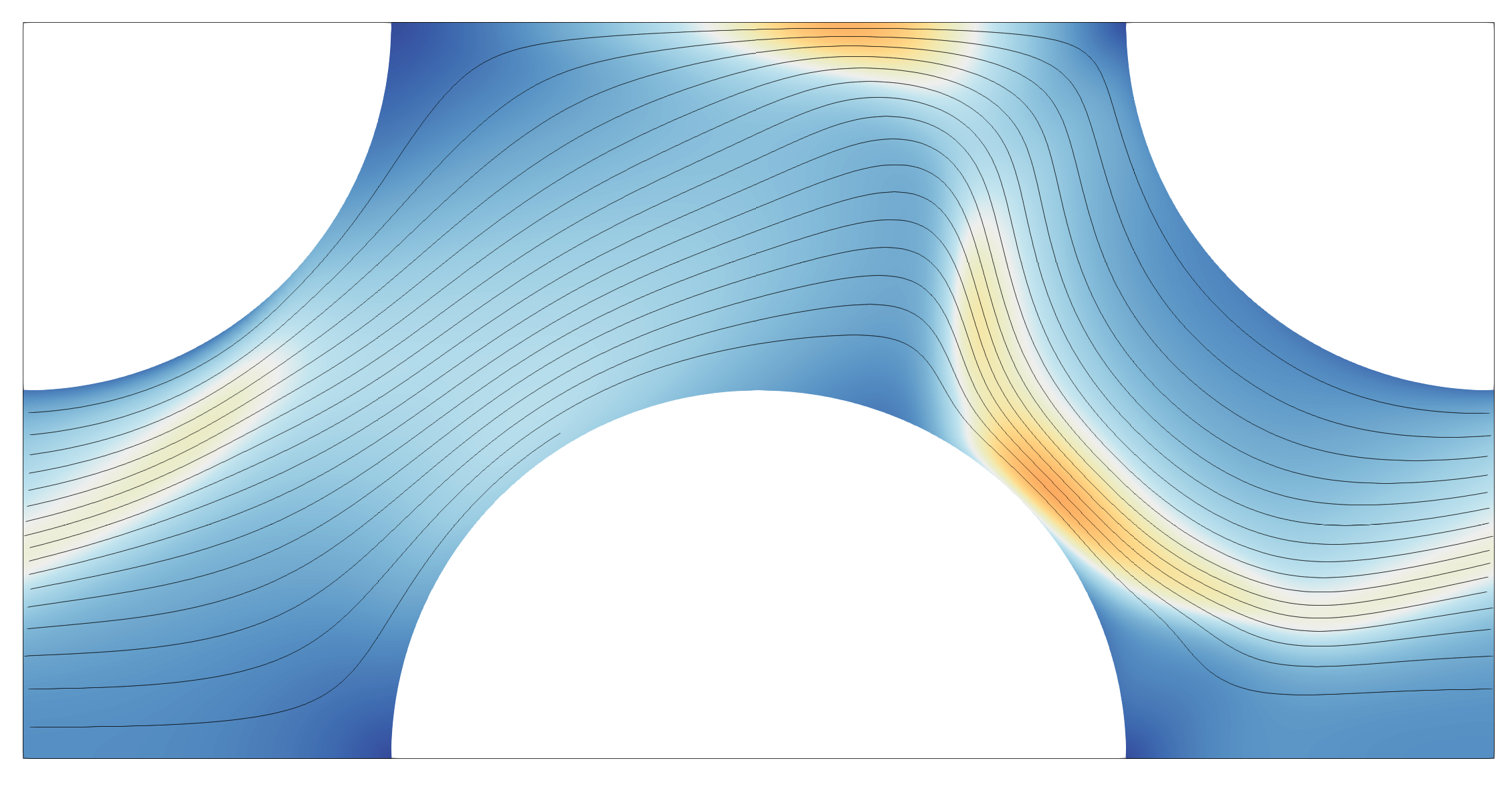}\label{fig:latticefirstcyclef}}
\\
\end{tabular}
    \caption{Time evolution of the phase field $\varphi$ (left) and velocity field $\u$ (right), showing the motion of the lamellae past the obstacles.}
\end{figure}

\begin{figure}\ContinuedFloat
\begin{tabular}{cc}
\centering
\subfloat[$t=13.4375\mu \rm s$]{\includegraphics[width=0.49\textwidth]{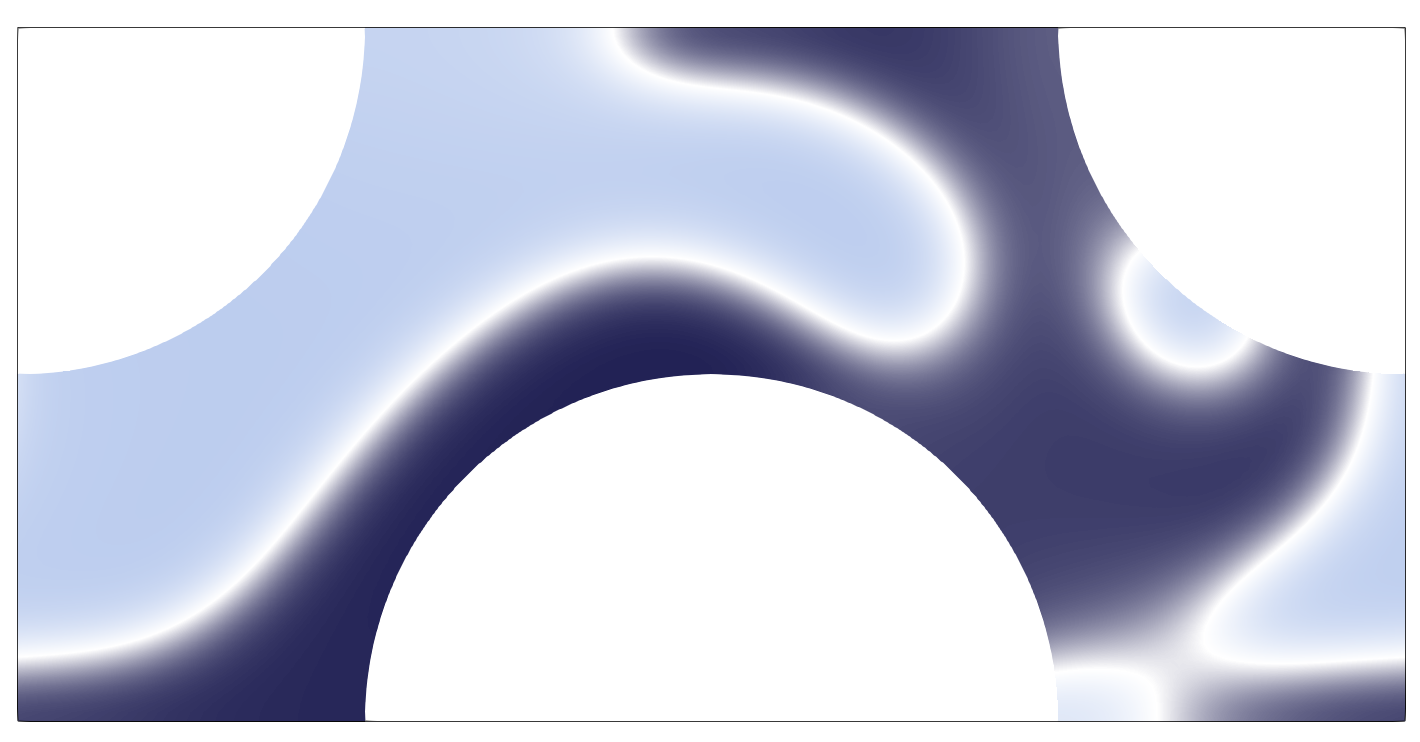}\label{fig:latticefirstcycleg}}
&
\subfloat[$t=13.4375\mu \rm s$]{\includegraphics[width=0.49\textwidth]{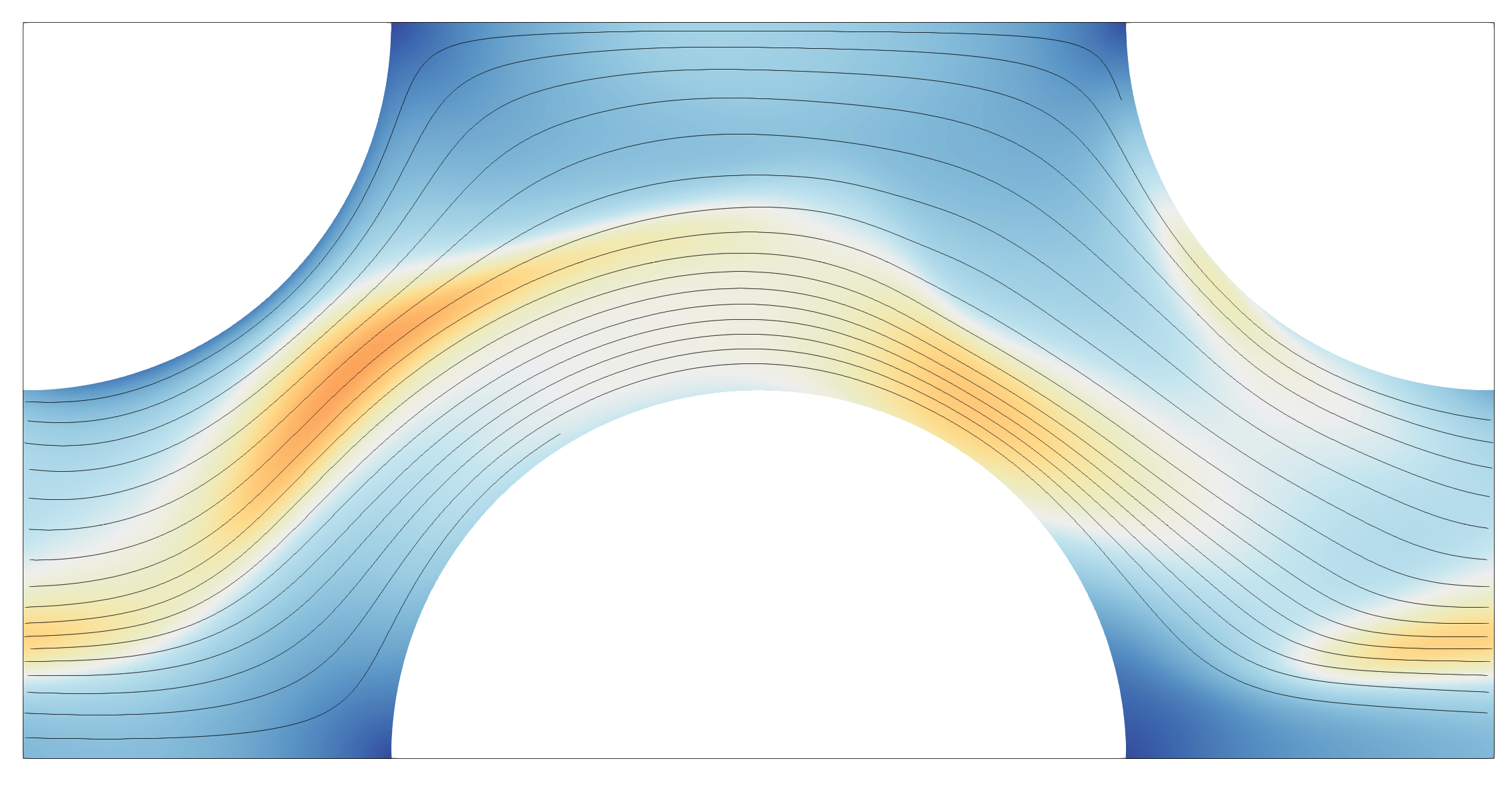}\label{fig:latticefirstcycleh}}
\\
\subfloat[$t=15.75\mu\rm s$]{\includegraphics[width=0.49\textwidth]{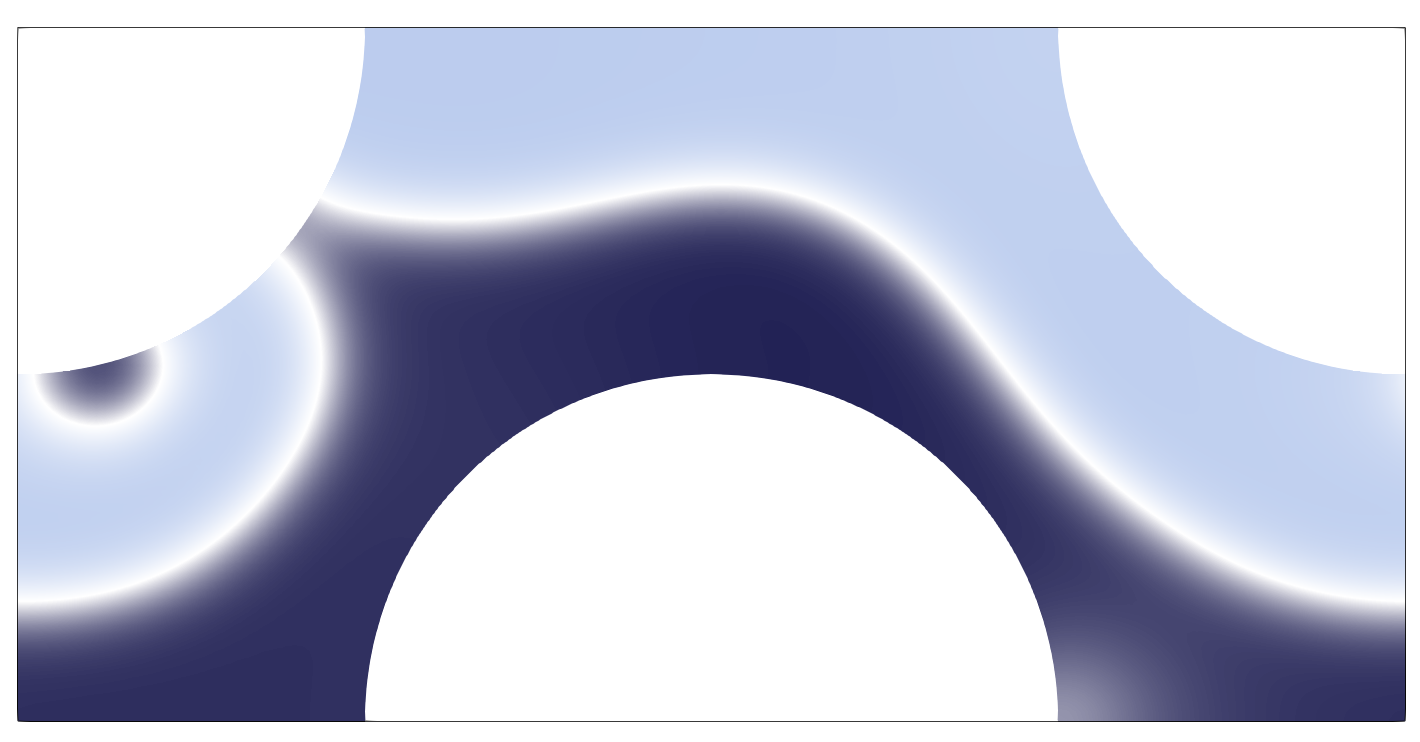}\label{fig:latticefirstcyclei}}
&
\subfloat[$t=15.75\mu\rm s$]{\includegraphics[width=0.49\textwidth]{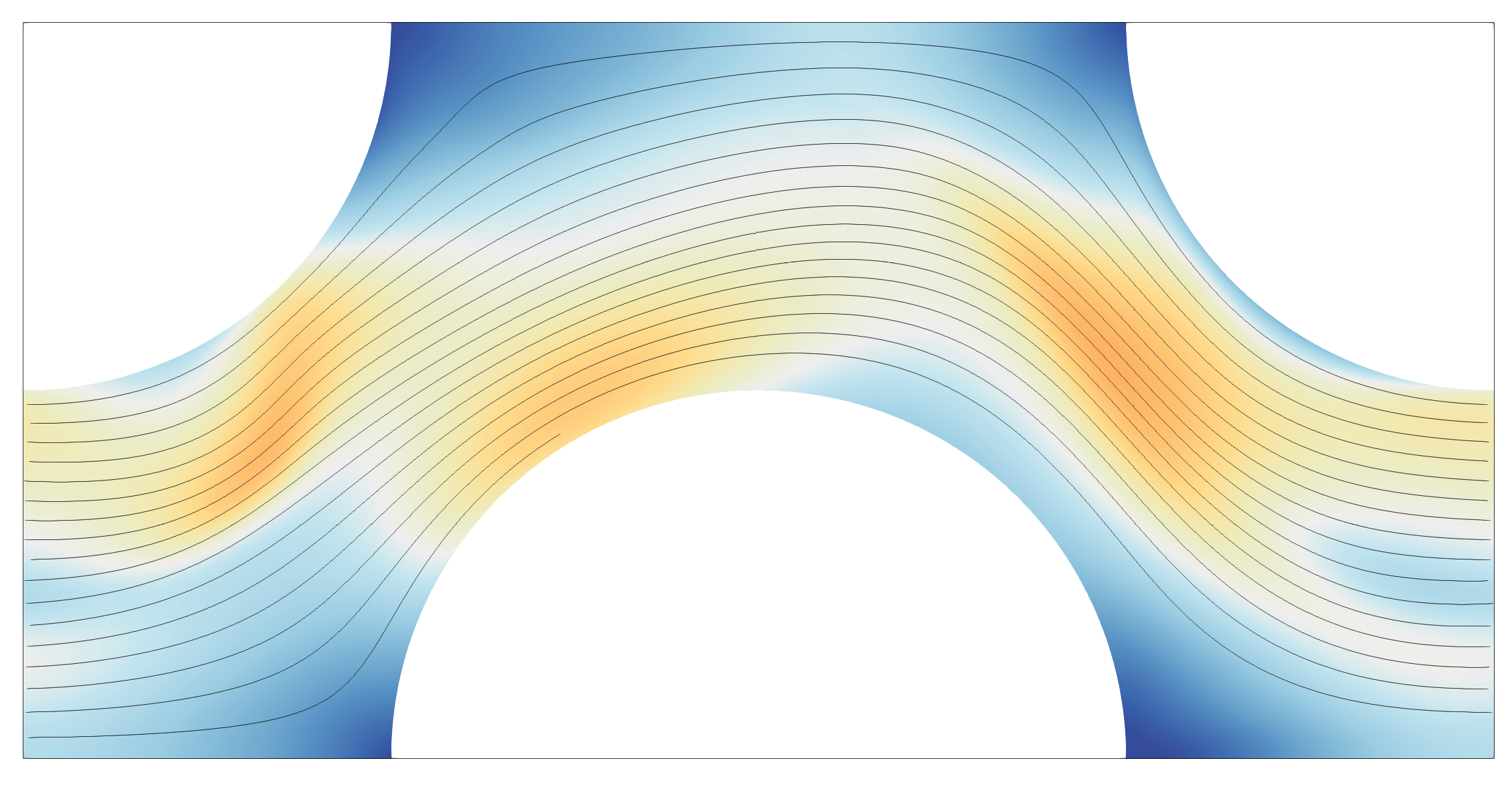}\label{fig:latticefirstcyclej}}
\\
\subfloat[$t=16.875\mu\rm s$]{\includegraphics[width=0.49\textwidth]{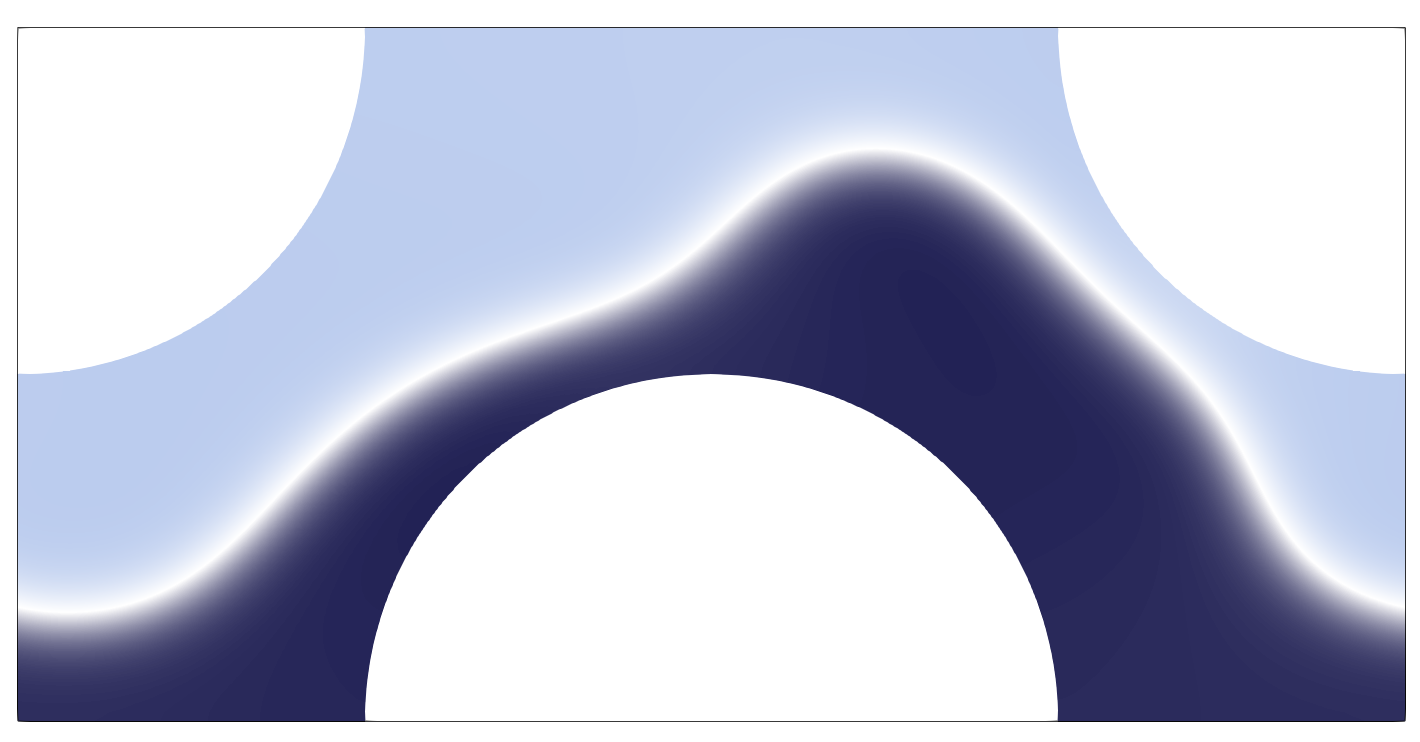}\label{fig:latticefirstcyclek}}
&
\subfloat[$t=16.875\mu\rm s$]{\includegraphics[width=0.49\textwidth]{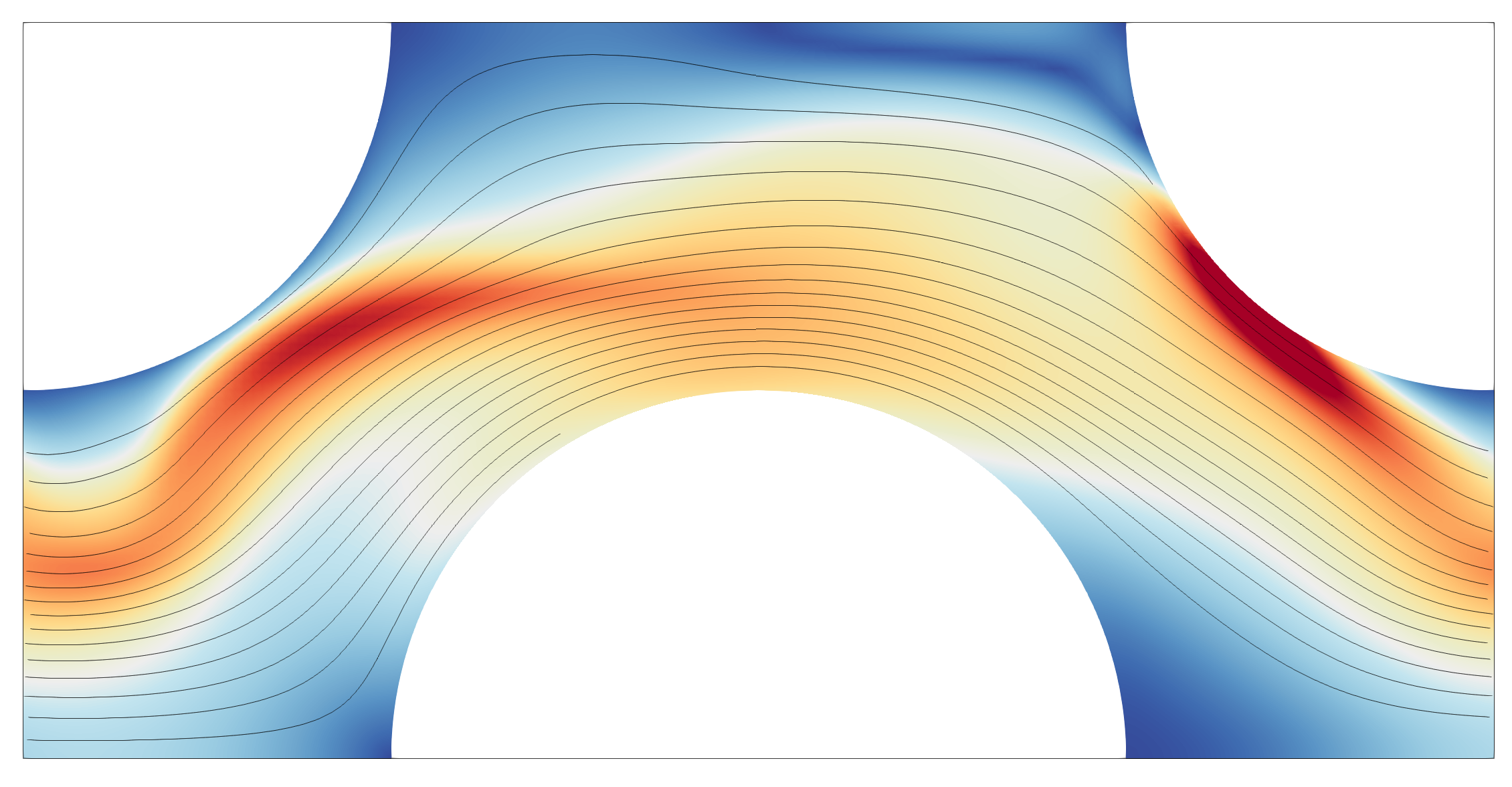}\label{fig:latticefirstcyclel}}
\\
\end{tabular}
    \caption{Time evolution of the phase field $\varphi$ (left) and velocity field $\u$ (right), showing the motion of the lamellae past the obstacles.}
    \label{fig:latticefirstcycle}
\end{figure}

\begin{figure}
\begin{tabular}{cc}
\centering
\subfloat[$t=\infty$]{\includegraphics[width=0.49\textwidth]{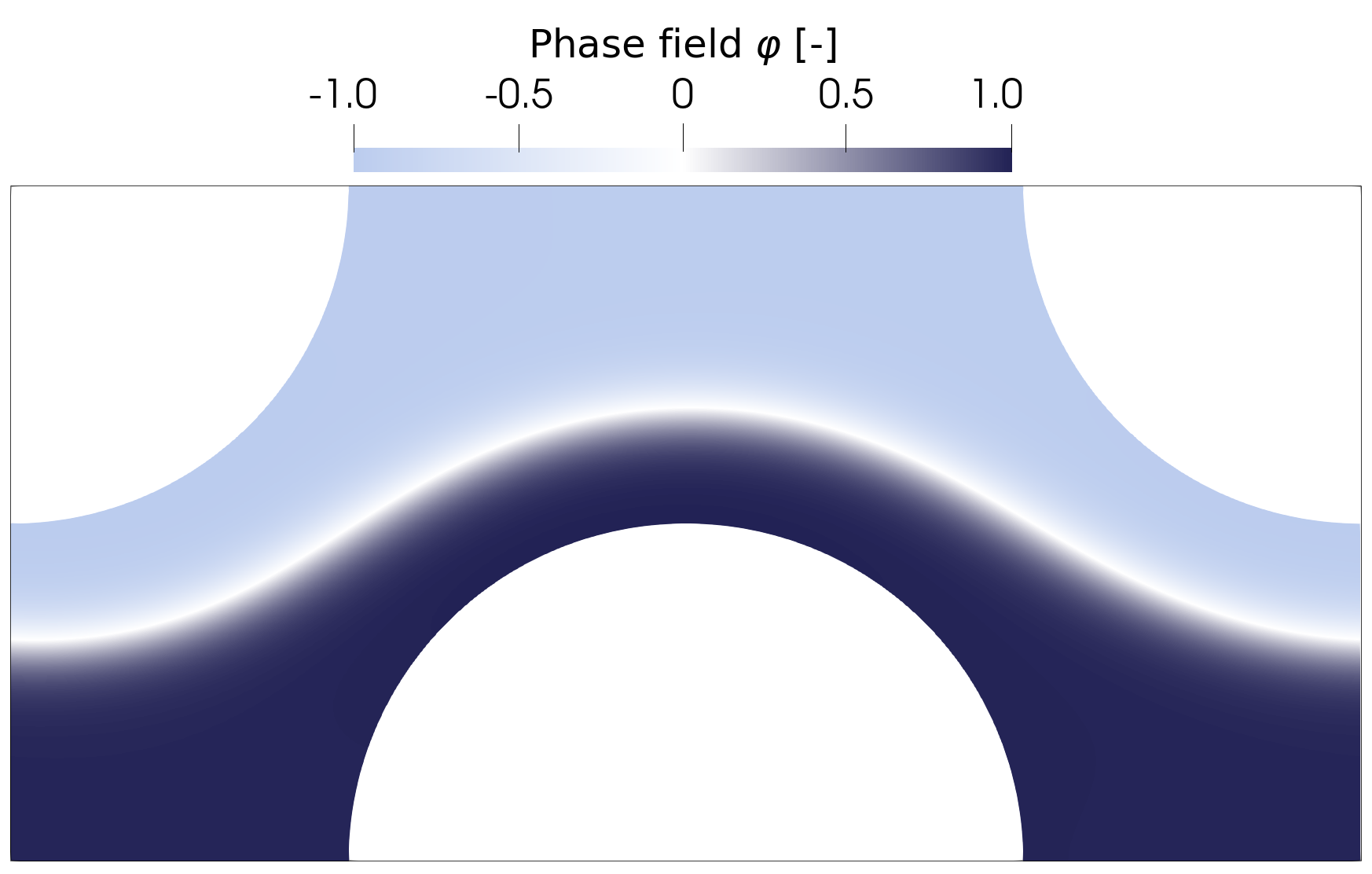}\label{fig:latticefirstcyclet25.5mus}}
&
\subfloat[$t=\infty$]{\includegraphics[width=0.49\textwidth]{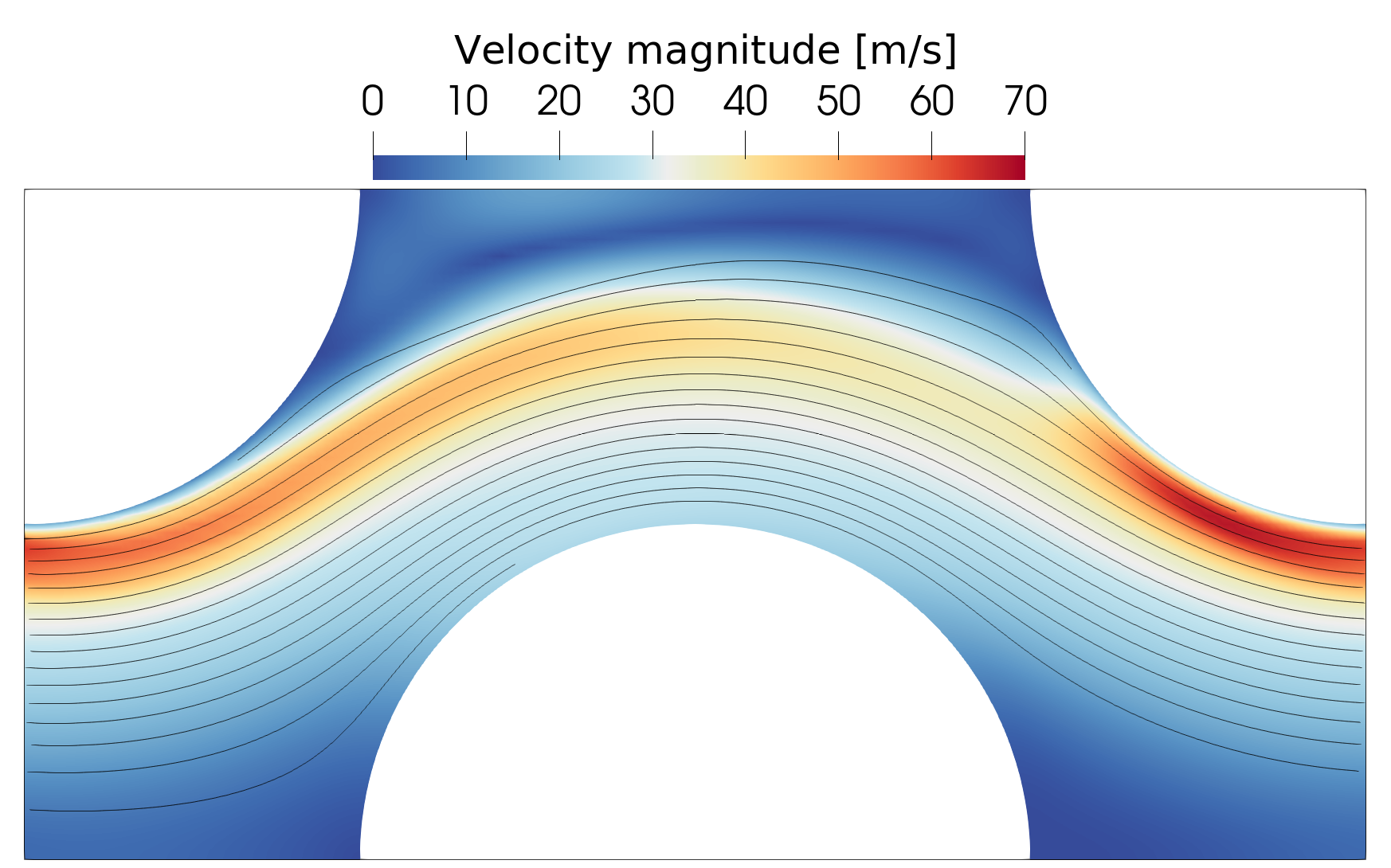}}
\\
\end{tabular}
    \caption{The phase field $\varphi$ (left) and velocity field $\u$ (right), showing the steady state of the flow in a porous medium with circular inclusions.}
    \label{fig:latticesteady}
\end{figure}

Fig.~\ref{fig:latticefirstcycle} illustrates the transient behavior of the lamellae passing through the unit cell. 
Initially, the water lamella (corresponding to $\varphi = 1$, shown in dark blue) enters the domain from the left. As the flow velocity is stalled by the obstacles, the fluids move fastest near the center of the free space between the inclusions. This causes the interfaces to curve, as observed already in Fig.~\ref{fig:latticefirstcyclea}. In this figure, the air-water interface approaches the upper-right obstacle. Due to the low velocities close to the obstacle boundary, a thin layer of water remains, eventually resulting in the trapped water droplet shown in Fig.~\ref{fig:latticefirstcyclec}. Due to the Ostwald ripening effect, this satellite droplet is quickly re-absorbed by the inflowing water lamella, see Fig.~\ref{fig:latticefirstcyclee}. Fig.~\ref{fig:latticefirstcycleg} shows the phase field after the lamellae have moved through two complete cycles of the periodic unit-cell. At this time instance, the water lamella coalesces with the detaching tail of the upstream lamella. After the isolation and re-absorption of a few more droplets, which may be observed in \cref{fig:latticefirstcycleg,fig:latticefirstcyclei}, the fully re-segregated and realigned state shown in Fig.~\ref{fig:latticefirstcyclek} is realized after approximately 5 cycles. The steady state, illustrated in Fig.~\ref{fig:latticesteady}, is reached shortly after.
In their new horizontal orientation, no further topological changes to the lamellae are induced by the obstacles.

\subsection{Porous medium}
\label{sec:porousmedium}
\begin{figure}[!b]
    \centering
    \includegraphics[width=0.95\textwidth]{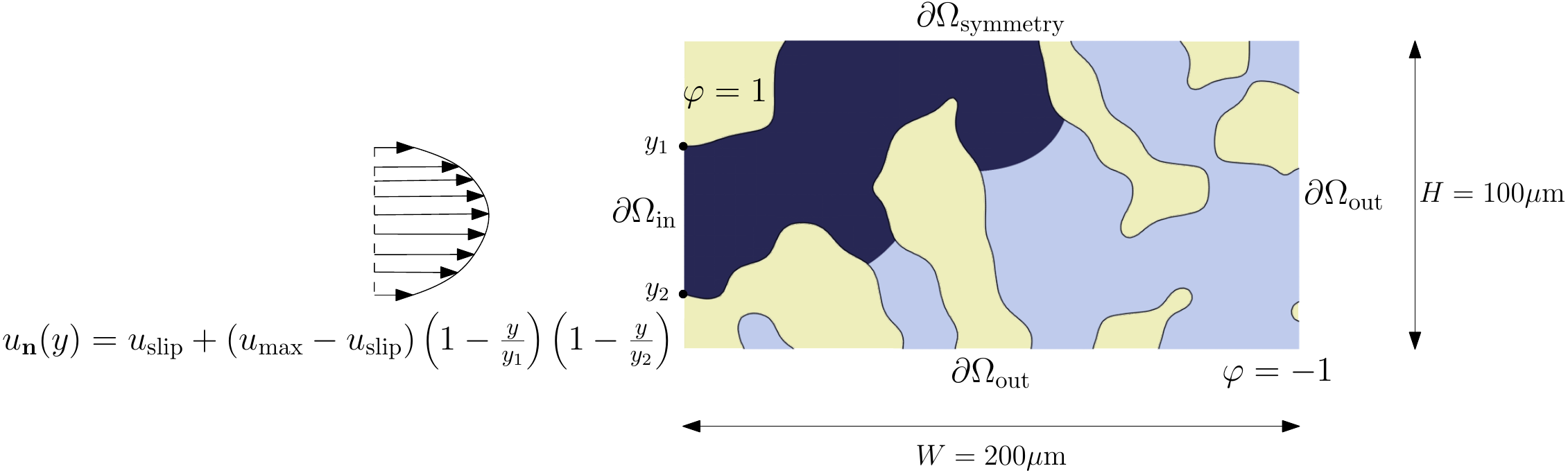}
    \caption{Illustration of the domains and inflow conditions for the binary flow through an example porous medium.\\}
    \label{fig:porousdomain}
\end{figure}

\noindent
To illustrate the capabilities of the immersed isogeometric analysis framework, we study the evolution of a water-air binary-fluid system propagating through the example porous medium depicted in Fig.~\ref{fig:porousdomain}. The width and height of the domain are $200\mu\rm m$ and $50\mu\rm m$, respectively. On the left boundary of the domain, the velocity is prescribed by a horizontal flow with a parabolic profile which satisfies the generalized Navier condition. The maximum velocity of the inflow profile is set to $5 {\rm m/s}$. On the inflow boundary, the phase field is constrained to $\varphi=1$, representing the (dark blue) water phase. On the top boundary, symmetry conditions are imposed, similar to the lattice problem discussed above. Outflow boundary conditions are considered on the bottom and right boundaries, meaning that the traction is set to zero. We make use of a mesh size of $h=0.625\mu\rm m$, leading to 185,618 degrees-of-freedom systems of equations, and a time step size of $\Delta t=0.05{\rm \mu s}$. 

\begin{figure}
    \begin{tabular}{cc}
\centering
\subfloat[$t=0\mu \rm s$]{\includegraphics[width=0.49\textwidth]{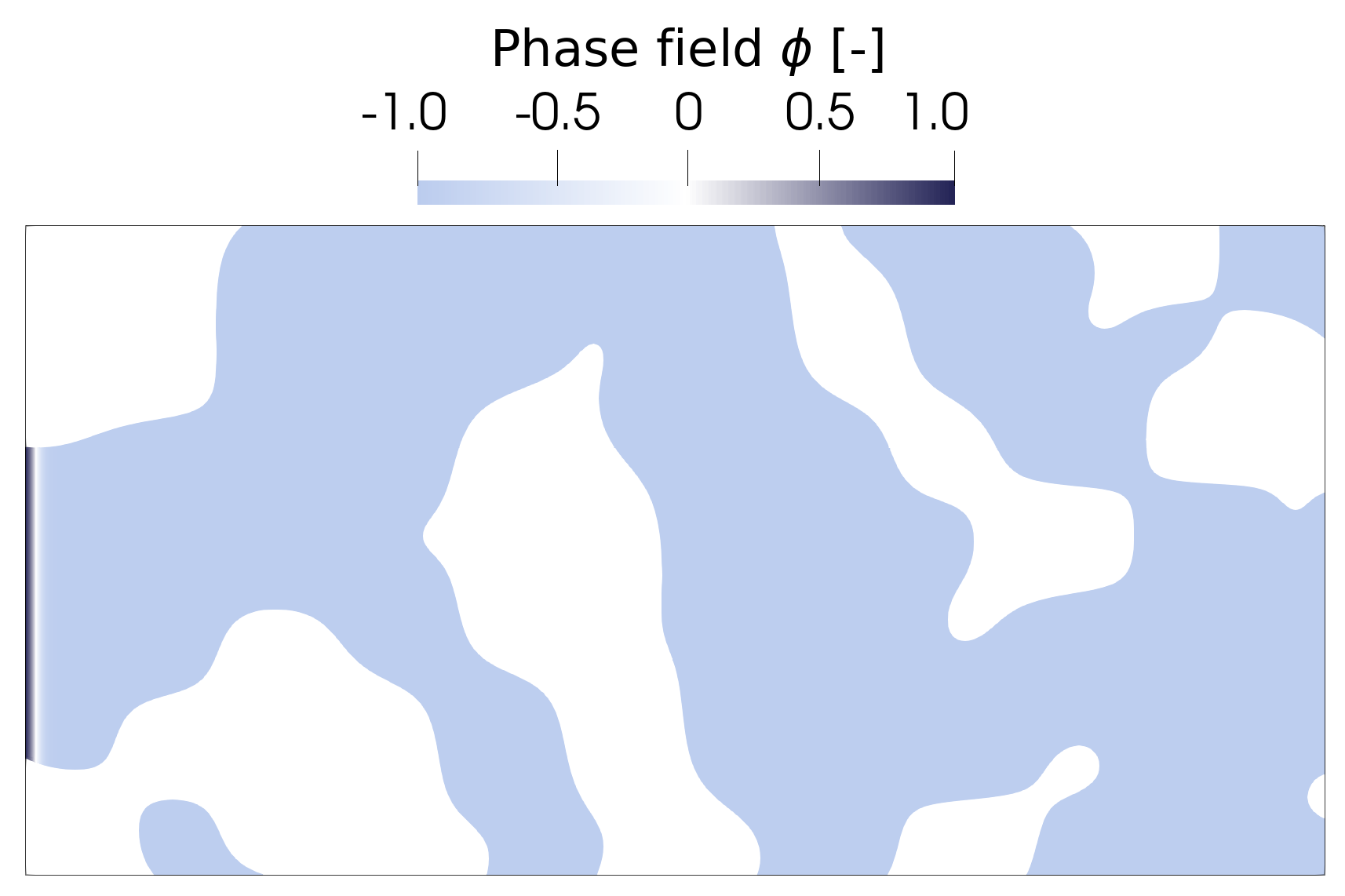}\label{fig:poroussnapshotsa}}
&
\subfloat[$t=0\mu \rm s$]{\includegraphics[width=0.49\textwidth]{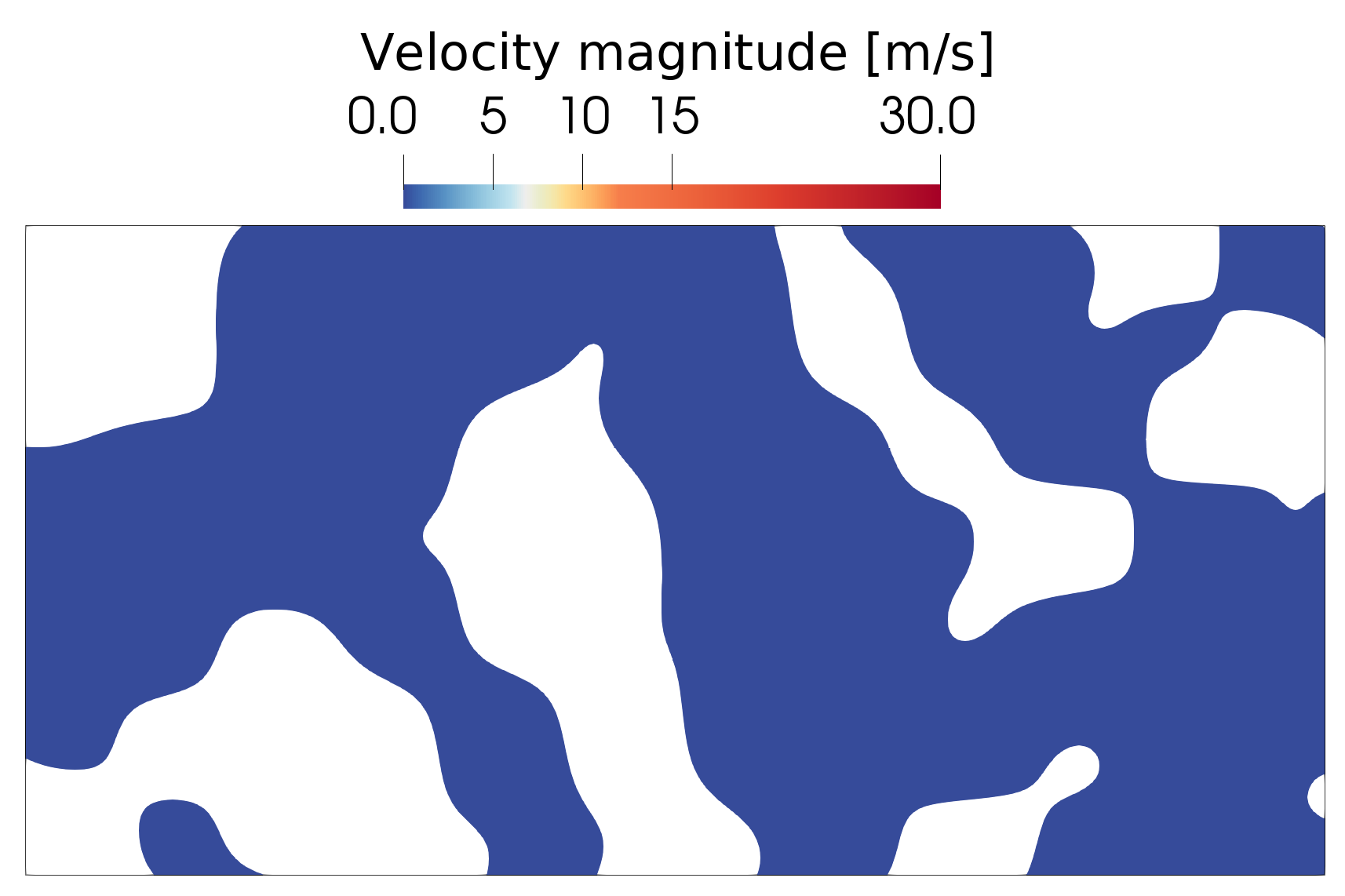}}
\\
\subfloat[$t=9\mu \rm s$]{\includegraphics[width=0.49\textwidth]{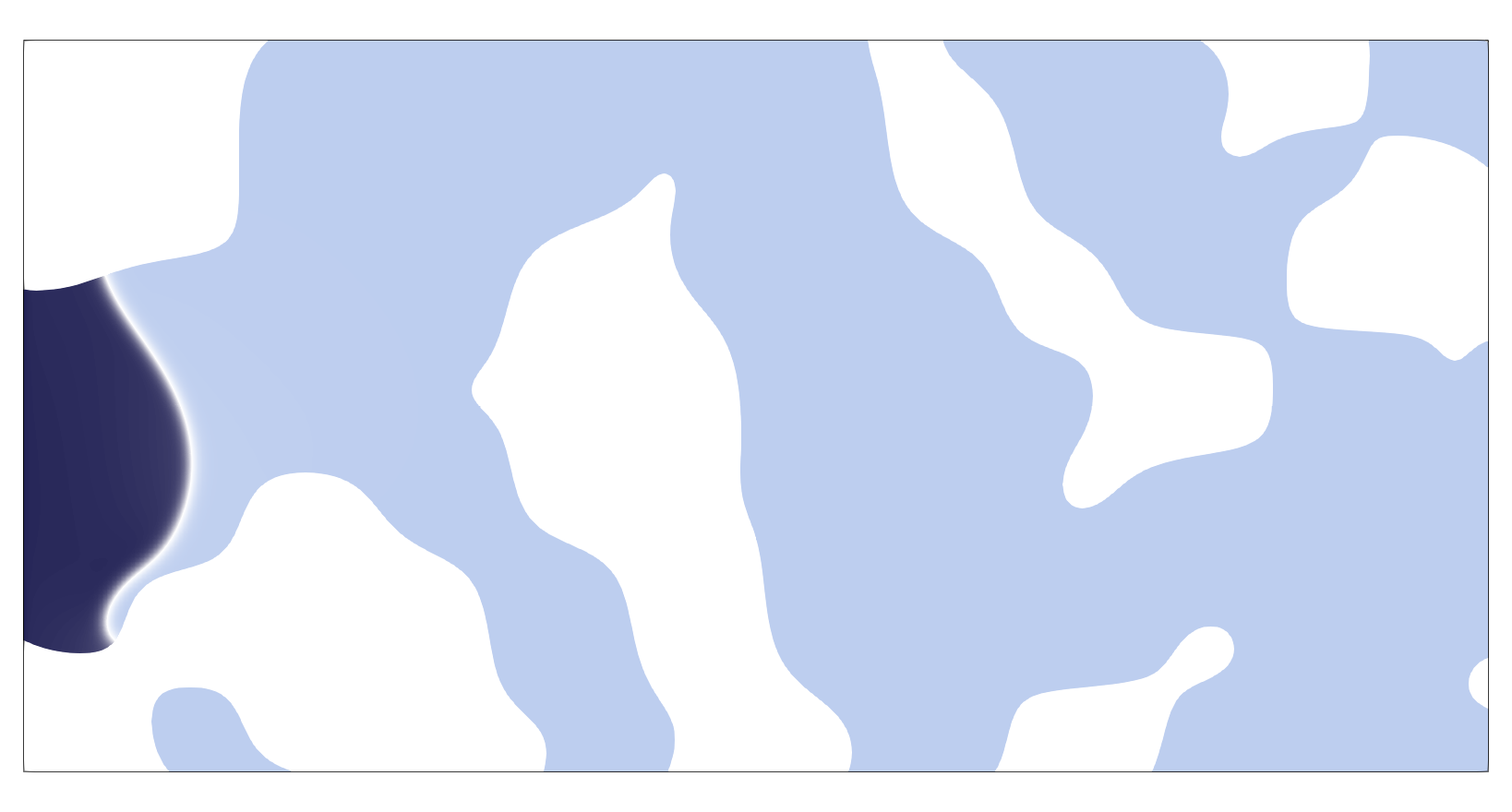}\label{fig:poroussnapshotsc}}
&
\subfloat[$t=9\mu \rm s$]{\includegraphics[width=0.49\textwidth]{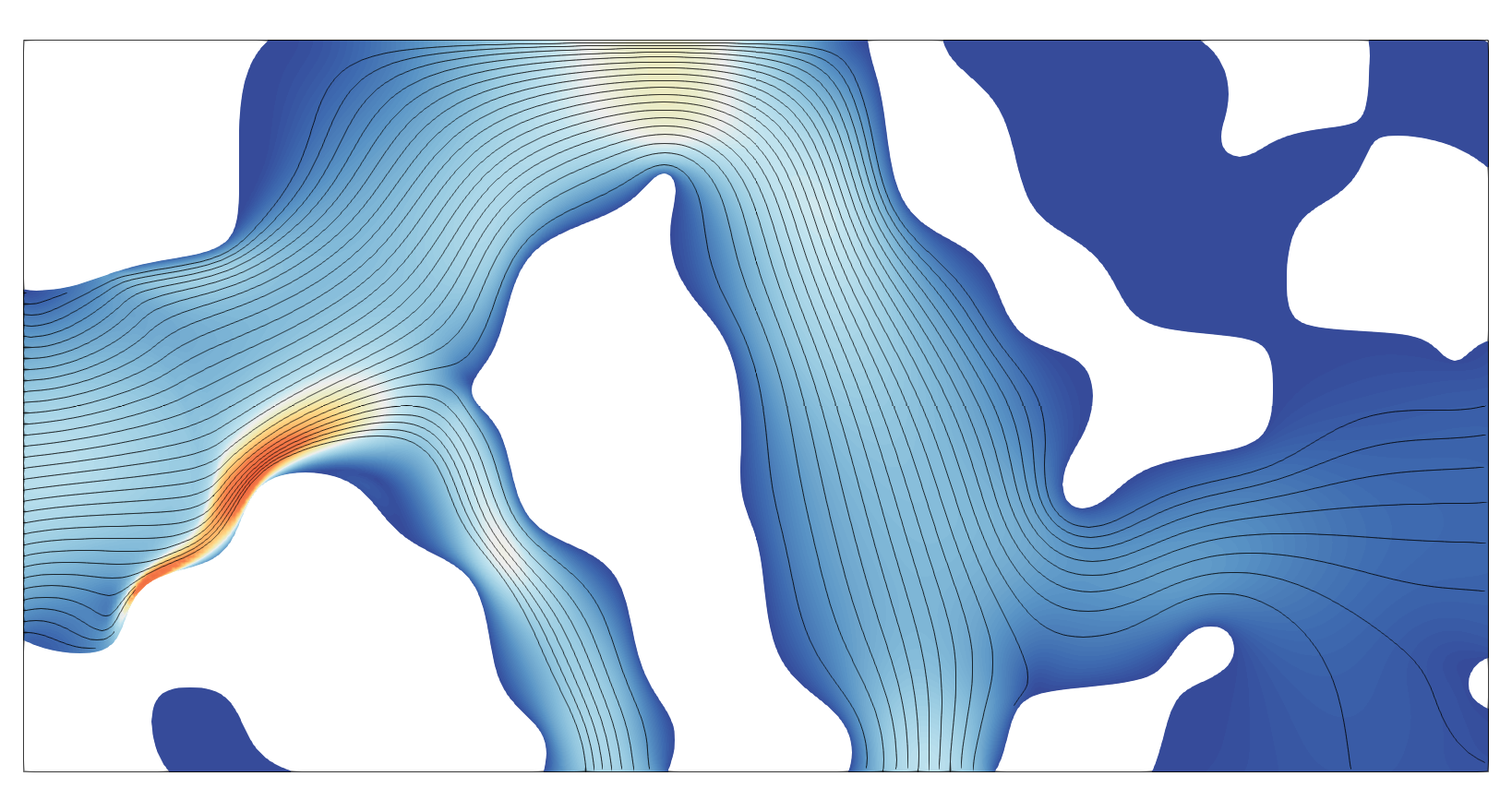}}
\\
\subfloat[$t=10\mu \rm s$]{\includegraphics[width=0.49\textwidth]{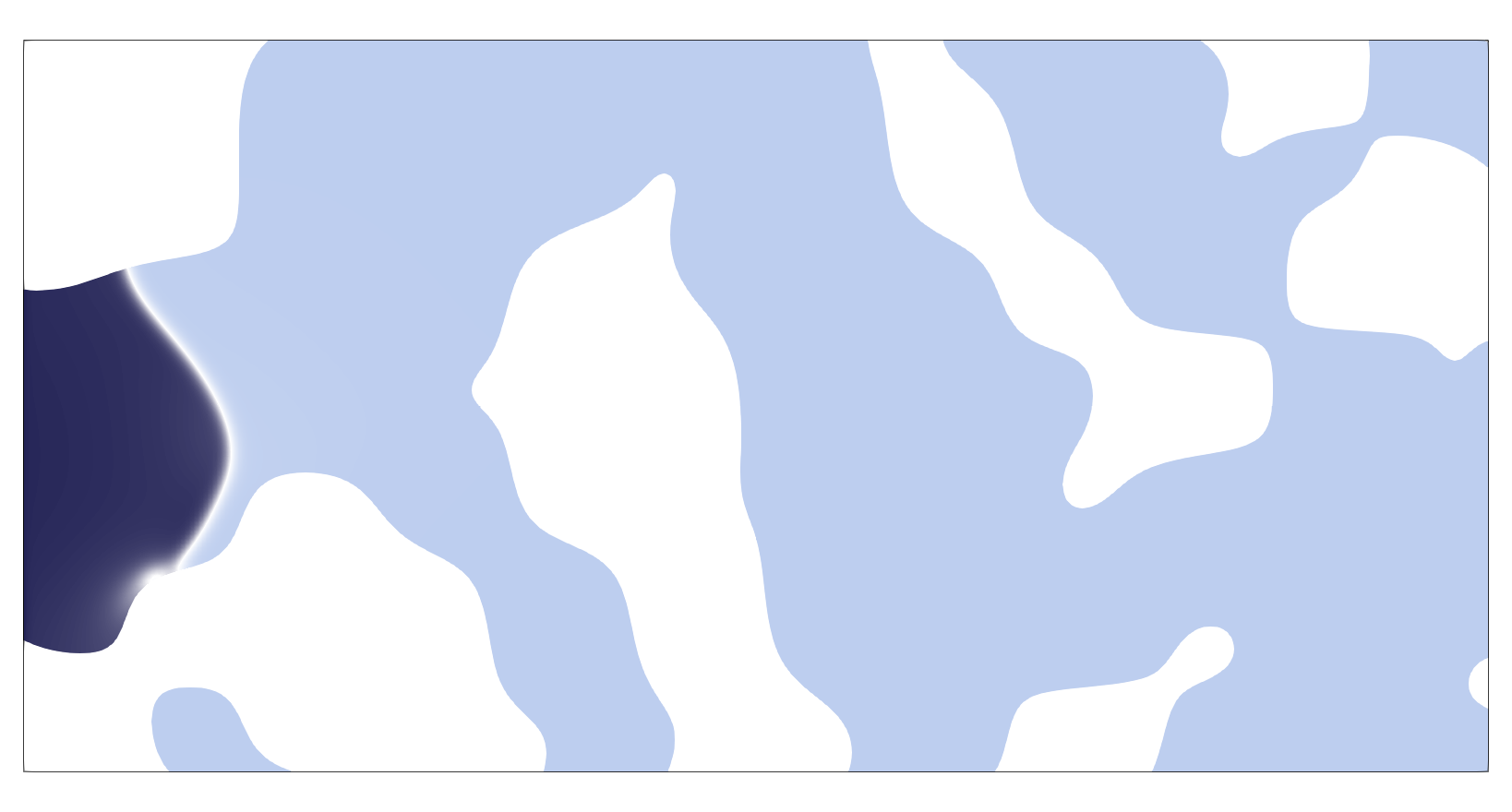}\label{fig:poroussnapshotse}}
&
\subfloat[$t=10\mu \rm s$]{\includegraphics[width=0.49\textwidth]{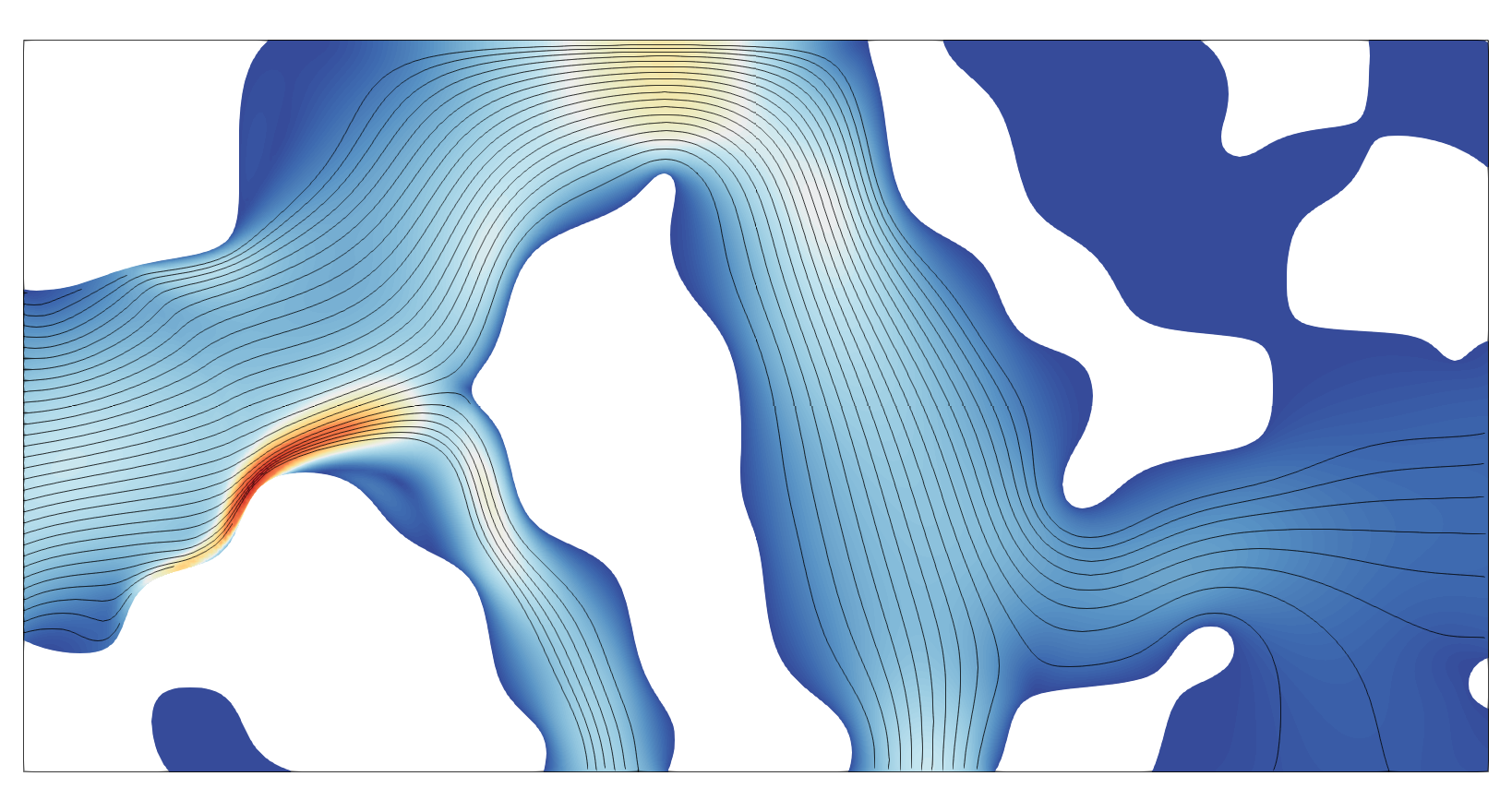}}
\\
\end{tabular}
    \caption{Time evolution of the phase field $\varphi$ (left) and velocity field $\u$ (right) for the porous medium test case.}
\end{figure}

\begin{figure}\ContinuedFloat
\begin{tabular}{cc}
\centering
\subfloat[$t=17.0\mu \rm s$]{\includegraphics[width=0.49\textwidth]{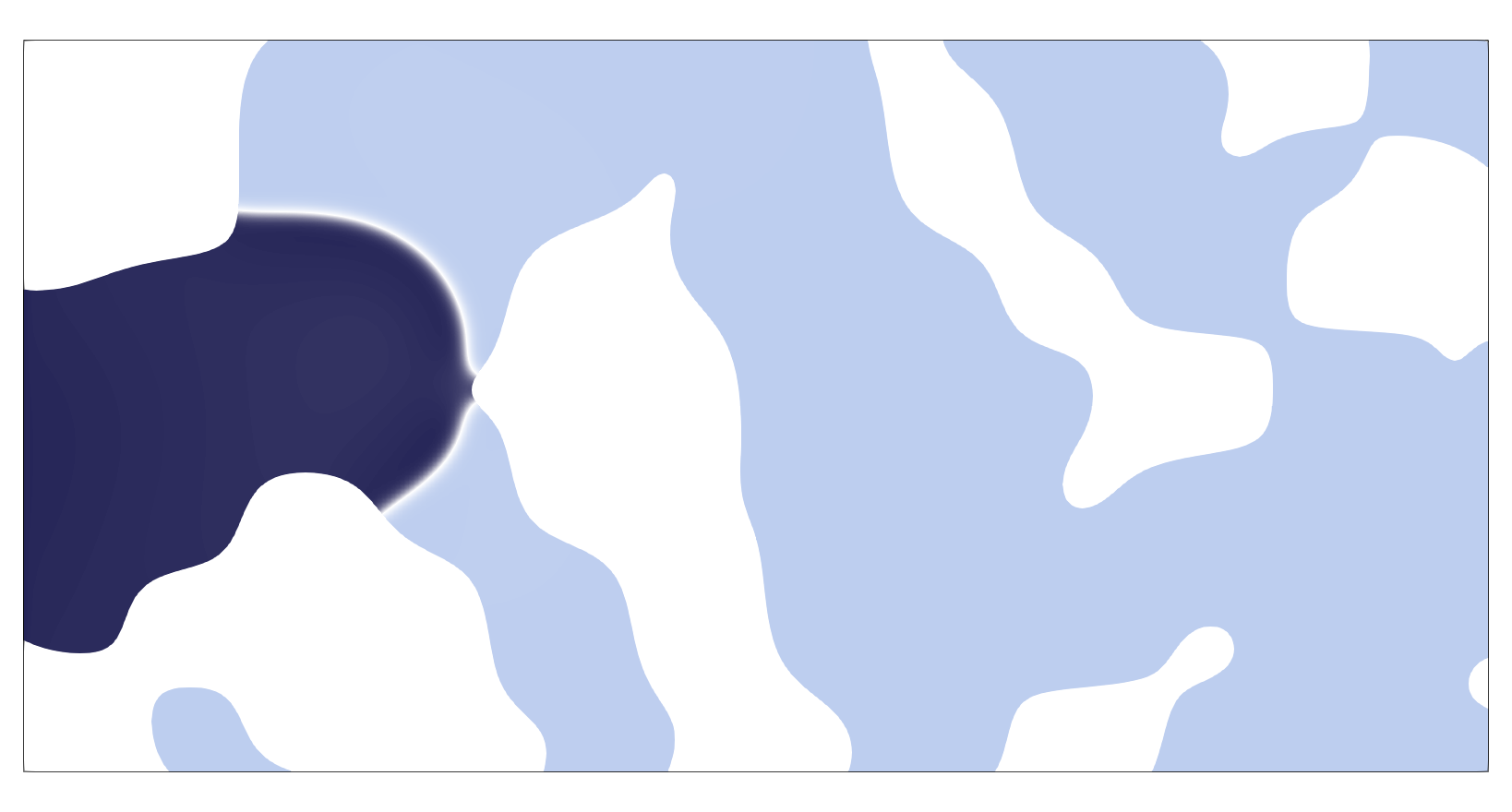}\label{fig:poroussnapshotsg}}
&
\subfloat[$t=17.0\mu \rm s$]{\includegraphics[width=0.49\textwidth]{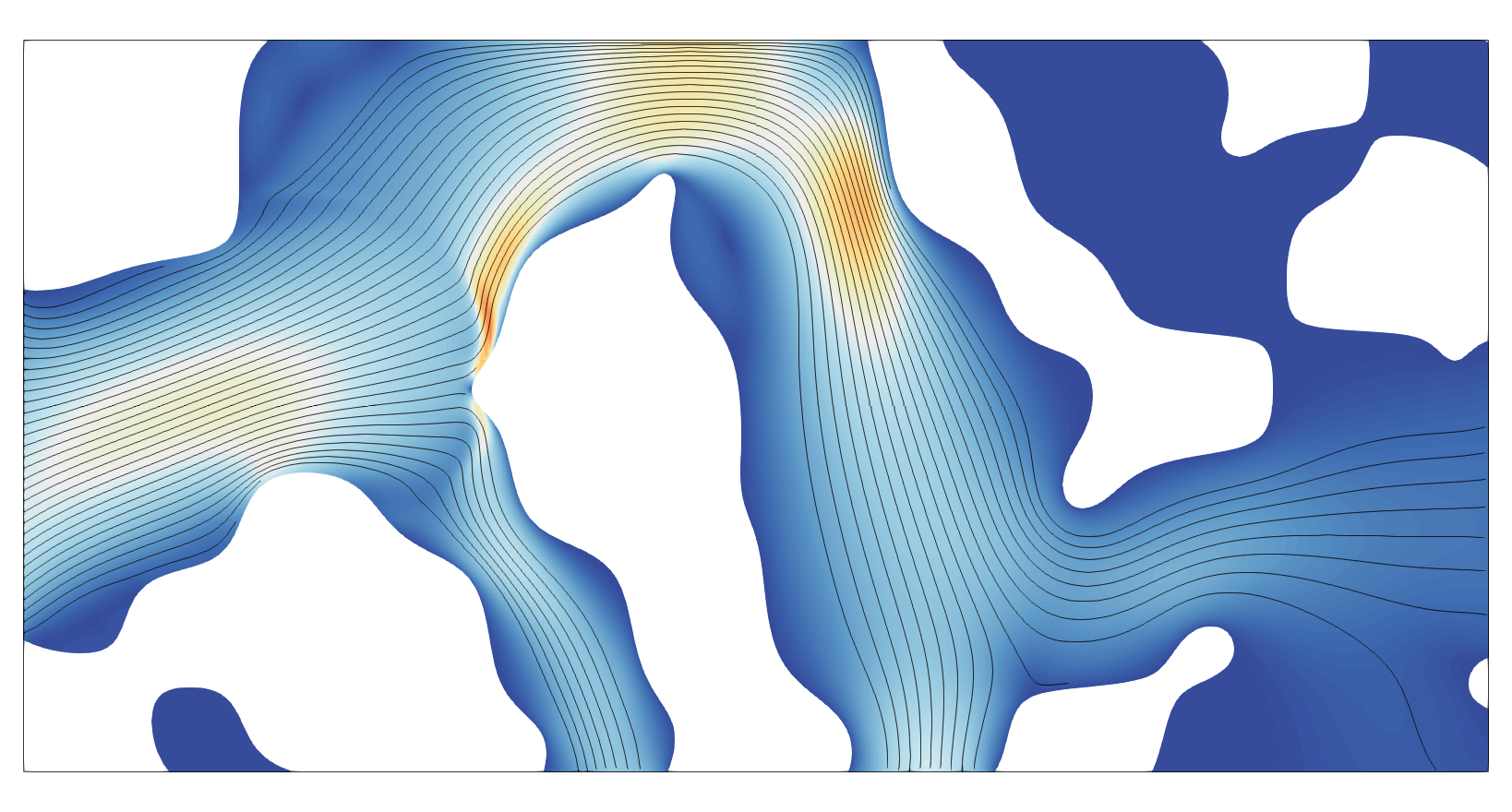}}
\\
\subfloat[$t=25.9\mu \rm s$]{\includegraphics[width=0.49\textwidth]{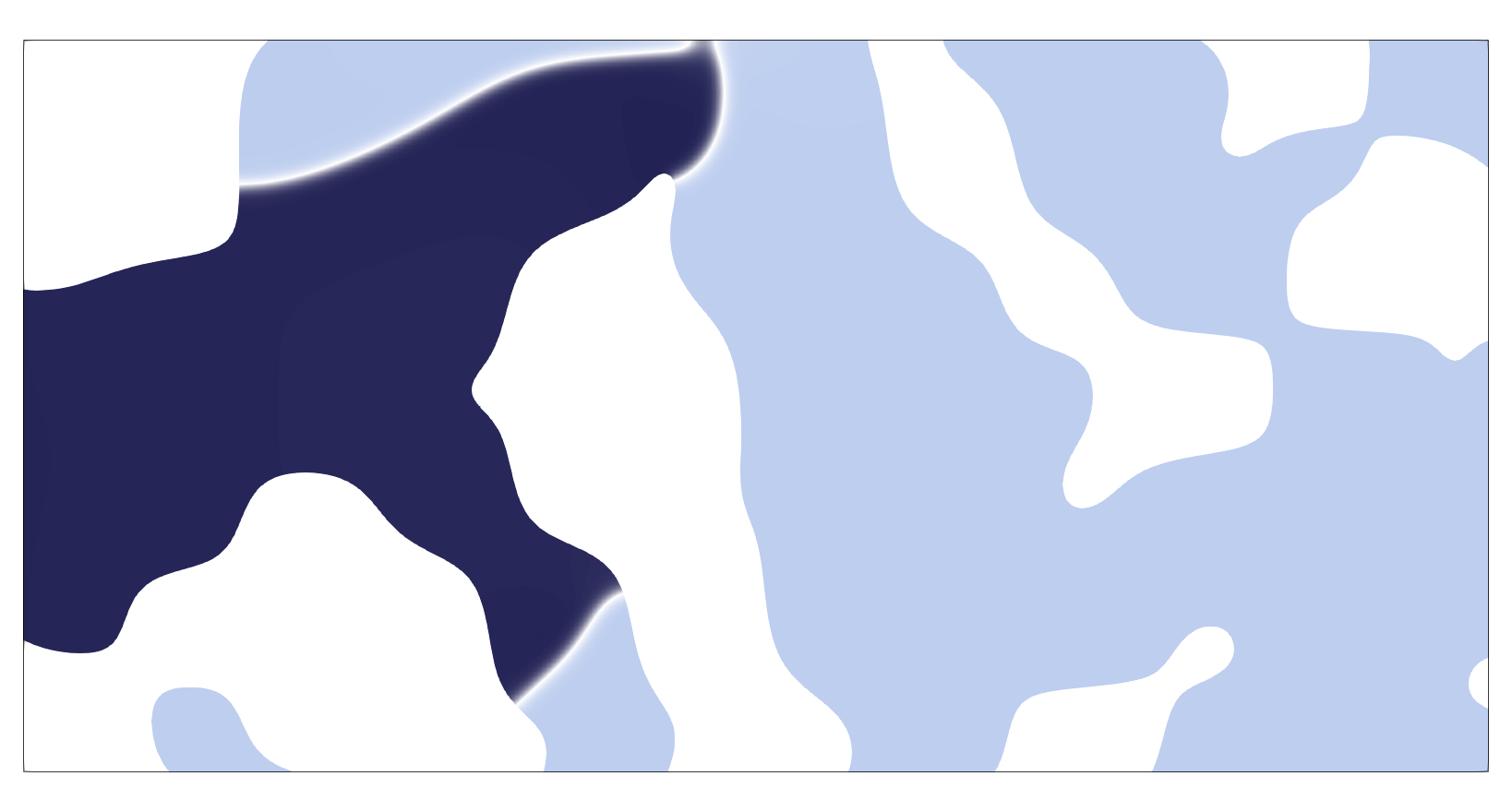}\label{fig:poroussnapshotsi}}
&
\subfloat[$t=25.9\mu \rm s$]{\includegraphics[width=0.49\textwidth]{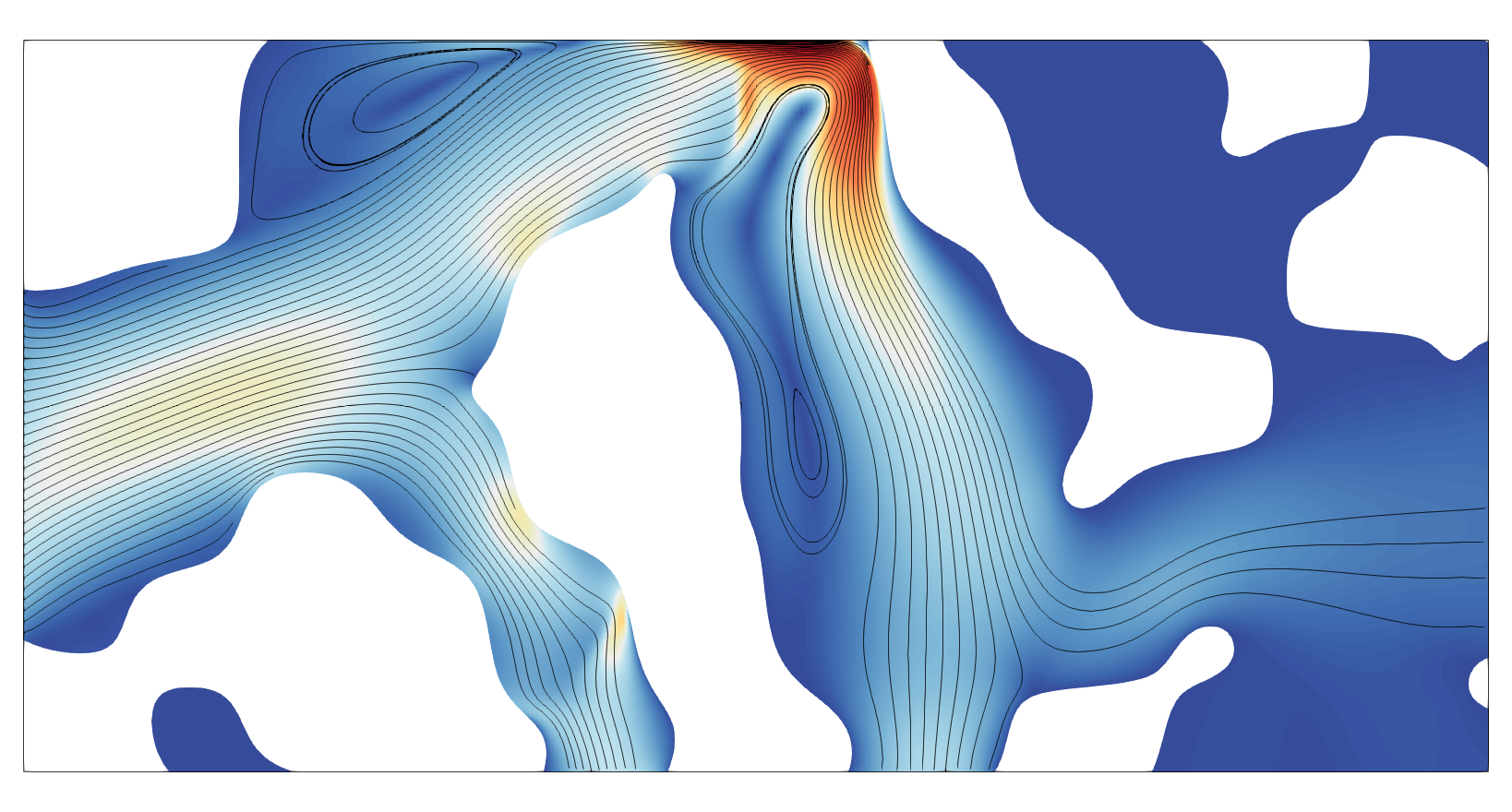}}
\\
\subfloat[$t=27.15\mu \rm s$]{\includegraphics[width=0.49\textwidth]{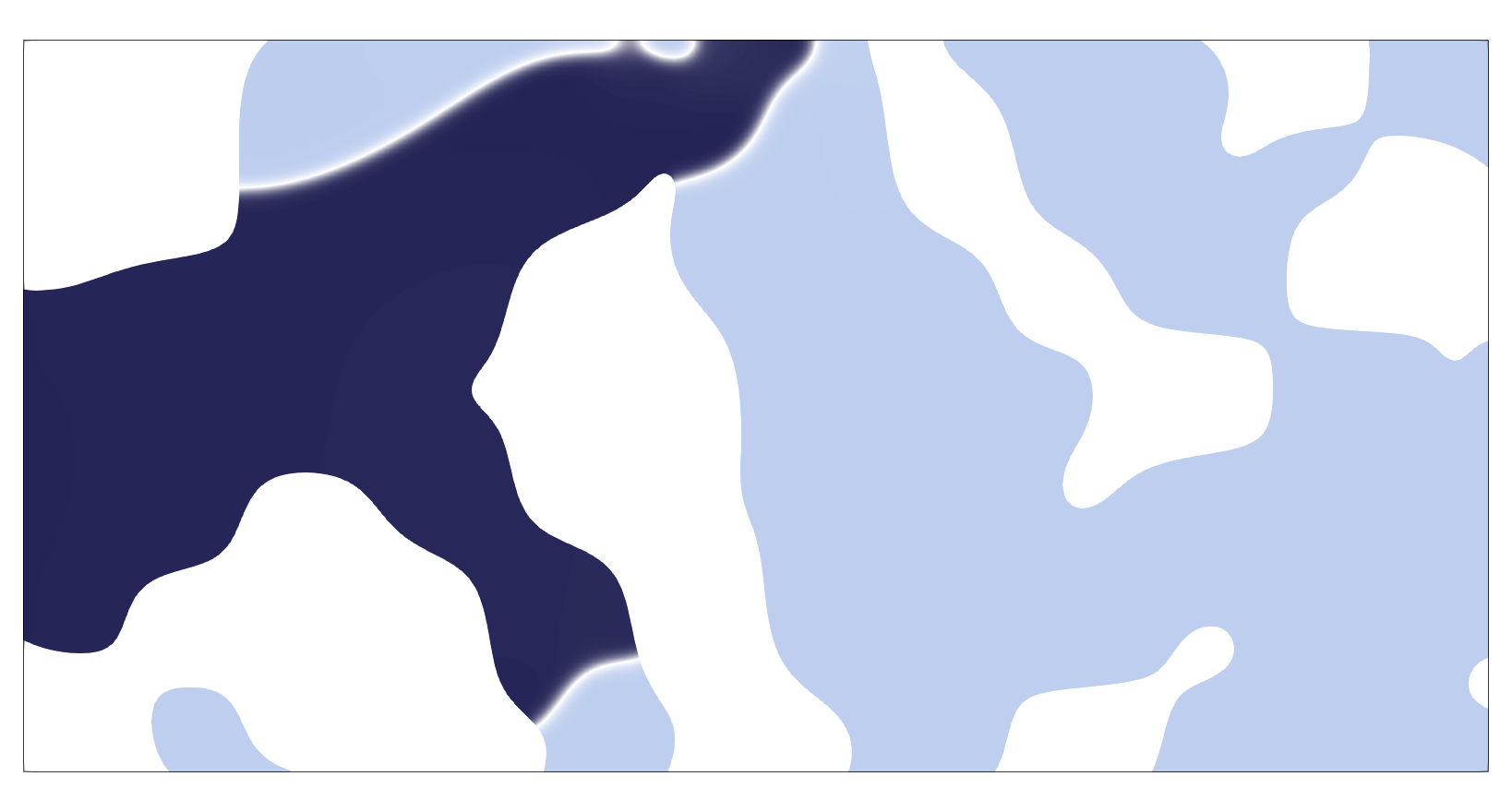}\label{fig:poroussnapshotsk}}
&
\subfloat[$t=27.15\mu \rm s$]{\includegraphics[width=0.49\textwidth]{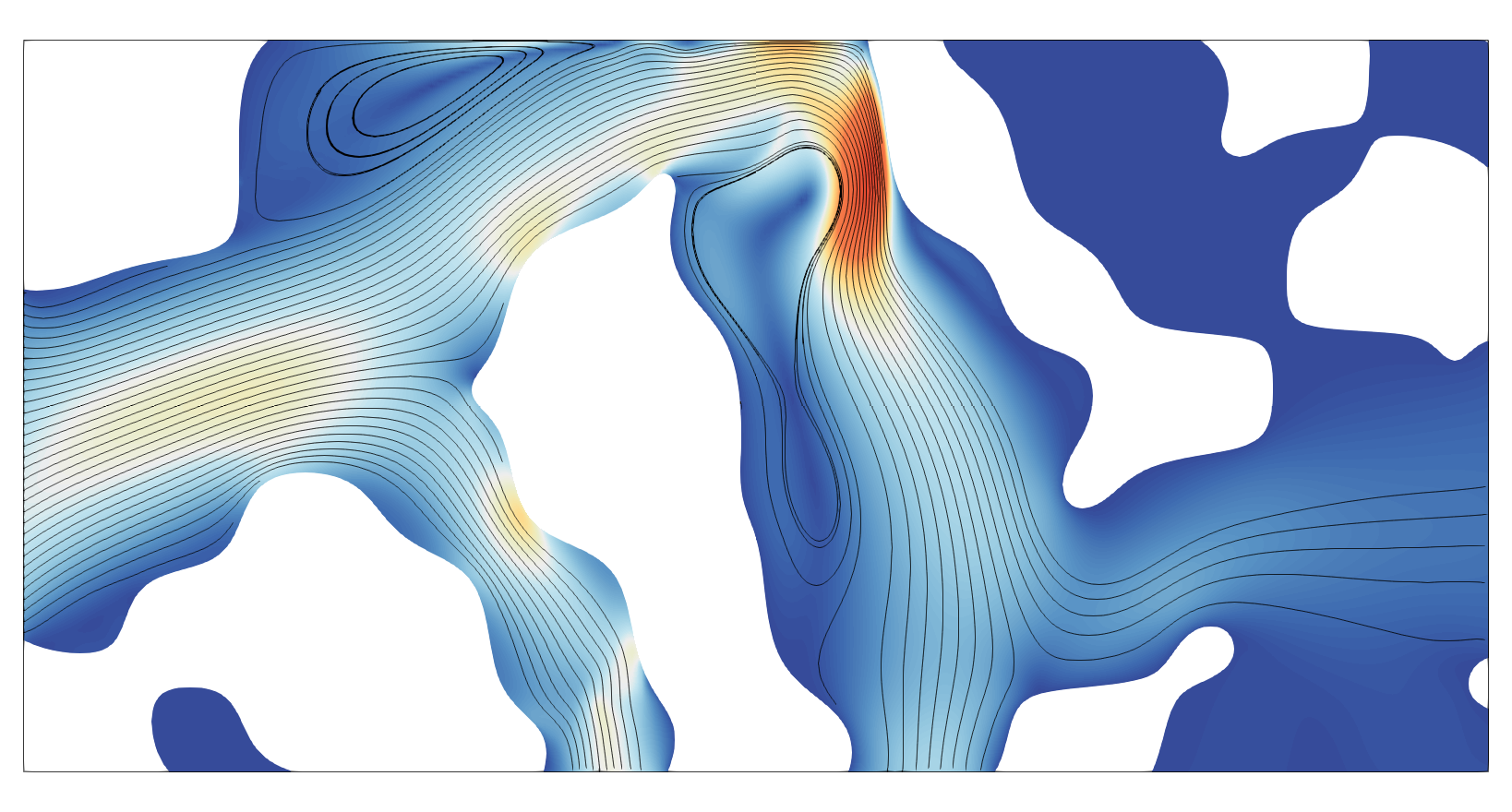}}
\\
\end{tabular}
    \caption{Time evolution of the phase field $\varphi$ (left) and velocity field $\u$ (right) for the porous medium test case.}
\end{figure}

\begin{figure}\ContinuedFloat
\begin{tabular}{cc}
\centering
\subfloat[$t=31.5\mu \rm s$]{\includegraphics[width=0.49\textwidth]{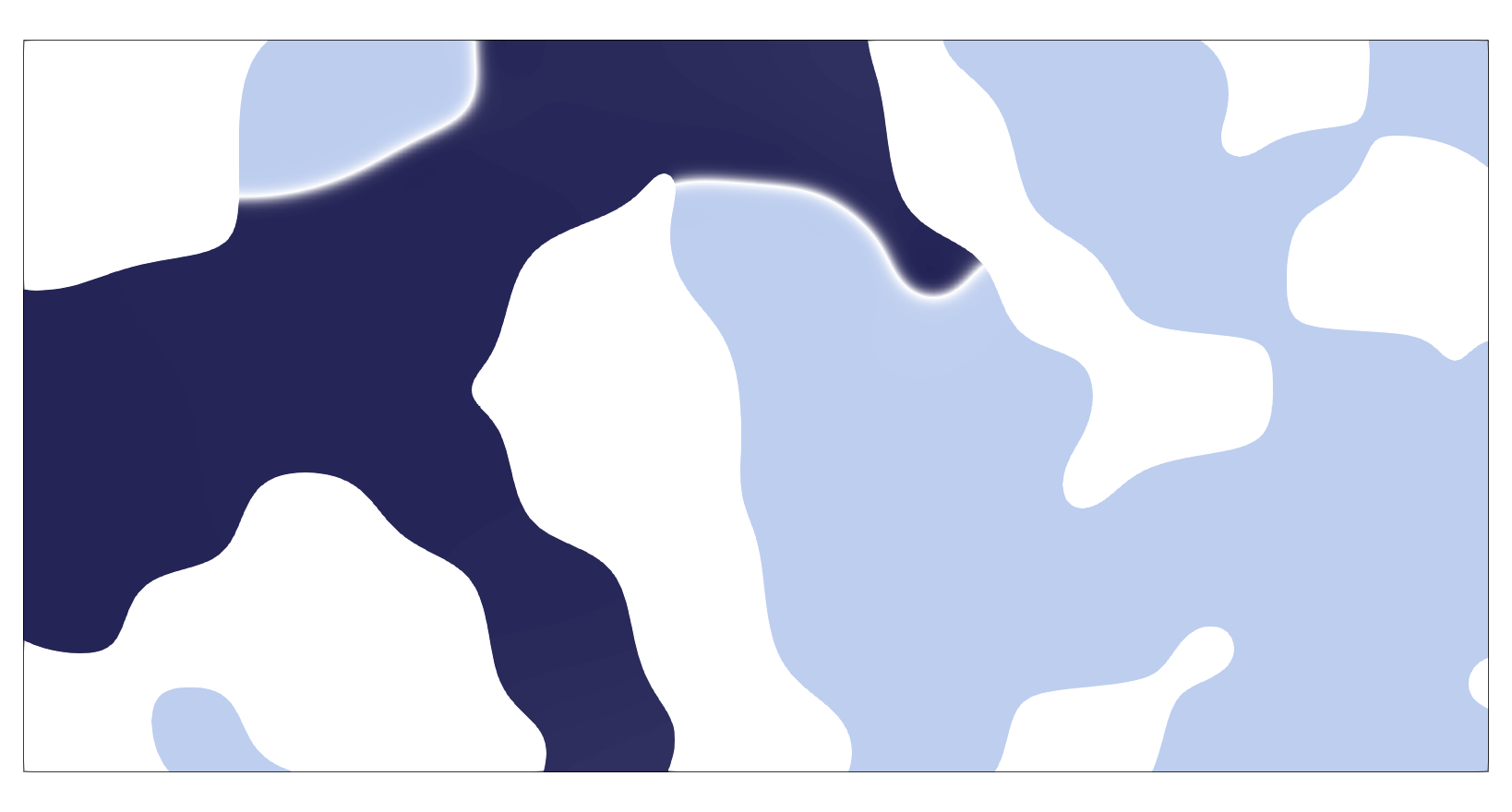}\label{fig:poroussnapshotsm}}
&
\subfloat[$t=31.5\mu \rm s$]{\includegraphics[width=0.49\textwidth]{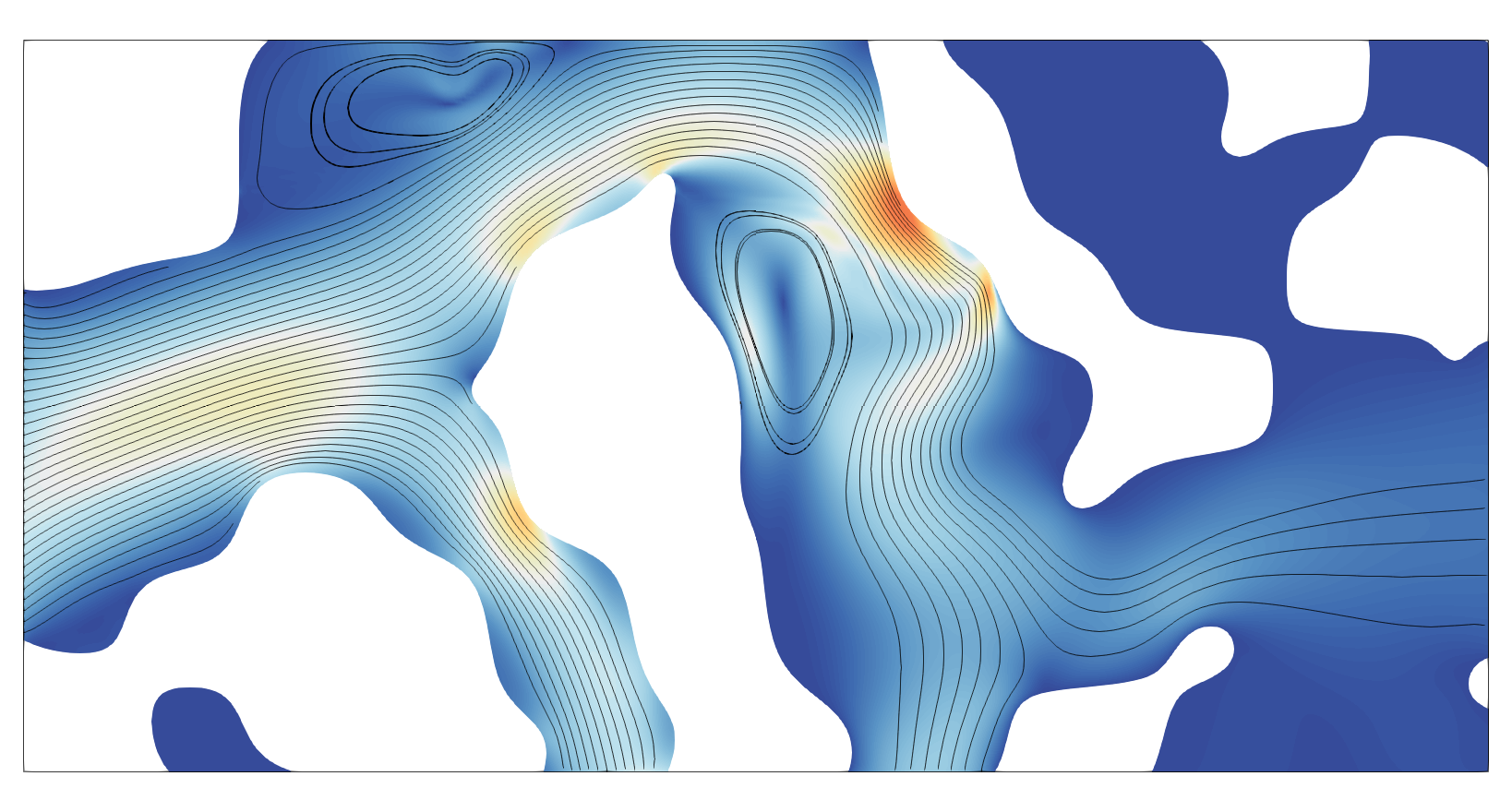}}
\\
\subfloat[$t=40.0\mu \rm s$]{\includegraphics[width=0.49\textwidth]{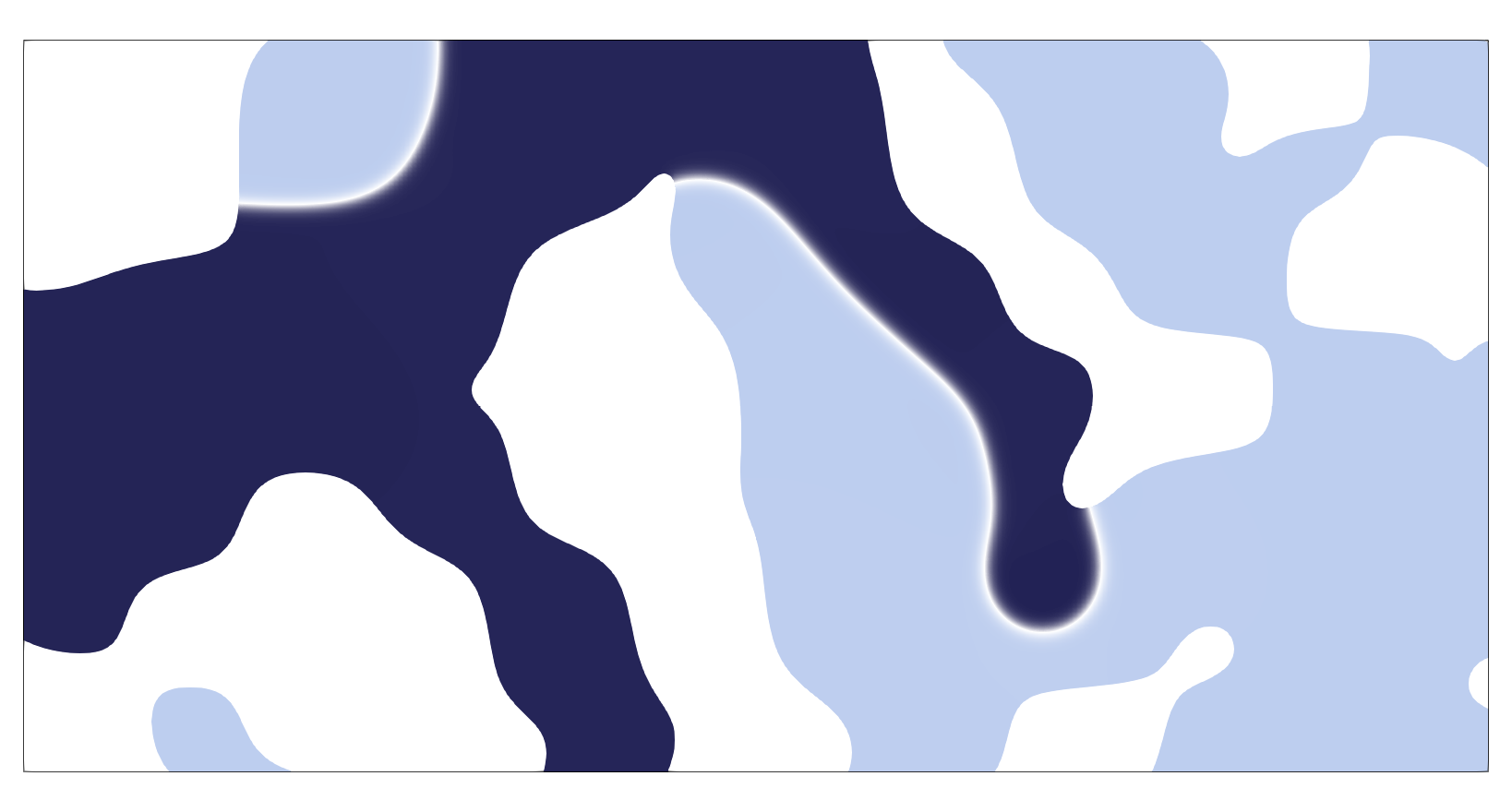}\label{fig:poroussnapshotso}}
&
\subfloat[$t=40.0\mu \rm s$]{\includegraphics[width=0.49\textwidth]{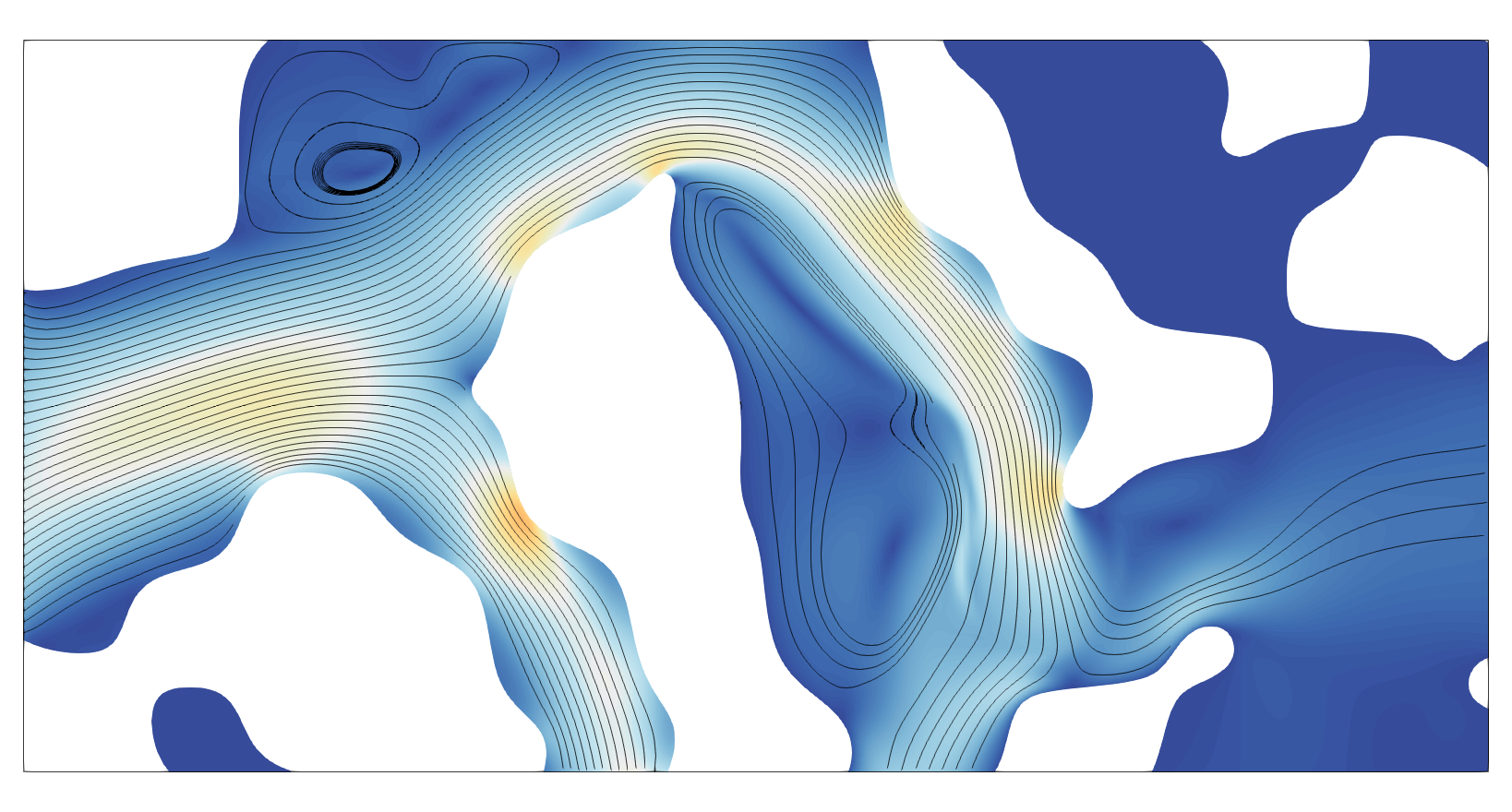}}
\\
\subfloat[$t=43.9\mu \rm s$]{\includegraphics[width=0.49\textwidth]{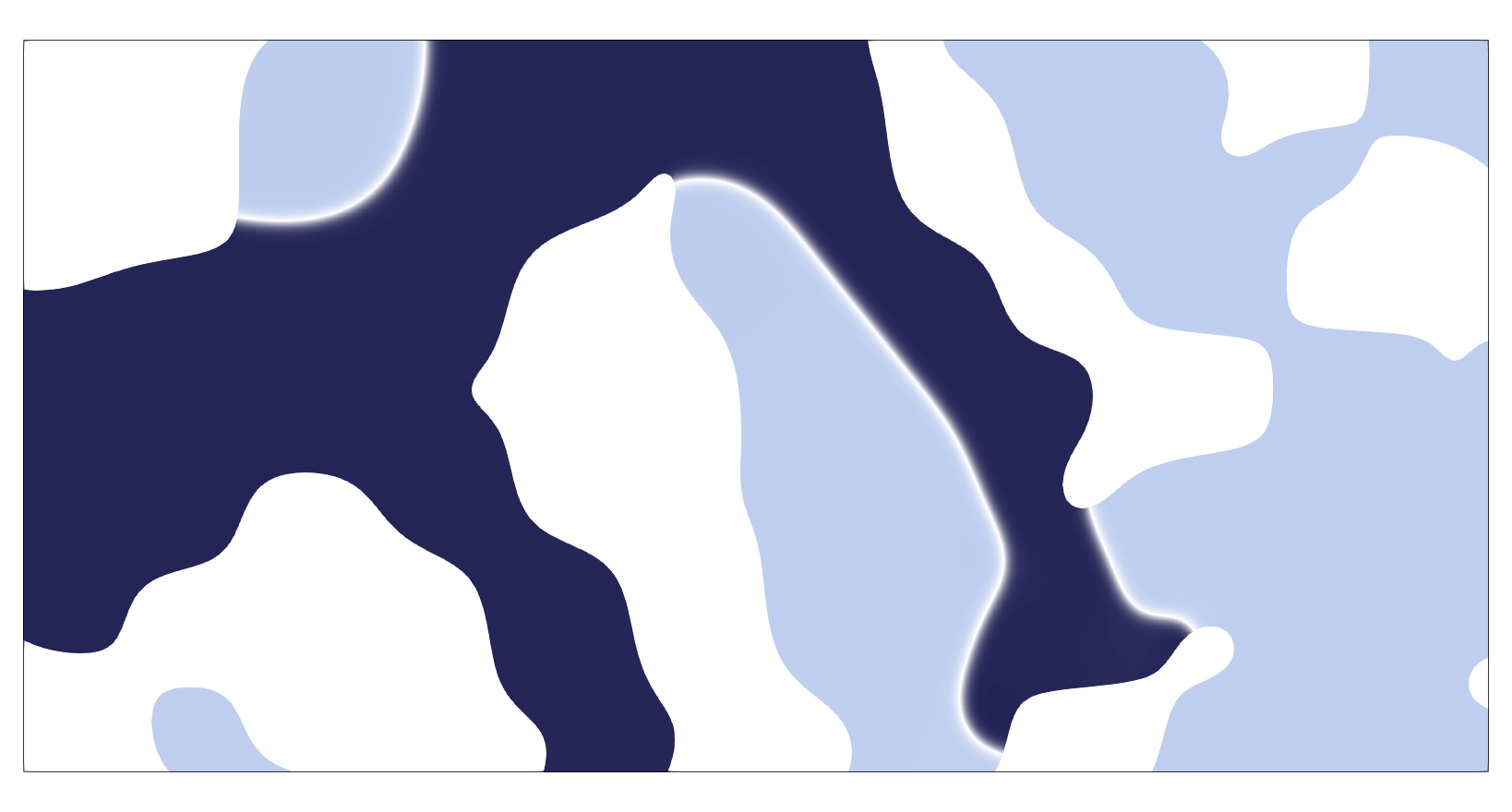}\label{fig:poroussnapshotsq}}
&
\subfloat[$t=43.9\mu \rm s$]{\includegraphics[width=0.49\textwidth]{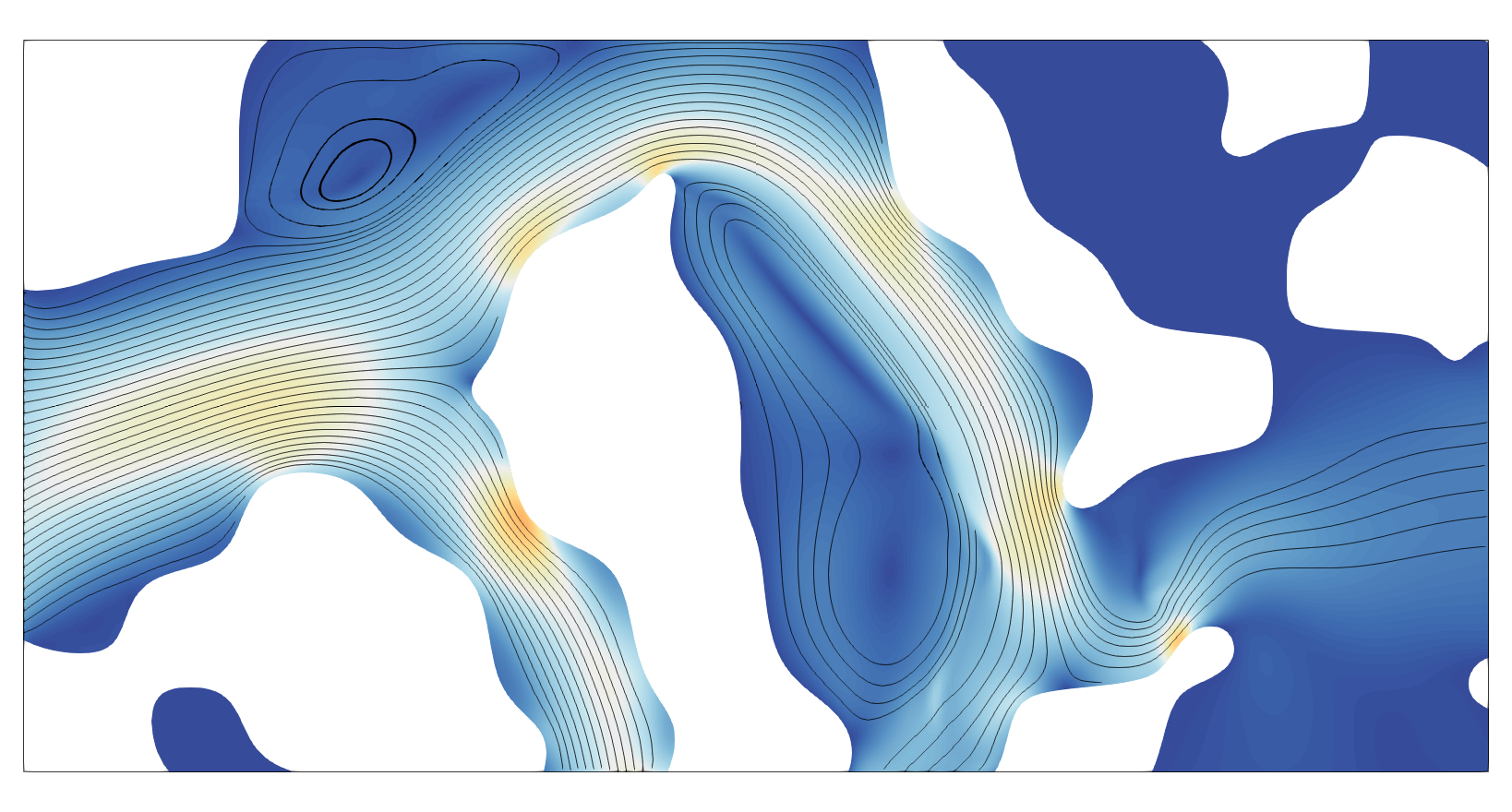}}
\\
\end{tabular}
    \caption{Time evolution of the phase field $\varphi$ (left) and velocity field $\u$ (right) for the porous medium test case.}
    \label{fig:poroussnapshots}
\end{figure}
Fig.~\ref{fig:poroussnapshotsa} shows the initial condition, where the diffuse interface is positioned at the inflow boundary. At $t=9{\rm \mu s}$ the water-front, $\varphi=1$, has moved inside the domain and meets the strong curvature of the bottom boundary of the solid domain. This stagnates the motion of the interface along that boundary. At $t=10{\rm \mu s}$, in \cref{fig:poroussnapshotse}, the inflow overcomes this stagnation, resulting in high velocities and rapid change of the interface shape. In \cref{fig:poroussnapshotsg}, the inflow phase branches between the leftmost outflow channel and the top channel, following the motion of the fluid. At $t=25.9{\rm \mu s}$ two trapped air bubbles (corresponding to $\varphi=-1$, shown in light blue) emerge at the top symmetry boundary. The small air bubble is quickly re-absorbed by the larger bubble, as illustrated in the progression snapshots of \cref{fig:poroussnapshotsi,fig:poroussnapshotsk,fig:poroussnapshotsm}. Since the velocity in this region is very small, this isolated region becomes stable in time. The high radius-of-curvature of the solid wall at the top effectively pins the water-front from \cref{fig:poroussnapshotsi} until \cref{fig:poroussnapshotsq}, where the second outflow is reached, creating a second entrapment of air in the porous domain.

\section{Conclusion}\label{sec:conclusion}

We have developed an immersed isogeometric analysis (IGA) framework to simulate binary-fluid flows on complex domains, such as encountered in, \emph{e.g.}, imbibition processes in porous media. A stabilized Galerkin formulation is proposed to robustly discretize the Navier-Stokes-Cahn-Hilliard diffuse-interface model, which describes the behavior of the binary-fluid flow by means of a velocity field, a pressure field, a phase field, and a chemical potential field. In this formulation, only the normal component of the velocity field is treated as an essential boundary condition, which is weakly enforced by Nitsche's method. For the tangential components of the velocity field, the use of a generalized Navier boundary condition is proposed. This model addresses the contact-line pinning problem related to the no-slip boundary condition and can be treated as a natural boundary condition in the immersed IGA formulation. All the boundary conditions for the phase field and chemical potential field are also natural conditions and do not require special treatment. 

The proposed framework uses optimal-regularity B-splines of the same order for all fields, resulting in a higher-order discretization with relatively few degrees of freedom compared to its Lagrange finite element counterpart. The construction of this basis is straightforward in the immersed setting, as the ambient mesh in which the computational domain is embedded is rectilinear. To obtain stable results when using this equal-order spline basis, the Galerkin formulation is amended with two forms of stabilization. First, to ensure inf-sup stability of the velocity-pressure pair, use is made of skeleton stabilization. This stabilization assigns a penalty to jumps in the higher-order non-vanishing normal derivative of the pressure field across the faces of the background mesh. Furthermore, to ameliorate problems associated with small or unfavorably cut elements, the remaining fields are stabilized through a ghost-penalty term. This ghost penalty takes the same form as the skeleton-stabilization penalty, but only needs to be applied on faces in the vicinity of the immersed boundary. The use of optimal-regularity B-splines reduced the number of stabilization parameters to two, \emph{i.e.}, one skeleton-stabilization parameter and one ghost-penalty parameter.

We have demonstrated the developed immersed isogeometric analysis framework using a range of test cases. The first test case focuses on benchmarking, by considering the well-understood binary-fluid Taylor-Couette flow between two parallel plates. The computational domain for this problem is rectangular and can hence straightforwardly be meshed, facilitating the comparison with a traditional mixed finite element formulation. The immersed framework is tested by rotation of the background mesh, which gives insights on the robustness of the method with respect to cut-element configurations. From the results we conclude that the immersed IGA framework reproduces the finite element benchmark results well, regardless of the encountered cut-element configurations, and despite its modest number of degrees of freedom. The results also show that the immersed framework has an increased error in the divergence of the velocity field, albeit that this error reduces under mesh refinement. For the considered cases, the divergence errors are not detrimental to the simulations.

Besides the Taylor-Couette flow, two porous medium test cases with complex immersed geometries are considered. The results for these cases show that the water-air flow behavior can be captured well by the immersed IGA framework. In particular, the motion of the diffuse interface along the solid boundary is observed to be influenced in a physically sound way by the geometry of the boundary. As the considered simulations pertain to uniform meshes, the overall number of degrees of freedom is strongly influenced by the interface thickness. This implies that the considered simulations, although two-dimensional, are computationally intensive. Extension to three-dimensional cases warrants the use of mesh-adaptivity with local refinements \cite{Demont:2022dk,divi_error-estimate-based_2020} and requires implementation in a high-performance computing environment.

The results presented in this work focus on demonstrating the binary-fluid flow modeling capabilities of the immersed IGA framework. In future work, the mesh dependency study of the current work should be detailed, including a study of the influence of the stability parameters on the accuracy. As future work, the ideas behind the developed stabilized formulation can be carried over to a broader class of multi-physics and phase-field problems.

\section*{Acknowledgement}
\noindent
S.K.F.\ Stoter, T.B.\ van Sluijs and T.H.B.\ Demont gratefully acknowledge the financial support through the Industrial Partnership Program {\it Fundamental Fluid Dynamics Challenges in Inkjet Printing\/} ({\it FIP\/}), a joint research program of Canon Production Printing, Eindhoven University of Technology, University of Twente, and the Netherlands Organization for Scientific Research (NWO). All simulations have been performed using the open source software package Nutils \cite{nutils}.


\begin{thebibliography}{55}
\providecommand{\natexlab}[1]{#1}
\providecommand{\url}[1]{\texttt{#1}}
\expandafter\ifx\csname urlstyle\endcsname\relax
  \providecommand{\doi}[1]{doi: #1}\else
  \providecommand{\doi}{doi: \begingroup \urlstyle{rm}\Url}\fi

\bibitem[Prosperetti and Tryggvason(2009)]{prosperetti2009computational}
A.~Prosperetti and G.~Tryggvason.
\newblock \emph{Computational methods for multiphase flow}.
\newblock Cambridge university press, 2009.

\bibitem[Sun and Beckermann(2007)]{Sun2007}
Y.~Sun and C.~Beckermann.
\newblock Sharp interface tracking using the phase-field equation.
\newblock \emph{Journal of Computational Physics}, 220, 2007.

\bibitem[Hohenberg and Halperin(1977)]{Hohenberg:1977hh}
P.C. Hohenberg and B.I. Halperin.
\newblock Theory of dynamic critical phenomena.
\newblock \emph{Reviews of Modern Physics}, 49:\penalty0 435--479, 1977.

\bibitem[Lowengrub and Truskinovsky(1998)]{Lowengrub:1998uq}
J.~Lowengrub and L.~Truskinovsky.
\newblock Quasi-incompressible {C}ahn-{H}illiard fluids and topological
  transitions.
\newblock \emph{Proceedings of the Royal Society of London A: Mathematical,
  Physical and Engineering Sciences}, 454:\penalty0 2617--2654, 1998.

\bibitem[Shokrpour~Roudbari et~al.(2018)Shokrpour~Roudbari, \c{S}im\c{s}ek,
  {van Brummelen}, and {van der Zee}]{Simsek:2018gb}
M.~Shokrpour~Roudbari, G.~\c{S}im\c{s}ek, E.H. {van Brummelen}, and K.G. {van
  der Zee}.
\newblock Diffuse-interface two-phase flow models with different densities: A
  new quasi-incompressible form and a linear energy-stable method.
\newblock \emph{Mathematical Models and Methods in Applied SciencesReviews of
  Modern Physics}, 28:\penalty0 733--770, 2018.

\bibitem[Abels et~al.(2012)Abels, Garcke, and Gr{\"u}n]{Abels:2012vn}
H.~Abels, H.~Garcke, and G.~Gr{\"u}n.
\newblock Thermodynamically consistent, frame indifferent diffuse interface
  models for incompressible two-phase flows with different densities.
\newblock \emph{Mathematical Models and Methods in Applied SciencesReviews of
  Modern Physics}, 22:\penalty0 1150013, 2012.

\bibitem[Navier(1823)]{Navier1823}
C.L.M.H. Navier.
\newblock M{\'e}moire sur les lois du mouvement des fluides.
\newblock \emph{M{\'e}moires de l’Acad{\'e}mie Royale des Sciences de
  l’Institut de France}, 6\penalty0 (1823):\penalty0 389--440, 1823.

\bibitem[Gerbeau and Leli{\`e}vre(2009)]{Gerbeau:2009jx}
J.-F. Gerbeau and T.~Leli{\`e}vre.
\newblock Generalized {N}avier boundary condition and geometric conservation
  law for surface tension.
\newblock \emph{Computer Methods in Applied Mechanics and Engineering},
  198\penalty0 (5):\penalty0 644--656, 2009.

\bibitem[Cox(1986)]{Cox:1986jq}
R.G. Cox.
\newblock The dynamics of the spreading of liquids on a solid surface. {P}art
  1. {V}iscous flow.
\newblock \emph{Journal of Fluid Mechanics}, 168:\penalty0 169--194, 1986.

\bibitem[Huh and Scriven(1971)]{huh1971hydrodynamic}
C.~Huh and L.~E. Scriven.
\newblock Hydrodynamic model of steady movement of a solid-liquid-fluid contact
  line.
\newblock \emph{Journal of Colloid and Interface Science}, 35:\penalty0
  85--101, 1971.

\bibitem[Jacqmin(2000)]{Jacqmin:2000kx}
D.~Jacqmin.
\newblock Contact-line dynamics of a diffuse fluid interface.
\newblock \emph{Journal of Fluid Mechanics}, 402:\penalty0 57--88, 2000.

\bibitem[Abels and Garcke(2018)]{Abels:2018ly}
H.~Abels and H.~Garcke.
\newblock \emph{Weak solutions and diffuse interface models for incompressible
  two-phase flows}, pages 1267--1327.
\newblock Springer International Publishing, Cham, 2018.

\bibitem[Yue et~al.(2010)Yue, Zhou, and Feng]{Yue:2010hq}
P.~Yue, C.~Zhou, and J.J. Feng.
\newblock Sharp-interface limit of the {C}ahn--{H}illiard model for moving
  contact lines.
\newblock \emph{Journal of Fluid Mechanics}, 645:\penalty0 279--294, 2010.

\bibitem[van Brummelen et~al.(2021)van Brummelen, Demont, and van
  Zwieten]{Brummelen:2021aw}
E.H. van Brummelen, T.H.B. Demont, and G.J. van Zwieten.
\newblock An adaptive isogeometric analysis approach to elasto-capillary
  fluid-solid interaction.
\newblock \emph{Int. J. Numer. Meth. Engng.}, 122\penalty0 (19):\penalty0
  5331--5352, 2021.

\bibitem[Hughes et~al.(2005)Hughes, Cottrell, and Bazilevs]{Hughes:2005it}
{T.J.R.} Hughes, {J.A.} Cottrell, and {Y.} Bazilevs.
\newblock Isogeometric analysis: {CAD}, finite elements, {NURBS}, exact
  geometry and mesh refinement.
\newblock \emph{Computer Methods in Applied Mechanics and Engineering},
  194:\penalty0 4135--4195, 2005.

\bibitem[Parvizian et~al.(2007)Parvizian, Düster, and
  Rank]{parvizian_finite_2007}
J.~Parvizian, A.~Düster, and E.~Rank.
\newblock Finite cell method: h- and p-extension for embedded domain problems
  in solid mechanics.
\newblock \emph{Computational Mechanics}, 41\penalty0 (1):\penalty0 121--133,
  2007.

\bibitem[Düster et~al.(2008)Düster, Parvizian, Yang, and
  Rank]{duster_finite_2008}
A.~Düster, J.~Parvizian, Z.~Yang, and E.~Rank.
\newblock The finite cell method for three-dimensional problems of solid
  mechanics.
\newblock \emph{Computer Methods in Applied Mechanics and Engineering},
  197\penalty0 (45):\penalty0 3768--3782, 2008.

\bibitem[Schillinger and Ruess(2015)]{schillinger_finite_2015}
D.~Schillinger and M.~Ruess.
\newblock The finite cell method: A review in the context of higher-order
  structural analysis of {CAD} and image-based geometric models.
\newblock \emph{Archives of Computational Methods in Engineering}, 22\penalty0
  (3):\penalty0 391--455, 2015.

\bibitem[Burman(2010)]{burman_ghost_2010}
E.~Burman.
\newblock Ghost penalty.
\newblock \emph{Comptes Rendus Mathematique}, 348\penalty0 (21):\penalty0
  1217--1220, 2010.

\bibitem[Burman and Hansbo(2012)]{burman_fictitious_2012}
E.~Burman and P.~Hansbo.
\newblock Fictitious domain finite element methods using cut elements: {II}.
  {A} stabilized {N}itsche method.
\newblock \emph{Applied Numerical Mathematics}, 62\penalty0 (4):\penalty0
  328--341, 2012.

\bibitem[Burman et~al.(2015{\natexlab{a}})Burman, Claus, Hansbo, Larson, and
  Massing]{burman_cutfem_2015}
E.~Burman, S.~Claus, P.~Hansbo, M.G. Larson, and A.~Massing.
\newblock {CutFEM}: Discretizing geometry and partial differential equations.
\newblock \emph{International Journal for Numerical Methods in Engineering},
  104\penalty0 (7):\penalty0 472--501, 2015{\natexlab{a}}.

\bibitem[Rank et~al.(2012)Rank, Ruess, Kollmannsberger, Schillinger, and
  Düster]{rank_geometric_2012}
E.~Rank, M.~Ruess, S.~Kollmannsberger, D.~Schillinger, and A.~Düster.
\newblock Geometric modeling, isogeometric analysis and the finite cell method.
\newblock \emph{Computer Methods in Applied Mechanics and Engineering},
  249-252:\penalty0 104--115, 2012.

\bibitem[Schillinger et~al.(2012)Schillinger, Dedè, Scott, Evans, Borden,
  Rank, and Hughes]{schillinger_isogeometric_2012}
D.~Schillinger, L.~Dedè, M.A. Scott, J.A. Evans, M.J. Borden, E.~Rank, and
  T.J.R. Hughes.
\newblock An isogeometric design-through-analysis methodology based on adaptive
  hierarchical refinement of {NURBS}, immersed boundary methods, and t-spline
  {CAD} surfaces.
\newblock \emph{Computer Methods in Applied Mechanics and Engineering},
  249-252:\penalty0 116--150, 2012.

\bibitem[Ruess et~al.(2013)Ruess, Schillinger, Bazilevs, Varduhn, and
  Rank]{ruess_weakly_2013}
M.~Ruess, D.~Schillinger, Y.~Bazilevs, V.~Varduhn, and E.~Rank.
\newblock Weakly enforced essential boundary conditions for {NURBS}-embedded
  and trimmed {NURBS} geometries on the basis of the finite cell method.
\newblock \emph{International Journal for Numerical Methods in Engineering},
  95\penalty0 (10):\penalty0 811--846, 2013.

\bibitem[Kamensky et~al.(2015)Kamensky, Hsu, Schillinger, Evans, Aggarwal,
  Bazilevs, Sacks, and Hughes]{kamensky_immersogeometric_2015}
D.~Kamensky, M.-C. Hsu, D.~Schillinger, J.A. Evans, A.~Aggarwal, Y.~Bazilevs,
  M.S. Sacks, and T.J.R. Hughes.
\newblock An immersogeometric variational framework for fluid–structure
  interaction: Application to bioprosthetic heart valves.
\newblock \emph{Computer Methods in Applied Mechanics and Engineering},
  284:\penalty0 1005--1053, 2015.

\bibitem[Hsu et~al.(2016)Hsu, Wang, Xu, Herrema, and
  Krishnamurthy]{hsu_direct_2016}
M.-C. Hsu, C.~Wang, F.~Xu, A.J. Herrema, and A.~Krishnamurthy.
\newblock Direct immersogeometric fluid flow analysis using b-rep {CAD} models.
\newblock \emph{Computer Aided Geometric Design}, 43:\penalty0 143--158, 2016.

\bibitem[Verhoosel et~al.(2015)Verhoosel, van Zwieten, van Rietbergen, and
  de~Borst]{verhoosel_image-based_2015}
C.V. Verhoosel, G.~van Zwieten, B.~van Rietbergen, and R.~de~Borst.
\newblock Image-based goal-oriented adaptive isogeometric analysis with
  application to the micro-mechanical modeling of trabecular bone.
\newblock \emph{Computer Methods in Applied Mechanics and Engineering},
  284:\penalty0 138--164, 2015.

\bibitem[Ruess et~al.(2012)Ruess, Tal, Trabelsi, Yosibash, and
  Rank]{ruess_finite_2012}
M.~Ruess, D.~Tal, N.~Trabelsi, Z.~Yosibash, and E.~Rank.
\newblock The finite cell method for bone simulations: verification and
  validation.
\newblock \emph{Biomechanics and Modeling in Mechanobiology}, 11\penalty0
  (3):\penalty0 425--437, 2012.

\bibitem[de~Prenter et~al.(2020)de~Prenter, Verhoosel, van Brummelen, Evans,
  Messe, Benzaken, and Maute]{de_prenter_multigrid_2020}
F.~de~Prenter, C.V. Verhoosel, E.H. van Brummelen, J.A. Evans, C.~Messe,
  J.~Benzaken, and K.~Maute.
\newblock Multigrid solvers for immersed finite element methods and immersed
  isogeometric analysis.
\newblock \emph{Computational Mechanics}, 65\penalty0 (3):\penalty0 807--838,
  2020.

\bibitem[Abedian et~al.(2013)Abedian, Parvizian, Düster, Khademyzadeh, and
  Rank]{abedian_performance_2013}
A.~Abedian, J.~Parvizian, A.~Düster, H.~Khademyzadeh, and E.~Rank.
\newblock Performance of different integration schemes in facing
  discontinuities in the finite cell method.
\newblock \emph{International Journal of Computational Methods}, 10\penalty0
  (3):\penalty0 1350002, 2013.

\bibitem[Divi et~al.(2020)Divi, Verhoosel, Auricchio, Reali, and van
  Brummelen]{divi_error-estimate-based_2020}
S.C. Divi, C.V. Verhoosel, F.~Auricchio, A.~Reali, and E.H. van Brummelen.
\newblock Error-estimate-based adaptive integration for immersed isogeometric
  analysis.
\newblock \emph{Computers \& Mathematics with Applications}, 80\penalty0
  (11):\penalty0 2481--2516, 2020.

\bibitem[Douglas and Dupont(1975)]{Douglas1975}
J.~Douglas and T.~Dupont.
\newblock {A Galerkin method for a nonlinear Dirichlet problem}.
\newblock \emph{Mathematics of Computation}, 29\penalty0 (131), 1975.

\bibitem[Burman et~al.(2015{\natexlab{b}})Burman, Claus, and
  Massing]{Burman2015}
E.~Burman, S.~Claus, and A.~Massing.
\newblock {A stabilized cut finite element method for the three field {S}tokes
  problem}.
\newblock \emph{SIAM Journal on Scientific Computing}, 37\penalty0 (4),
  2015{\natexlab{b}}.

\bibitem[de~Prenter et~al.(2023)de~Prenter, Verhoosel, van Brummelen, Larson,
  and Badia]{deprenter2023stability}
F.~de~Prenter, C.V. Verhoosel, E.H. van Brummelen, M.G. Larson, and S.~Badia.
\newblock Stability and conditioning of immersed finite element methods:
  analysis and remedies.
\newblock \emph{Archives of Computational Methods in Engineering}, pages 1--40,
  2023.

\bibitem[Stoter et~al.(2023)Stoter, Divi, van Brummelen, Larson, de~Prenter,
  and Verhoosel]{Stoter2023a}
S.K.F. Stoter, S.C. Divi, E.H. van Brummelen, M.G. Larson, F.~de~Prenter, and
  C.V. Verhoosel.
\newblock {Critical time-step size analysis and mass scaling by ghost-penalty
  for immersogeometric explicit dynamics}.
\newblock \emph{Computer Methods in Applied Mechanics and Engineering},
  412:\penalty0 116074, 2023.

\bibitem[Divi et~al.(2022{\natexlab{a}})Divi, Verhoosel, Auricchio, Reali, and
  van Brummelen]{divi_topology-preserving_2022}
S.C. Divi, C.V. Verhoosel, F.~Auricchio, A.~Reali, and E.H. van Brummelen.
\newblock Topology-preserving scan-based immersed isogeometric analysis.
\newblock \emph{Computer Methods in Applied Mechanics and Engineering},
  392:\penalty0 114648, 2022{\natexlab{a}}.

\bibitem[Karatzas and Rozza(2021)]{karatzas2021reduced}
E.N. Karatzas and G.~Rozza.
\newblock A reduced order model for a stable embedded boundary parametrized
  {C}ahn--{H}illiard phase-field system based on cut finite elements.
\newblock \emph{Journal of Scientific Computing}, 89\penalty0 (1):\penalty0 9,
  2021.

\bibitem[Hoang et~al.(2018)Hoang, Verhoosel, Auricchio, van Brummelen, and
  Reali]{hoang_skeleton-stabilized_2018}
T.~Hoang, C.V. Verhoosel, F.~Auricchio, E.H. van Brummelen, and A.~Reali.
\newblock Skeleton-stabilized {IsoGeometric} analysis: High-regularity
  interior-penalty methods for incompressible viscous flow problems.
\newblock \emph{Computer Methods in Applied Mechanics and Engineering},
  337:\penalty0 324--351, 2018.

\bibitem[Hoang et~al.(2019)Hoang, Verhoosel, Qin, Auricchio, Reali, and van
  Brummelen]{hoang_skeleton-stabilized_2019}
T.~Hoang, C.V. Verhoosel, C.-Z. Qin, F.~Auricchio, A.~Reali, and E.H. van
  Brummelen.
\newblock Skeleton-stabilized immersogeometric analysis for incompressible
  viscous flow problems.
\newblock \emph{Computer Methods in Applied Mechanics and Engineering},
  344:\penalty0 421--450, 2019.

\bibitem[Demont et~al.(2022)Demont, van Zwieten, Diddens, and van
  Brummelen]{Demont:2022dk}
T.H.B. Demont, G.J. van Zwieten, C.~Diddens, and E.H. van Brummelen.
\newblock A robust and accurate adaptive approximation method for a
  diffuse-interface model of binary-fluid flows.
\newblock \emph{Computer Methods}, 400:\penalty0 115563, 2022.

\bibitem[Bonart et~al.(2019)Bonart, Kahle, and Repke]{Bonart:2019re}
H.~Bonart, C.~Kahle, and J.-U. Repke.
\newblock Comparison of energy stable simulation of moving contact line
  problems using a thermodynamically consistent {C}ahn--{H}illiard
  {N}avier--{S}tokes model.
\newblock \emph{Journal of Computational Physics}, 399:\penalty0 108959, 2019.

\bibitem[Arrhenius(1887)]{Arrhenius:1887xr}
S.~Arrhenius.
\newblock {\"U}ber die innere {R}eibung verd{\"u}nnter w{\"a}sseriger
  {L}{\"o}sungen.
\newblock \emph{Zeitschrift für Physikalische Chemie}, 1U\penalty0
  (1):\penalty0 285--298, 1887.

\bibitem[Yue and Feng(2011)]{Yue:2011uq}
P.~Yue and J.J. Feng.
\newblock Wall energy relaxation in the {C}ahn--{H}illiard model for moving
  contact lines.
\newblock \emph{Physics of Fluids}, 23:\penalty0 012106, 2011.

\bibitem[Shokrpour~Roudbari et~al.(2016)Shokrpour~Roudbari, {van Brummelen},
  and Verhoosel]{Shokrpour-Roudbari:2016dp}
M.~Shokrpour~Roudbari, E.H. {van Brummelen}, and C.V. Verhoosel.
\newblock A multiscale diffuse-interface model for two-phase flow in porous
  media.
\newblock \emph{Computers \& Fluids}, 141:\penalty0 212--222, 2016.

\bibitem[{van Brummelen} et~al.(2016){van Brummelen}, Shokrpour~Roudbari, and
  {van Zwieten}]{Brummelen:2016qa}
{E.H.} {van Brummelen}, M.~Shokrpour~Roudbari, and {G.J.} {van Zwieten}.
\newblock Elasto-capillarity simulations based on the
  {N}avier-{S}tokes-{C}ahn-{H}illiard equations.
\newblock In \emph{Advances in Computational Fluid-Structure Interaction and
  Flow Simulation}, Modeling and Simulation in Science, Engineering and
  Technology, pages 451--462. Birkh{{\"a}}user, 2016.

\bibitem[Cottrell et~al.(2009)Cottrell, Hughes, and Bazilevs]{Cottrell:2009ad}
J.A. Cottrell, T.J.R. Hughes, and Y.~Bazilevs.
\newblock \emph{Isogeometric Analysis: Toward Integration of {CAD} and {FEA}}.
\newblock Wiley, Chichester, 2009.

\bibitem[Giannelli et~al.(2012)Giannelli, J\"uttler, and
  Speleers]{Giannelli:2012rr}
C.~Giannelli, B.~J\"uttler, and H.~Speleers.
\newblock {THB}-splines: The truncated basis for hierarchical splines.
\newblock \emph{Computer Aided Geometric Design}, 29:\penalty0 485--498, 2012.

\bibitem[Divi et~al.(2022{\natexlab{b}})Divi, van Zuijlen, Hoang, de~Prenter,
  Auricchio, Reali, van Brummelen, and Verhoosel]{divi_residual-based_2022}
S.C. Divi, P.H. van Zuijlen, T.~Hoang, F.~de~Prenter, F.~Auricchio, A.~Reali,
  E.H. van Brummelen, and C.V. Verhoosel.
\newblock Residual-based error estimation and adaptivity for stabilized
  immersed isogeometric analysis using truncated hierarchical {B}-splines.
\newblock \emph{Journal of Mechanics}, 38:\penalty0 204--237,
  2022{\natexlab{b}}.

\bibitem[Verhoosel et~al.(2023)Verhoosel, van Brummelen, Divi, and
  de~Prenter]{verhoosel2022scan}
C.V. Verhoosel, E.H. van Brummelen, S.C. Divi, and F.~de~Prenter.
\newblock \emph{Frontiers in Computational Fluid-Structure Interaction and Flow
  Simulation: Research from Lead Investigators under 40}, chapter Scan-based
  immersed isogeometric flow analysis.
\newblock Springer Nature Switzerland AG, 2023.
\newblock accepted for publication (in press), preprint arXiv:2208.14994.

\bibitem[de~Prenter et~al.(2018)de~Prenter, Lehrenfeld, and
  Massing]{de_prenter_note_2018}
F.~de~Prenter, C.~Lehrenfeld, and A.~Massing.
\newblock A note on the stability parameter in nitsche’s method for unfitted
  boundary value problems.
\newblock \emph{Computers \& Mathematics with Applications}, 75\penalty0
  (12):\penalty0 4322--4336, 2018.

\bibitem[Badia et~al.(2018)Badia, Martin, and Verdugo]{badia_mixed_2018}
S.~Badia, A.F. Martin, and F.~Verdugo.
\newblock Mixed aggregated finite element methods for the unfitted
  discretization of the stokes problem.
\newblock \emph{{SIAM} Journal on Scientific Computing}, 40\penalty0
  (6):\penalty0 B1541--B1576, 2018.

\bibitem[Burman and Hansbo(2006)]{burman_edge_2006}
E.~Burman and P.~Hansbo.
\newblock Edge stabilization for the generalized {S}tokes problem: {A}
  continuous interior penalty method.
\newblock \emph{Computer Methods in Applied Mechanics and Engineering},
  195\penalty0 (19):\penalty0 2393--2410, 2006.

\bibitem[Badia et~al.(2022)Badia, Neiva, and Verdugo]{badia_linking_2022}
S.~Badia, E.~Neiva, and F.~Verdugo.
\newblock Linking ghost penalty and aggregated unfitted methods.
\newblock \emph{Computer Methods in Applied Mechanics and Engineering},
  388:\penalty0 114232, 2022.

\bibitem[{van Zwieten} et~al.(2022){van Zwieten}, {van Zwieten}, and
  Hoitinga]{nutils}
{G.J.} {van Zwieten}, T.M. {van Zwieten}, J., and {W.} Hoitinga.
\newblock Nutils (ver. 7.0).
\newblock https://dx.doi.org/10.5281/zenodo.6006701, 2022.

\bibitem[Jacqmin(1999)]{Jacqmin:1999fk}
D.~Jacqmin.
\newblock Calculation of two-phase {N}avier-{S}tokes flows using phase-field
  modeling.
\newblock \emph{Journal of Computational Physics}, 155:\penalty0 96--127, 10
  1999.

\end{thebibliography}
\end{document}